\documentclass[12pt,a4paper]{article}
\usepackage[T2A]{fontenc}
\usepackage[cp1251]{inputenc}
\usepackage[russian]{babel}
\usepackage{amsmath,amssymb,amsthm,amsfonts,amscd}
\begin{document}

\title{Геометрический подход к стабильным гомотопическим группам
сфер.I. Инвариант Хопфа}
\author{П.М.Ахметьев  \thanks{Работа автора поддержена  грантами INTAS 00-0259, РФФИ 08-01-00663-а.}}
\date{посвящается памяти проф. М.М.Постникова}

\sloppy \theoremstyle{plain}
\newtheorem{theorem}{Теорема}
\newtheorem*{main*}{Основная Теорема}
\newtheorem*{theorem*}{Теорема}
\newtheorem{lemma}[theorem]{Лемма}
\newtheorem{proposition}[theorem]{Предложение}
\newtheorem{corollary}[theorem]{Следствие}
\newtheorem{conjecture}[theorem]{Гипотеза}
\newtheorem{problem}[theorem]{Проблема}

\theoremstyle{definition}
\newtheorem{definition}[theorem]{Определение}
\newtheorem{remark}[theorem]{Замечание}
\newtheorem*{remark*}{Замечание}
\newtheorem*{example*}{Пример}
\newtheorem{example}[theorem]{Пример}

\def\Z{{\Bbb Z}}
\def\R{{\Bbb R}}
\def\RP{{\Bbb R}\!{\rm P}}
\def\N{{\Bbb N}}
\def\C{{\Bbb C}}
\def\A{{\bf A}}
\def\D{{\bf D}}
\def\Q{{\bf Q}}
\def\i{{\bf i}}
\def\j{{\bf j}}
\def\k{{\bf k}}
\def\E{{\bf E}}
\def\F{{\bf F}}
\def\J{{\bf J}}
\def\G{{\bf G}}
\def\I{{\bf I}}
\def\e{{\bf e}}
\def\f{{\bf f}}
\def\d{{\bf d}}
\def\H{{\bf H}}
\def\fr{{\operatorname{fr}}}
\def\st{{\operatorname{st}}}
\def\mod{{\operatorname{mod}\,}}
\def\cyl{{\operatorname{cyl}}}
\def\dist{{\operatorname{dist}}}
\def\sf{{\operatorname{sf}}}
\def\dim{{\operatorname{dim}}}
\def\dist{\operatorname{dist}}

\maketitle

\begin{abstract}
В работе развивается геометрический подход к стабильным
гомотопическим группам сфер, основанный на конструкции
Понтрягина-Тома. В рамках этого подхода получено доказательство
теоремы Адамса об инварианте Хопфа для всех размерностей, исключая
$15,31$. Доказывается, что при $n
> 31$ в стабильной гомотопической группе сфер $\Pi_n$ не
существует элемента с инвариантом Хопфа 1. Новое доказательство
основано на методах геометрической топологии. Используется теорема
Понтрягина-Тома в форме Уоллеса о представлении стабильных
гомотопических групп $\lim_{k \to +\infty} \pi_{k+n}(\Sigma^k(\RP^{\infty}))$
вещественного проективного бесконечномерного
пространства (которые по теореме Кана-Придди эпиморфно
отображаются на 2-компоненты стабильных гомотопических групп сфер)
классами кобордизма погружений (вообще говоря, неориентируемых)
$n-1$--мерных многообразий в коразмерности 1. Инвариант Хопфа выражаеся
характеристическим классом $n-1$--мерного многообразия
двукратных точек самопересечения погруженного $n$--мерного многообразия,
представляющего заданный элемент в стабильной гомотопической
группе $\Pi_n$.
\end{abstract}

\section*{Введение}

Пусть $\pi_{n+m}(S^m)$ -- гомотопические группы сфер. При условии
$m \ge n+2$ эта группа не зависит от $m$ и обозначается $\Pi_n$.
Она называется стабильной гомотопической группой сфер в
размерности $n$. Проблема вычисления стабильных гомотопических
групп является одной из основных нерешенных проблем в топологии.
Эта проблема имеет важное прикладное значение при изучении пространства вещественных функций
(см. [В], гл.3) и при изучении задачи аппроксимации отображений вложениями [Ме2].

При вычислении элементов стабильной гомотопической группы сфер
изучают те алгебраические инварианты, определение которых дается
сразу для всех размерностей (или для некоторой бесконечной
последовательности размерностей). Тем не менее, эти инварианты,
как правило, оказываются тривиальными и не вырождаются лишь в
исключительных случаях, подробнее см. [M].

Основным алгебраическим инвариантом на гомотопических группах сфер
является  стабильный инвариант Хопфа, который представляет собой гомоморфизм
$$ h: \Pi_{2k-1} \to \Z/2, $$
и который называется также инвариантом Стинрода-Хопфа, см.,
например, [M-T] по поводу определения и основных свойств.
 Стабильный инвариант
Хопфа изучается в работе.

Доказательство следующей теоремы было получено Адамсом в [A].

\begin{theorem*}{Адамс [A]}

Стабильный инвариант Хопфа $h: \Pi_n \to \Z/2$, $n \equiv 1 \pmod{2}$,
является тривиальным тогда и только тогда, когда $n
\ne 1,3,7$.
\end{theorem*}

\begin{remark*}
Теорема Адамса в простом случае $n \ne 2^k-1$ была доказана Адемом
[Ad]. Cлучай $n=15$ был доказан Тодой (см. [M-T] гл.18).
\end{remark*}

Впоследствии Адамс и Атья предложили альтернативный подход к
изучению инвариантов Хопфа, основанный на результатах из
$K$-теории и на теореме Ботта о периодичности см. [A-A]. В
последующих работах этот подход также был обобщен. Простое
доказательство Теоремы Адамса, близкое к доказательству
Адамса-Атьи, было предложено В.М.Бух\-шта\-бе\-ром в [B].

Определение стабильного инварианта Хопфа переформулируется на языке
групп кобордизмов погруженных многообразий. В основе этого метода
лежат результаты работ [E1, K2, K-S1, K-S2, La]. Стабильный инвариант Хопфа
равен характеристическому числу многообразия двукратных точек
самопересечения исходного многообразия, представляющего элемент в
исходной группе кобордизма погружений.
Это, например, явно сформулировано в работе [E1], Лемма 3.1. Формулировка указанной леммы стандартным способом
переводится на язык теории погружений при помощи конструкции Понтрягина-Тома.

 В случае $n \not\equiv 3 \pmod{4}$
доказательство теоремы Адамса, использующее теорию погружений,
было получено А.Сючем в работе [Sz]. Следующий по сложности
случай (теорема Адема) возникает при $n \ne 2^l-1$. Теорема Адамса
в этом случае была передоказана геометрическими методами в
совместной работе автора с А. Сючем [A-Sz].

Далее в работе всюду предполагается, что $n=2^l-1$. Основной
результат настоящей работы обобщает предыдущий и состоит в
следующем.

\subsubsection*{Основная Теорема}
Пусть $g: M^{\frac{3n+3}{4}+7} \looparrowright \R^n$--
произвольное гладкое погружение замкнутого, вообще говоря,
несвязного многообразия $M$, $\dim(M)=\frac{3n+3}{4}+7$, причем
нормальное расслоение $\nu(g)$ погружения $g$ изоморфно сумме
Уитни $(\frac{n+1}{4}-8)$ экземпляров линейного расслоения
$\kappa$ над $M$, т.е. $\nu(g)=(\frac{n+1}{4}-8)\kappa$ (в
частности, $w_1(M)=0$, поскольку коразмерность погужения $g$ четна
и $w_1(M)=(\frac{n+1}{4}-8)w_1(\kappa)=0$). Тогда в предположении
$n \ge 63$ (т.е. при $l \ge 6$) справедливо равенство $\langle
w_1(\kappa)^{ \dim(M)};[M] \rangle=0$, где $w_1$ -- первый
характеристический класс Штифеля-Уитни.
\[  \]

Новое в доказательстве Основной Теоремы состоит в применении
принципа геометрического контроля, который позволяет в классе
кобордизма погружения находить погружение, у которого
подмногообразие кратных точек самопересечения имеет более простую
структурную группу нормального расслоения. Напомню, что
структурная группа нормального расслоения к многообразию кратных
точек самопересечения отвечает за перестановку листов исходного
многообразия и изменение векторов скошенного оснащения при обходе
по замкнутому пути на многообразии кратных точек.

При помощи стандартных рассуждений теории погружений из Основной
Теоремы можно вывести следующее.

\subsubsection*{Основное Следствие}

Пусть $g: M^{n-1} \looparrowright \R^n$ -- произвольное гладкое
погружение замкнутого многообразия $M^{n-1}$ (которое, вообще
говоря, ориентируемым и связным не предполагается). Тогда в
предположении $n=2^l-1$, $n \ge 63$ (т.е. при $l \ge 6$),
справедливо равенство $\langle w_1(M)^{n-1};[M] \rangle =0$.
\[  \]

\begin{remark*}
 Эквивалентность предыдущего утверждения  и теоремы Адамса (при
 $l \ge 4$) доказана в [E1,La].
\end{remark*}

Отметим, что в топологии существуют теоремы, близкие по
формулировке к Теореме Адамса. Как правило, эти теоремы являются
следствиями Теоремы Адамса. Иногда сами эти теоремы могут быть
доказаны альтернативными и более простыми методами. Как замечает
С.П.Новиков в обзоре [N], к таким теоремам относится Теорема
Милнора-Ботта о том, что касательное $n$-мерное расслоение к
стандартной сфере $S^n$ тривиально тогда и только тогда, когда
$n=1,3,7$. Это теорема была открыта в работе [B-M]. Элегантная
модификация известного доказательства была недавно получена в [F].

Остановимся на структуре работы. В разделе 1 напоминаются основные
определения и конструкции теории погружений. Результаты этого
раздела формально являются новыми, но легко получаются известными
методами. В разделе 2 Основная Теорема
переформулирована с использованием обозначений из раздела 1
(Теорема 5) и представлены основные этапы ее доказательства. В
разделе 3 проводится основной этап доказательства, а именно,
доказывается первая Основная Лемма 24, которая была высказана
С.А.Мелиховым как гипотеза и обсуждалась на семинаре проф. О.Саеки
в 2006 году, а также аналогичная ей вторая Основная Лемма 25.

Работа была написана при неоднократных обсуждениях на семинаре по
алгебраической топологии под руководством проф. М.М.Постникова. В
настоящей версии работы исправлена ошибка в формулировке Основной
Леммы 3 из [Akh], которая приводилась без доказательства.

Автор благодарит за многочисленные обсуждения проф.
В.М.Бухштабера, проф. В.А.Васильева, проф. П. Ландвебера, проф.
А.С.Мищенко, проф. О.Саеки, проф. А.Б.Скопенкова, проф.
Ю.П.Соловьёва, проф. А.В.Чернавского, проф. Е.В.Щепина, проф.
П.Дж.Экклза, Н.Бродского, С.А.Мелихова, Р.Р.Садыкова,
М.Б.Скопенкова.

\section{Предварительные редукции}
Напомним определения групп кобордизма оснащенных погружений  в
евклидово пространство, которые являются частными случаями более
общей конструкции, изложенной в книге [K1] на стр. 55 и в разделе
10. Связь с конструкцией Понтрягина-Тома объясняется в [A-E].

Пусть $f:  M^{n-1}  \looparrowright  \R^{n}$--  гладкое
погружение, причем $(n-1)$-мерное многообразие $M^{n-1}$ замкнуто
но, вообще говоря, неориентировано и несвязно. На пространстве
таких погружений введем отношение кобордантности. Скажем, что два
погружения $f_0$, $f_1$  связаны кобордизмом, $f_0  \sim  f_1$,
если существует погружение
 $ \Phi  :(  W^{n},  \partial  W  =  M^{n-1}_0  \cup
M^{n-1}_1) \looparrowright  (\R^{n}  \times  [0;1];\R^n \times
\{0,1\})$    такое, что выполнены граничные условия $f_i = \Phi
\vert_{ M^{n-1}_i} : M^{n-1}_i \looparrowright \R^{n} \times \{i\}
$, $i=0,1$ и, кроме того, требуется, чтобы погружение $\Phi$ было ортогонально к $\R \times \{0,1\}$.

Множество классов кобордантности погружений образует
абелеву группу относительно операции дизъюнктного объединения
погружений. Например, нулевой элемент этой группы представляется пустым погружением, а
элемент, обратный к данному погружению $f_0$ представлен композицией $S \circ f_0$, где $S$--зеркальная симметрия
пространства $\R^n$.
Эту группу обозначим через $Imm^{sf}(n-1,1)$. По теореме Уоллеса [Wa] указанная группа
изоморфна стабильной гомотопической группе $\lim_{k \to \infty} \pi_{k+n}(\Sigma^k\RP^{\infty})$.

Погружение $f$ определяет изоморфизм нормального расслоения
многообразия $M^{n-1}$ и ориентирующего линейного расслоения
$\kappa$, т.е. изоморфизм $D(f): T(M^{n-1}) \oplus \kappa =
n\varepsilon$, где $\varepsilon$--тривиальное линейное расслоение
над $M^{n-1}$. В похожих конструкциях теории хирургии гладких
многообразий всегда рассматривается стабильный изоморфизм
нормального расслоения многообразия $M^{n-1}$ и ориентирующего
линейного расслоения $\kappa$, т.е. изоморфизм $ T(M^{n-1}) \oplus
\kappa \oplus N\varepsilon= (n+N)\varepsilon$, где $N>>n$.
Пользуясь теоремой Хирша, легко проверить, что если двум
погружениям $f_1, f_2$ отвечают изоморфизмы $D(f_1)$, $D(f_2)$,
которые при стабилизации принадлежат одному классу (стабильного)
изоморфизма, то рассматриваемые погружения $f_1$, $f_2$
оказываются регулярно кобордантными и даже регулярно
конкордантными (но, вообще говоря, могут оказаться не регулярно
гомотопными).

Нам потребуется также группа $Imm^{sf}(n-k,k)$. Элемент этой
группы представлен тройкой $(f, \kappa, \Xi)$, где $f: M^{n-k}
\looparrowright \R^{n}$ -- погружение замкнутого многообразия,
$\kappa: E(\kappa) \to M^{n-k}$ линейное расслоение (для
сокращения обозначений далее линейное (одномерное) расслоение и
его характеристический класс в $H^1(M^{n-k};\Z/2)$ мы обозначаем
той же буквой), $\Xi$ скошенное оснащение нормального расслоения
погружения при помощи расслоения $\kappa$, т.е. изоморфизм
нормального расслоения погружения $f$ и расслоения $k \kappa$.
Тройка $(f, \kappa, \Xi)$ будет называться скошенно-оснащенным погружением.
В случае нечетного $k$ линейное расслоение $\kappa$ оказывается
ориентирующим над $M^{n-k}$, и обязательно  $\kappa=w_1(M^{n-k})$.

Два элемента группы кобордизма, представленные тройками
$(f_1,\kappa_1,\Xi_1)$ $(f_2,\kappa_2,\Xi_2)$, равны, если
погружения $f_1$, $f_2$ кобордантны, при этом требуется, чтобы
погружение кобордизма было скошенно-оснащено и требуется
согласование скошенного оснащения на кобордизме с заданными
скошенными оснащениями на компонентах границы. Заметим, что при
$k=1$ новое определение групп $Imm^{sf}(n-k,k)$ совпадает с
первоначальным.

Определим гомоморфизм: $$J^{sf}: Imm^{sf}(n-1,1) \to
Imm^{sf}(n-k,k),$$ который называется гомоморфизмом перехода в
коразмерность $k$. Рассмотрим многообразие $M'^{n-1}$, погружение
которого $f: M'^{n-1} \looparrowright \R^n$ представляет элемент
первой группы, и рассмотрим классифицирующее отображение $\kappa':
M'^{n-1} \to \RP^a$ в вещественное проективное пространство
высокой размерности ($a = n$ достаточно), представляющее
когомологический класс $w_1(M'^{n-1})$. Рассмотрим стандартное
подпространство $\RP^{a-k+1} \subset \RP^{a}$ коразмерности $k-1$.
Предположим, что отображение $\kappa'$ трансверсально вдоль
выбранного подпространства и определим подмногообразие $M^{n-k}
\subset M'^{n-1}$ как полный прообраз этого подпространства при
нашем отображении, $M^{n-k} = \kappa'^{-1}(\RP^{a-k+1})$.
Определим погружение $f: M^{n-k} \looparrowright \R^n$ как
ограничение погружения $f'$ на заданное подмногообразие. Заметим,
что погружение $f: M^{n-k} \looparrowright \R^n$ допускает
естественное скошенное оснащение. Действительно, нормальное
расслоения подмногообразия $M^{n-k} \subset M'^{n-1}$ естественно
изоморфно расслоению $(k-1) \kappa$, где $\kappa =\kappa'
\vert_{M}$ (здесь и далее, если обозначение многообразия
используется в нижнем индексе, то верхний индекс размерности
многообразия опускается). Указанный изоморфизм $\Xi$ определяется
стандартным скошенным оснащением нормального расслоения к
подмногообразию $\RP^{a-k+1}$ в многообразии $\RP^a$, которое
переносится на подмногообразие $M^{n-k} \subset M'^{n-1}$,
поскольку предполагалось, что $\kappa'$  регулярно вдоль
$\RP^{a-k+1}$. Еще одно прямое слагаемое в скошенном оснащении
$\Xi$ нормального расслоения погружения $f$ соответствует
нормальному одномерному расслоению погружения $f'$. Это расслоение
служит ориентирующим расслоением для $M'^{n-1}$, поэтому его
ограничение на $M^{n-k}$ совпадает с $\kappa$. Гомоморфизм
$J^{sf}$ переводит элемент, представленный погружением $f'$ в
элемент, представленный тройкой $(f,\kappa,\Xi)$. Из элементарных
геометрических соображений, использующих лишь понятие
трансверсальности вытекает, что гомоморфизм $J^{sf}$ корректно
определен, поскольку произвол в его конструкции приводит к тройке
из того же класса нормального кобордизма.

 При условии, что
погружение $f: M^{n-k} \looparrowright \R^n$ общего положения,
подмножество в $\R^n$ точек самопересечения погружения $f$
обозначается через $\Delta = \Delta(f)$, $\dim(\Delta)=n-2k$. Это
подмножество определяется по формуле:
\begin{eqnarray}\label{Delta}
\Delta=\{ x \in \R^n : \exists x_1, x_2 \in M^{n-k}, x_1 \ne
x_2, f(x_1)=f(x_2)=x\},
\end{eqnarray}
оно является замкнутым. Определим $\bar
\Delta \subset M^{n-k}$ по формуле $\bar \Delta = f^{-1}(\Delta)$.

Напомним стандартное определение многообразия двукратных точек
самопересечения данного погружения общего положения и параметризующего погружения
многообразия двукратных точек самопересечения.
Пространство $N$ определяется формулой
\begin{eqnarray}\label{N}
N=\{[(x_1,x_2)] \in (M^{n-k} \times M^{n-k})/T' : x_1 \ne x_2,
f(x_1)=f(x_2)\}
\end{eqnarray}
 ($T'$ -- инволюция, переставляющая координаты
сомножителей), а его каноническое накрывающее определяется
формулой
\begin{eqnarray}\label{barN}
\bar N=\{ (x_1,x_2) \in M^{n-k} \times M^{n-k} : x_1 \ne x_2,
f(x_1)=f(x_2)\}.
\end{eqnarray}
 В предположении о том, что погружение $f$ общего
положения, пространство $N$ является гладким многообразием
размерности $\dim(N)=n-2k$. Это многообразие обозначается через
$N^{n-2k}$.

Погружение $\bar g: \bar N^{n-2k} \looparrowright M^{n-k}$,
параметризующее $\bar \Delta$, определяется формулой $\bar g = \pi
\vert_{\bar N}$, где $\pi: M^{n-k} \times M^{n-k} \to M^{n-k}$ --
стандартная проекция на первый сомножитель. Погружение $g:
N^{n-2k} \looparrowright \R^n$, параметризующее $\Delta$,
определяется формулой $g([x_1,x_2])=f(x_1)$. Заметим, что
параметризующие погружения $g, \bar g$ не являются, вообще говоря,
погружениями общего положения. Определено двулистное накрытие $p:
\bar N^{n-2k} \to N^{n-2k}$, при этом $g \circ p = f \circ \bar g$. Это
накрытие назовем каноническим накрытием над многообразием точек
самопересечения.

\begin{definition}\label{def1}
Пусть $(f, \kappa, \Xi)$ представляет элемент из
$Imm^{sf}(n-k,k)$. Определим гомоморфизм:
$$ h_k : Imm^{sf}(n-k,k) \to \Z/2,$$
называемый инвариантом Хопфа, по формуле:
$$h_k([f,\kappa,\Xi]) =
\langle \kappa^{n-k};[\bar M^{n-k}] \rangle. $$
\end{definition}

Определения стабильного инварианта Хопфа (в смысле Определения 1)
при различных значениях $k$ согласованы между собой. Сформулируем
это в виде отдельного утверждения.

\begin{proposition}\label{prop2}
При гомоморфизме $J^{sf}: Imm^{sf}(n-1,1) \to Imm^{sf}(n-k,k)$
стабильный инвариант Хопфа сохраняется, т.е. инвариант $h_1:
Imm^{sf}(n-1,1) \to \Z/2$ и инвариант $h_k: Imm^{sf}(n-k,k) \to
\Z/2$ связаны между собой по формуле:
\begin{eqnarray}\label{1}
 h_1 = h_k \circ J^{sf}.
\end{eqnarray}
\end{proposition}

\subsubsection*{Доказательство Предложения $\ref{prop2}$}

Пусть $f: M^{n-k} \looparrowright \R^n$ --  погружение со
скошенным оснащением $\Xi$ и c характеристическим классом $\kappa
\in H^1(M^{n-k};\Z/2)$, представляющее элемент из
$Imm^{sf}(n-k,k)$, причем $J^{sf}([f'])=[f,\kappa,\Xi]$ для
некоторого элемента $[f'] \in Imm^{sf}(n-1,1)$. По определению
$h_k([f,\kappa,\Xi])=\langle \kappa^{n-k};[M^{n-k}] \rangle$.

C другой стороны, $M^{n-k} \subset M'^{n-1}$ представляет цикл,
двойственный в смысле Пуанкаре коциклу $\kappa'^{k-1} \in
H^{k-1}(M'^{n-1};\Z/2)$. Формула $\ref{1}$ справедлива, поскольку
$\langle \kappa'^{n-1};[M'^{n-1}]\rangle =\langle
\kappa^{n-k};[M^{n-k}]\rangle $. Действительно, образ
характеристического отображения $\kappa': M'^{n-1} \to
\RP^{\infty}$ ($\kappa: M^{n-k} \to \RP^{\infty}$), не ограничивая
общности, лежит в остове классифицирующего пространства
размерности $n-1$ ($n-k$). Характеристическое число $\langle
\kappa'^{n-1};[M'^{n-1}]\rangle$ ($\langle
\kappa^{n-k};[M^{n-k}]\rangle$) совпадает со степенью
$\deg(\kappa')$ ($\deg(\kappa)$) классифицирующего отображения
$\kappa': M'^{n-1} \to \RP^{n-1}$ ($\kappa: M^{n-k} \to
\RP^{n-k}$), которая рассматривается по модулю 2 и определяется
как четность числа прообразов регулярного значения отображения.
При этом значение $\deg(\kappa')$ ($\deg(\kappa)$) не зависит от
выбора отображения $\kappa'$ ($\kappa$) в указанный остов.
(Напомню, что классы когомологий и их характеристические
отображения обозначаются  одинаково.) Степени $\deg(\kappa')$ и
$\deg(\kappa)$ совпадают, поскольку регулярное значение можно
выбрать общим для рассматриваемых отображений. Предложение 2
доказано.
\[  \]

Сформулируем другое (эквивалентное) определение инварианта Хопфа
(в предположении $n-2k>0$).

\begin{definition}\label{proba}
Пусть $(f, \kappa, \Xi)$ представляет элемент из
$Imm^{sf}(n-k,k)$, $n-2k>0$.
 Пусть $N^{n-2k}$ -- многообразие двукратных точек погружения $f: M^{n-k} \looparrowright \R^n$,
 $\bar N$--каноническое двулистное накрытие над
$N$, $\kappa_{\bar N} \in H^1(\bar N;\Z/2)$ индуцирован из $\kappa
\in H^1(M^{n-k};\Z/2)$ погружением $\bar g: \bar N^{n-2k}
\looparrowright M^{n-k}$.

Определим значение гомоморфизма $h_k : Imm^{sf}(n-k,k) \to \Z/2$
по формуле:
$$h_k([f,\kappa,\Xi]) = \langle \kappa^{n-2k}_{\bar
N};[\bar N^{n-2k}] \rangle. $$
\end{definition}
\[  \]

Следующее утверждение устанавливает эквивалентность Определений 1
и $\ref{proba}$.

\begin{proposition}\label{prop4}
В условиях Определения 3 справедлива формула:
\begin{eqnarray}\label{1.400}
\langle \kappa^{n-2k}_{\bar N};[\bar N^{n-2k}] \rangle =\langle
\kappa^{n-k};[M^{n-k}]\rangle.
\end{eqnarray}
\end{proposition}

\subsubsection*{Доказательство Предложения $\ref{prop4}$}
Пусть $f: M^{n-k} \looparrowright \R^n$ --  погружение со
скошенным оснащением $\Xi$ и c характеристическим классом $\kappa
\in H^1(M^{n-k};\Z/2)$, представляющее элемент из
$Imm^{sf}(n-k,k)$. Пусть $N^{n-2k}$--многообразие двукратных точек
самопересечения погружения $f$, $g: N^{n-2k} \looparrowright
\R^n$--параметризующее погружение, $\bar N^{n-2k} \to
N^{n-2k}$--каноническое двулистное накрытие. Рассмотрим образ
фундаментального класса $\bar g_{\ast}([\bar N^{n-2k}]) \in
H_{n-2k}(M^{n-k};\Z/2)$ при погружении $\bar g: \bar N^{n-2k}
\looparrowright M^{n-k}$ и обозначим через $m \in
H^{k}(M^{n-k};\Z/2)$--когомологический класс, двойственный в
смысле Пуанкаре классу гомологий $\bar g_{\ast}([\bar N^{n-2k}])$.
Рассмотрим также когомологический эйлеров класс нормального
расслоения погружения $f$, который обозначим через $e \in
H^k(M^{n-k};\Z/2)$.

По Теореме Герберта (см. [E-G], Теорема 1.1) справедлива формула:
\begin{eqnarray}\label{1.500}
e=m.
\end{eqnarray}
Поскольку эйлеров класс $e$ нормального расслоения $k\kappa$
погружения $f$ равен $\kappa^k$ (классы когомологий и
соответствующие линейные расслоения обозначаются одинаково), то
цикл $\bar g_{\ast}([\bar N]) \in H_{n-2k}(M^{n-k};\Z/2)$,
двойственен в смысле Пуанкаре коциклу $\kappa^{k} \in
H^{k}(M^{n-k};\Z/2)$. Поэтому формула $(\ref{1.400})$ и
Предложение $\ref{prop4}$ доказано.
\[  \]

Более удобно переформулировать Предложение 4 (в несколько более общей форме) на языке коммутативных диаграмм.
Приступим к формулировке соответствующих определений.

Пусть $g: N^{n-2} \looparrowright \R^n$-- погружение многообразия точек самопересечения
погружения $f: M^{n-1} \looparrowright \R^n$ коразмерности $1$. Определено нормальное
двумерное расслоение погружения $g$,  которое мы будем обозначать через $\nu_N :
E(\nu_N) \to N^{n-2}$. (Пространство дискового  расслоения, ассоциированного с
векторным расслоением $\nu_N$, диффеоморфно регулярной замкнутой трубчатой окрестности
погруженния $g$.)

Это расслоение снабжено дополнительной
структурой по сравнению с произвольным векторным расслоением, а
именно, его структурная группа как $O(2)$-расслоения допускает
редукцию к дискретной группе диэдра, которую мы обозначим через
$\D_4$. Это группа восьмого порядка, она определяется как группа
перемещений плоскости, которые переводят стандартную пару
координатных осей  в себя (возможно, с изменением ориентации и
порядка).

В стандартном копредставлении группа $\D_4$ задается двумя
образующими $a,b$, которые связаны между собой соотношениями
$\{a^4=b^2=1, [a;b]=a^2\}$. Образующая $a$ представлена поворотом
плоскости на угол $\frac{\pi}{2}$, образующая $b$ представлена
симметрией относительно биссектриссы первого координатного угла.
Заметим, что элементу $ba$ (произведение будем записывать по
правилу композиции $b \circ a$ преобразований из $O(2)$)
соответствует симметрия относительно первой координатной оси.

\subsubsection*{Замечание}

В книге [A-M] описана структура когомологий по модулю 2
классифицирующего пространства группы $\D_4$ как модуля над
алгеброй Стинрода. Инвариант Хопфа, который определяется ниже как
характеристическое число на многообразии точек самопересечения
скошенно-оснащенного погружения, это
вычисление не использует.

\subsubsection*{Структурная группа нормального расслоения многообразия
точек самопересечения для погружения $f: M^{n-k} \looparrowright
\R^n$, при $k=1$.}

Используем условие трансверсальности самопересечения  погружения
$f: M^{n-1} \looparrowright \R^n$. Пусть $N^{n-2}$--многообразие самопересечения погружения $f$,
$g: N^{n-2} \looparrowright \R^n$--параметризующее погружение.
В слое $E(\nu_N)_x$ нормального
расслоения $\nu_N$ над точкой $x \in N$ фиксирована
неупорядоченная пара координатных осей. Эти оси образованы
касательными к кривым пересечения слоя $E(\nu_N)_x$ с двумя
вложенными листами погруженного многообразия, пересекающимися
трансверсально в окрестности этой точки. По построению расслоение
$\nu_N$ имеет структурную группу $\D_4 \subset O(2)$.

 Над пространством $K(\D_4,1)$
определено универсальное 2-мерное $\D_4$--расслоение, которое
обозначим через $\psi: E(\psi) \to K(\D_4,1)$. Скажем, что
отображение $\eta: N \to K(\D_4,1)$ является классифицирующим для
расслоения $\nu_N$, если определен изоморфизм $\Xi: \eta^{\ast}(\psi) = \nu_N$
обратного
образа $\eta^{\ast}(\psi)$ расслоения $\psi$ и расслоения нормального расслоения $\nu_N$ погружения $g$.
В дальнейшем само расслоение и его классифицирующее отображение
будут обозначаться одинаково и в рассматриваемом случае можно
записать $\eta=\nu_N$. Изоморфизм $\Xi$ будем называть $\D_4$--оснащением погружения $g$, а отображение $\eta$ будем называть характеристическим классом $\D_4$--оснащения $\Xi$.

\begin{remark*}

В действительности, мы описали лишь часть более общей конструкции.
Структурная группа $s$-мерного нормального расслоения к
подмногообразию $N_s$ точек самопересечения кратности $s$
погружения $f$ допускает редукцию к структурной группе $\Z/2 \int
\Sigma(s)$ -- сплетению циклической группы $\Z/2$ с группой
подстановок множества из $s$ элементов (см., например, [E1]).
\end{remark*}
\[  \]

 Пусть тройка $(f,\kappa,\Xi)$, где $f: M^{n-k}
\looparrowright \R^n$ -- погружение, $\Xi$ -- скошенное оснащение
нормального расслоения погружения $f$ c характеристическим классом $\kappa \in H^1(M;\Z/2)$, представляет элемент группы $Imm^{sf}(n-k,k)$.
 Нам потребуется обобщение предыдущего построения, которое относилось к случаю
 $k=1$. Напомним описание (см., например, [E1], [E-G])
структурной группы нормального расслоения $\nu_N$ над
многообразием точек самопересечения произвольного погружения
общего положения произвольной коразмерности $k$, $k \le
[\frac{n}{2}]$.

\begin{proposition}\label{prop5}
Нормальное $2k$-мерное расслоение $\nu_{N}$ погружения $g$
представимо в виде прямой суммы $k$ изоморфных копий двумерного
расслоения $\eta$ над $N^{n-2k}$, причем каждое двумерное
расслоение имеет структурную группу $\D_4$ и классифицируется
отображением $\eta: N^{n-2k} \to K(\D_4,1)$ (аналогичное
предложение доказано в {\rm[Sz2]}).
\end{proposition}

\subsubsection*{Доказательство Предложения $\ref{prop5}$}
Пусть $x \in N^{n-2k}$ точка на многообразии двукратных точек
погружения $f: M^{n-k} \looparrowright \R^n$. Обозначим через
$\bar x_1, \bar x_2 \in \bar N^{n-2k} \looparrowright M^{n-k}$ два
прообраза этой точки на каноническом накрывающем многообразии.
Ортогональное дополнение в пространстве $T_(g(x))(\R^n)$ к
подпространству $g_{\ast}(T_x(N^{n-2k}))$ является слоем
нормального расслоения $E(\nu_N)$ погружения $g$ над точкой $x \in
N^{n-2k}$. Этот слой представлен как прямая сумма двух линейных
пространств $E(\nu_N)_{x} = \bar E_{x,1} \oplus \bar E_{x,2}$, где
подпространство $\bar E_{x,i} \subset E(\nu_N)_x$ является слоем
нормального расслоения погружения $f$ в точке $\bar x_i$.

Подпространство $\bar E_{x,i}$ слоя канонически представлено в
виде прямой суммы $k$ упорядоченных линейных подпространств $\bar
E(\kappa_{x,j,i})$, $j=1, \dots, k$, $i=1,2$, причем фиксировано
семейство попарных изоморфизмов между этими пространствами,
поскольку нормальное расслоение погружения $f$ было снабжено
скошенным оснащением в коразмерности $k$. Сгруппируем слои с
одинаковыми номерами $j$ в двумерный подслой слоя $E(\nu_N)_x$.
Получим разложение слоя $E(\nu_N)_x$ над каждой точкой $x \in
N^{n-2k}$ в прямую сумму $k$ экземпляров попарно изоморфных
двумерных подпространств. Проведенная конструкция непрерывно
зависит от выбора точки $x$ и может быть проведена одновременно
для каждой точки базы $N^{n-2k}$. В результате получится требуемое
разложение расслоения $\nu_N$ в сумму Уитни некоторого числа
канонически изоморфных между собой двумерных подрасслоений. Каждое
двумерное прямое слагаемое классифицируется отображением $\eta: N
\to K(\D_4,1)$, что доказывает Предложение $\ref{prop5}$.
\[  \]

\begin{definition}\label{def6}
Определим группу кобордизма погружений
$Imm^{\D_4}(n-2k,2k)$, в предположении $n>2k$.
Пусть $(g,\Psi,\eta)$-- тройка, определяющая $\D_4$--оснащенное
погружение в коразмерности $2k$. Здесь $g: N^{n-2k} \looparrowright
\R^n$ -- погружение, $\eta: N^{n-2k} \to K(\D_4,1)$--
характеристический класс $\D_4$--оснащения $\Psi$. Отношение кобордизма троек стандартно.
\end{definition}

\begin{lemma}\label{lemma7}
В предположении $k_1 < k$, $2k < n$, определена коммутативная диаграмма гомоморфизмов групп:
\begin{eqnarray}\label{1.5}
\begin{array}{ccccc}
Imm^{sf}(n-k_1,k_1) & \stackrel{J^{sf}}{\longrightarrow} & Imm^{sf}(n-k,k) & \stackrel{h_k}{\longrightarrow} & \Z/2\\
\downarrow \delta_{k_1}& &\downarrow \delta_{k} && \| \\
Imm^{\D_4}(n-2k_1,2k_1) & \stackrel{J^{\D_4}}{\longrightarrow} &
Imm^{\D_4}(n-2k,2k) & \stackrel{h^{\D_4}_k}{\longrightarrow} &
\Z/2.
\end{array}
\end{eqnarray}
\end{lemma}

\subsubsection*{Доказательство Леммы $\ref{lemma7}$}
Определим гомоморфизмы в диаграмме $(\ref{1.5})$. Гомоморфизм
$$Imm^{sf}(n-k_1,k_1)  \stackrel{J^{sf}}{\longrightarrow}
Imm^{sf}(n-k,k)$$
 определяется дословно также, как был определен
гомоморфизм перехода в коразмерность $k$ для случая $k_1=1$.
Гомоморфизм $$Imm^{sf}(n-k,k)
\stackrel{\delta_{k}}{\longrightarrow}  Imm^{\D_4}(n-2k,2k)$$
переводит класс кобордизма тройки $(f,\kappa,\Xi)$ в класс
кобордизма тройки $(g,\eta,\Psi)$, где $g: N^{n-2k}
\looparrowright \R^n$ -- параметризующее погружение многообразия
точек самопересечения погружения $f$ (предполагается, что
погружение $f$ самопересекается трансверсально), $\Psi$--это
$\D_4$--оснащение нормального расслоения погружения $g$,
$\eta$--характеристический класс $\D_4$-- оснащения $\Psi$.

Приступим к определению гомоморфизма
$$Imm^{\D_4}(n-2k,2k)
\stackrel{h^{\D_4}_k}{\longrightarrow}  \Z/2,$$
который назовем
диэдральным инвариантом Хопфа. Определим подгруппу
\begin{eqnarray}\label{1.6}
\I_c \subset \D_4,
\end{eqnarray}
порожденную преобразованиями плоскости, сохраняющими собственное
подпространство каждого базисного вектора. Группа $\I_c$ является
элементарной абелевой группой ранга 2. Определим гомоморфизм
\begin{eqnarray}\label{1.7}
l^{[2]}: \I_c \to \Z/2,
\end{eqnarray}
отвечающий за преобразование в собственном пространстве первого
базисного вектора. Подгруппа $(\ref{1.6})$ имеет индекс 2 и
определено 2-листное накрытие:
\begin{eqnarray}\label{1.8}
 K(\I_c,1) \to K(\D_4,1),
\end{eqnarray}
индуцированное этой подгруппой.

Обозначим через
\begin{eqnarray}\label{1.9}
\bar N^{n-2k} \to N^{n-2k}
\end{eqnarray}
2-листное накрытие, индуцированное характеристическим классом
$\eta: N^{n-2k} \to K(\D_4,1)$ из накрытия $(\ref{1.8})$.
Определен характеристический класс:
$$  \bar \eta^{sf} = l \circ \bar \eta^{\I_c} : \bar N^{n-2k} \to K(\Z/2,1), $$
где $\bar \eta^{\I_c}: \bar N^{n-2k} \to  K(\I_c,1)$-- двулистное
накрытие над классифицирующим отображением $\eta: N^{n-2k} \to
K(\D_4,1)$, индуцированное накрытиями $(\ref{1.8})$, $(\ref{1.9})$
над пространствами обаза и прообраза отображения $\eta$
соответственно.

Определим значение гомоморфизма $Imm^{\D_4}(n-2k,2k)
\stackrel{h^{\D_4}_k}{\longrightarrow}  \Z/2$ по формуле:
\begin{eqnarray}\label{1.99}
h^{\D_4}_k([g,\eta,\Psi]) = \left\langle  (\bar
\eta^{sf})^{n-2k};[\bar N^{n-2k}] \right\rangle.
\end{eqnarray}

Диаграмма $(\ref{1.5})$ определена. Коммутативность правого
квадрата диаграммы вытекает из прямых геометрических рассуждений.
Коммутативность правого квадрата диаграммы вытекает при $k_1=1$ из
Предложения 5, а при произвольном $k_1$ доказательство аналогично.
Лемма $\ref{lemma7}$ доказана.
\[  \]

Нам потребуется эквивалентное определение диэдрального инварианта Хопфа.
Рассмотрим подгруппу в группе $O(4)$ преобразований, переводящих
множество векторов $(\e_1,\e_2,\e_3,\e_4)$ стандартного базиса в себя, быть может, изменяя направления некоторых
векторов и, кроме того, сохраняющих пару 2-мерных подпространств, порожденных
базисными векторами $(\e_1,\e_2)$, $(\e_3,\e_4)$.
При этом указанные 2-мерные подпространства могут, вообще говоря, переставляться между собой.
Обозначим подгруппу таких преобразований через $\Z/2^{[3]}$.
Определим цепочку подгрупп индекса 2:
\begin{eqnarray}\label{1.10}
\H_{\bar c} \subset \H_c \subset \Z/2^{[3]}.
\end{eqnarray}
Подгруппа $\H_c$ определена как подгруппа преобразований,
инвариантных на каждом 2-мерном подпространстве, порожденном
парами векторов $(\e_1,\e_2)$, $(\e_3,\e_4)$. Эта подгруппа
изоморфна прямой сумме двух экземпляров группы $\D_4$, каждое
прямое слагаемое инвариантно действует в соответствующем 2-мерном
подпространстве. Подгруппа $ \H_{\bar c} \subset \H_c$ определена
как подгруппа преобразований, инвариантных на каждом линейном
подпространстве, порожденных векторами $\e_1$, $\e_2$.

Пусть $(g, \eta, \Psi)$ представляет элемент из
$Imm^{\D_4}(n-2k,2k)$ в предположении, что $n-4k>0$ и что $g$ является погружением общего положения.
 Пусть $L^{n-4k}$ -- многообразие двукратных точек погружения $g: N^{n-2k} \looparrowright \R^n$. Определена
башня 2-листных накрытий:
\begin{eqnarray}\label{1.12}
L_{\H_{\bar c}}^{n-4k} \to L_{\H_c}^{n-4k} \to L^{n-4k}
\end{eqnarray}
в результате следующей конструкции. Определено параметризующее
погружение $h: L^{n-4k} \looparrowright \R^n$. Нормальное
расслоение погружения $h$ обозначим через $\nu_L$. Это расслоение
крассифицируется отображением
 $\zeta: L^{n-4k} \to K(\Z/3^{[3]},1)$.
Цепочка подгрупп $(\ref{1.10})$ индуцирует башню 2-листных накрытий классифицирующих пространств:
\begin{eqnarray}\label{1.11}
K(\H_{\bar c},1) \subset K(\H_c,1) \subset K(\Z/2^{[3]},1)
\end{eqnarray}
над пространством-образом классифицирующего отображения $\zeta$ и
башню 2-листных накрытий $(\ref{1.12})$ над
пространством-прообразом отображения $\zeta$. Накрытие
$L_{\H_{\bar c}} \to  L^{n-4k}$, определенное формулой
$(\ref{1.12})$, будем называть каноническим 4-листным накрытием
над многообразием точек самопересечения погружения $g$.

Определено отображение классифицирующих пространств:
\begin{eqnarray}\label{1.41}
K(\H_{\bar c},1) \to K(\Z/2,1),
\end{eqnarray}
индуцированное гомоморфизмом
\begin{eqnarray}\label{1.51}
l^{[3]}: \H_{\bar c} \to \Z/2,
\end{eqnarray}
который отвечает за преобразование линейного собственного
подпространства вектора $\e_1$. Определены классифицирующее
отображение $\bar \zeta^{\H_{\bar c}}: \bar L_{\H_{\bar c}}^{n-4k}
\to K(\H_{\bar c},1)$ в результате перехода к 4-листному накрытию
над отображением $\zeta$ и характеристическое отображение $\bar
\zeta^{sf}: L^{n-4k}_{\H_{\bar c}} \to K(\Z/2,1)$ в результате
композиции классифицирующего отображение $\bar \zeta^{\H_{\bar
c}}$ с отображением $(\ref{1.41})$.

\begin{proposition}\label{prop8}
Предположим, что $\D_4$--оснащенное погружение $(g, \eta, \Psi)$,
представляющее элемент из группы $Imm^{\D_4}(n-2k,2k)$, лежит в
образе гомоморфизма
$$ \delta_k: Imm^{sf}(n-k,k) \to Imm^{\D_4}(n-2k,2k).$$
Тогда справедлива формула:
\begin{eqnarray}\label{1.60}
\langle (\bar \zeta^{sf})^{n-4k}; [\bar L_{\H_{\bar c}}^{n-4k}]
\rangle =\langle (\bar \eta^{sf})^{n-2k};[\bar N^{n-2k}]\rangle.
\end{eqnarray}
\end{proposition}

\subsubsection*{Доказательство Предложения $\ref{prop8}$}
Пусть $g: N^{n-2k} \looparrowright \R^n$ --  погружение с
$\D_4$--оснащением $\Psi$ и c характеристическим классом $\eta:
N^{n-2k} \to K(\D_4,1)$, представляющее элемент из
$Imm^{\D_4}(n-2k,2k)$ из образа $\delta_k$. Пусть
$L^{n-4k}$--многообразие двукратных точек самопересечения
погружения $g$, $h: L^{n-4k} \looparrowright
\R^n$--параметризующее погружение, $\bar L^{n-4k} \to
L^{n-4k}$--каноническое двулистное накрытие. Рассмотрим образ
фундаментального класса $\bar h_{\ast}([\bar L^{n-4k}]) \in
H_{n-4k}(N^{n-2k};\Z/2)$ при погружении $\bar h: \bar L^{n-4k}
\looparrowright N^{n-2k}$ и обозначим через $m \in
H^{2k}(N^{n-2k};\Z/2)$--когомологический класс, двойственный в
смысле Пуанкаре классу гомологий $\bar h_{\ast}([\bar L^{n-4k}])$.
Рассмотрим также когомологический эйлеров класс нормального
расслоения погружения $g$, который обозначим через $e \in
H^{2k}(N^{n-2k};\Z/2)$.

По Теореме Герберта (см. [E-G], Теорема 1.1) справедлива формула
$(\ref{1.500})$.
Обозначим через $\bar e \in H^{2k}(\bar N^{n-2k}_{sf};\Z/2)$,
$\bar m \in H^{2k}(\bar N^{n-2k}_{sf};\Z/2)$ -- образы коциклов
$e$, $m$ соответственно при каноническом двулистном накрытии $\bar
N^{n-2k}_{sf} \to N^{n-2k}$. Справедлива также формула:
$$ \bar e= \bar m, $$
в частности, справедлива формула:
\begin{eqnarray}\label{1.61}
\langle (\bar \eta^{sf})^{n-4k} \bar m;[\bar N_{sf}^{n-2k}]
\rangle = \langle (\bar \eta^{sf})^{n-4k} \bar e;[\bar
N_{sf}^{n-2k}] \rangle.
\end{eqnarray}
 Правая часть формулы $(\ref{1.61})$ равна
$h^{\D_4}_{2k}(J^{\D_4}([g,\eta,\Psi]))$. Это характеристическое
число в силу коммутативности диаграммы $(\ref{1.5})$ совпадает с
$h^{\D_4}_{k}([g,\eta,\Psi]) = \langle (\bar
\eta^{sf})^{n-2k};[\bar N^{n-2k}]\rangle $ в правой части формулы
$(\ref{1.60})$. Поскольку характеристические классы $\bar
\zeta^{sf}: \bar L_{sf}^{n-4k} \to K(\Z/2,1)$ и $\bar \eta^{sf}
\vert_{\bar L_{sf}}$ совпадают, левая часть формулы $(\ref{1.61})$
равна характеристическому числу $\langle (\bar \zeta^{sf})^{n-4k};
[\bar L_{\H_{\bar c}}^{n-4k}] \rangle$ из левой части формулы
$(\ref{1.60})$. Предложение $\ref{prop8}$ доказано.
\[  \]

\section{Доказательство основной теоремы}

Переформулируем Основную Теорему.

\begin{theorem}\label{th9}
При $l \ge 6$ гомоморфизм $h_{2^{l-2}-8}: Imm^{sf}(3 \cdot
2^{l-2}+7,2^{l-2}-8) \to \Z/2$, заданный Определением 3,
тривиален.
\end{theorem}

\subsubsection*{Вывод Основного Следствия из Теоремы $\ref{th9}$}
Рассмотрим гомоморфизм $J^{sf}: Imm^{sf}(n-1,1) \to Imm^{sf}(3
\cdot 2^{l-2}+7,2^{l-2}-8)$. Согласно Предложению 4,
$h_1=h_{2^{l-2}-8} \circ J^{sf}$. Пусть элемент в группе
$Imm^{sf}(n-1,1)$ представлен погружением $f: M^{n-1}
\looparrowright \R^n$, $w_1(M)=\kappa$. Значение $h_k(J^{sf}(f))$,
где $k=2^{l-2}-8$, совпадает с характеристическим числом $\langle
\kappa^{n-1};[M^{n-1}]\rangle$. Используя Теорему $\ref{th9}$,
получаем требуемое.
\[  \]

Для доказательства Теоремы $\ref{th9}$ потребуются ввести
дополнительные определения и обозначения.

\subsubsection*{Определение подгрупп $\I_d \subset \I_a \subset \D_4$, $\I_b \subset \D_4$}
 Обозначим через
 $\I_a = \subset
 \D_4$ циклическую подгруппу порядка 4 индекса 2,  содержащую нетривиальные элементы
 $a,a^2,a^3
\in \D_4$ (т.е. порожденную преобразованием поворота, меняющим
координатные оси). Обозначим через $\I_d
\subset \I_a$ -- подгруппу c нетривиальным элементом $a^2$ (
преобразование центральной симметрии).
Обозначим через $\I_b
\subset \D_4$ -- подгруппу c образующими $a^2$, $ab$ (преобразование симметрии относительно
биссектриссы первой координатной четверти).

Определены
гомоморфизмы включения подгрупп: $i_{d,a}: \I_d \subset
\I_a$  $i_{d,b}: \I_d \subset
\I_b$. В случае, если образ совпадает со всей группой
$\D_4$ соответствующий индекс при гомоморфизмах включений опускается:
$i_d: \I_d \subset \D_4$, $i_a: \I_a \subset \D_4$, $i_b: \I_b \subset \D_4$.
\[  \]

\subsubsection*{Определение подгруппы $i_{\Q_a}: \Q_a \subset \Z/2^{[3]}$}

Обозначим через $\Q_a$ -- группу кватернионов порядка 8.
Эта группа задана копредставлением $\{ \i,\j,\k \vert \i\j = \k = -\j\i,
\j\k=\i = -\k\j, \k\i = \j = -\i\k, \i^2 =\j^2 = \k^2 =-1 \}$. Определено
стандартное представление $\chi_+: \Q_a \to SO(4)$. Представление
$\chi_+$ (матрица действует слева на вектор-столбец) переводит
единичные кватернионы $\i,\j,\k$ в матрицы
\begin{eqnarray}\label{Q a1}
\left(
\begin{array}{cccc}
0 & -1 & 0 & 0 \\
1 & 0 & 0 & 0 \\
0 & 0 & 0 & -1 \\
0 & 0 & 1 & 0 \\
\end{array}
\right),
\end{eqnarray}
\begin{eqnarray}\label{Q a2}
\left(
\begin{array}{cccc}
0 & 0 & -1 & 0 \\
0 & 0 & 0 & 1 \\
1 & 0 & 0 & 0 \\
0 & -1 & 0 & 0 \\
\end{array}
\right),
\end{eqnarray}
\begin{eqnarray}\label{Q a3}
\left(
\begin{array}{cccc}
0 & 0 & 0 & -1 \\
0 & 0 & -1 & 0 \\
0 & 1 & 0 & 0 \\
1 & 0 & 0 & 0 \\
\end{array}
\right).
\end{eqnarray}
Указанное представление определяет подгруппу $i_{\Q_a}: \Q_a \subset
\Z/2^{[3]} \subset O(4)$.

\subsubsection*{Определение подгрупп $\I_d \subset \I_a \subset \Q_a$}

Обозначим через
$i_{\I_d,\Q_a}: \I_d \subset \Q_{a}$ -- центральную подгруппу в кватернионной группе,
которая оказывается также центральной и во всей группе $\Z/2^{[3]}$.

Обозначим через
$i_{\I_a,\Q_a}: \I_a \subset \Q_{a}$ -- подгруппу в кватернионной группе,
порожденную комплексным кватернионом $\i$.

Определены также гомоморфизмы включения
$i_{\I_d}: \I_d \subset \Z/2^{[3]}$,
$i_{\Q_a}: \Q_a \subset \Z/2^{[3]}$
\[  \]

\begin{definition}\label{defcycl}
Структурное отображение $\eta: N^{n-2k} \to K(\D_4,1)$ называется
циклическим, если оно является композицией некоторого отображения
$\mu_a:  N^{n-2k} \to K(\I_a,1)$ и включения $i_a: K(\I_a,1)
\subset K(\D_4,1)$.
\end{definition}

\begin{definition}\label{defQ}
Структурное отображение $\zeta: L^{n-4k} \to K(
\Z/2^{[3]},1)$ называется кватернионным, если оно является композицией
некоторого отображения $\lambda_a:  L^{n-4k} \to K(\Q_a,1)$ и
включения $K(\Q_a,1) \subset K(\Z/2^{[3]},1)$.
\end{definition}

 Ниже нам потребуется конструкция пространств
Эйленберга-Маклейна $K(\I_a,1)$, $K(\Q_a,1)$ и описание
конечномерных остовов этих пространств, которое мы напомним.
Рассмотрим бесконечномерную сферу $S^{\infty}$ (это стягиваемое
пространство), которая определяется как прямой предел
бесконечной цепочки вложений стандартных сфер нечетной
размерности,
$$S^{\infty}=
\lim_{\longrightarrow} \quad (S^1 \subset S^3 \subset \dots \subset
S^{2j-1} \subset S^{2j+1} \subset \dots  \quad) .$$

При этом $S^{2j-1}$ определяется по формуле $S^{2j-1}=\{(z_1,
\dots, z_j) \in \C^j \vert \vert z_1 \vert ^2+ \dots + \vert
z_j\vert ^2=1 \}$. Пусть $\i(z_1, \dots, z_j)=(\i z_1, \dots, \i z_j)$.
 Тем самым,
пространство $S^{2j-1}/\i$ служит $(2j-1)$-мерным остовом
пространства $S^{\infty}/i$, который  называется $(2j-1)$-мерным
линзовым пространством над $\Z/4$. Само пространство
$S^{\infty}/\i$ является пространством Эйленберга-Маклейна
$K(\I_a,1)$.

Определим пространство $K(\Q_a,1)$.
Рассмотрим
бесконечномерную сферу $S^{\infty}$ (это стягиваемое
пространство), которая определяется как прямой предел
бесконечной цепочки вложений стандартных сфер:
$$S^{\infty}=
\lim_{\longrightarrow} \quad (S^3 \subset S^7 \subset \dots \subset
S^{4j-1} \subset S^{4j+3} \subset \dots  \quad).$$

 Определено покоординатное действие $\Q_a
\times (\C^{2})^j \to (\C^{2})^j$ на каждом прямом слагаемом $\C^2$ в
соответствии с формулами $\ref{Q a1}, \ref{Q a2}, \ref{Q a3}$.
 Тем самым,
пространство $S^{4j-1}/\Q_a$ служит $(4j-1)$-мерным остовом
пространства $S^{\infty}/\Q_a$ и  называется
$(4j-1)$-мерным линзовым пространством над $\Q_a$. Само
пространство $S^{\infty}/\Q_a$ является пространством
Эйленберга-Маклейна $K(\Q_a,1)$.

\subsubsection*{Определение характеристического числа $h_{\mu_a,k}$}

Пусть в предположении $n>2k$ на многообразии точек самопересечения
$N^{n-2k}$ скошенно-оснащенного погружения определено произвольное
отображение $\mu_a: N^{n-2k} \to K(\I_a,1)$. Определим
характеристическое число $h_{\mu_a,k}$ по формуле:

\begin{eqnarray}\label{mu_a}
 h_{\mu_a,k}=\langle \bar \mu_a^{\ast}x;[\bar
N_a^{n-2k}]\rangle,
\end{eqnarray}
где $\bar \mu_a: \bar N_a^{n-2k} \to K(\I_d,1)$ -- двулистное
накрывающее над отображением $\mu_a: N_a^{n-2k} \to K(\I_a,1)$, индуцированное накрытием
$K(\I_d,1) \to K(\I_a,1)$, $x \in  H^{n-2k}(K(\I_d,1);\Z/2)$ --
образующая,  $[\bar N_a^{n-2k}]$ -- фундаментальный класс многообразия
$\bar N_a^{n-2k}$. (Многообразие $\bar N_a^{n-2k}$ совпадает с
каноническим 2-листным накрывающим $\bar N^{n-2k}$, если структурное отображение
$\eta: N^{n-2k} \to K(\D_4,1)$
является циклическим.)

\subsubsection*{Определение характеристического числа $h_{\lambda_a,k}$}

Пусть в предположении $n>4k$ на многообразии точек самопересечения
$L^{n-4k}$ $\D_4$--оснащенного погружения определено произвольное
отображение $\lambda_a: L^{n-4k} \to K(\Q_a,1)$. Определим
характеристическое число $h_{\lambda_a,k}$ по формуле:

\begin{eqnarray} \label{lambda_a}
 h_{\lambda_a,k}=\langle \bar
\lambda_a^{\ast}y;[\bar L_a]\rangle,
\end{eqnarray}
 где
$y \in H^{n-4k}(K(\H_e,1);\Z/2)$ -- образующая, $\bar \lambda_a:
\bar L_a^{n-4k} \to K(\H_e,1)$ -- 4-листное накрытие над
отображением $\lambda_a$, индуцированное из накрытия $K(\H_e,1)
\to K(\Q_a,1)$, $[\bar L_a]$-- фундаментальный класс многообразия
$\bar L_a^{n-4k}$. (Многообразие $\bar L_a^{n-4k}$ не совпадает с
каноническим 4-листным накрывающим $\bar L^{n-4k}$, если
структурное отображение $\zeta$ не является кватернионным.)
где $\bar \lambda_a: \bar L_a^{n-4k} \to K(\H_e,1)$ -- 4-листное
накрытие над отображением $\lambda_a: L_a^{n-4k} \to K(\Q_a,1)$, индуцированное накрытием
$K(\Q_e,1) \to K(\Q_a,1)$, $y \in  H^{n-4k}(K(\H_e,1);\Z/2)$ --
образующая,  $[\bar L_a^{n-4k}]$ -- фундаментальный класс многообразия
$\bar L_a^{n-4k}$. (Многообразие $\bar L_a^{n-4k}$ совпадает с
каноническим 4-листным накрывающим $L_{\H_{\bar c}}^{n-4k}$, если структурное отображение $\zeta: L^{n-4k}
\to K(\Z/2^{[3]},1)$
является кватернионным.)

\begin{lemma}\label{lemma12}
Для произвольного скошенно-оснащенного погружения $(f: M^{n-k}
\looparrowright \R^n, \kappa, \Xi)$ с многообразием
самопересечения $N^{n-2k}$, для которого классифицирующее
отображение $\eta$ нормального расслоения оказывается циклическим,
справедливо равенство
$$h_k(f,\kappa, \Xi) = h_{\mu_a,k},$$
где характеристическое число в правой части вычислено для
такого отображения $\mu_a$, что $\eta = i_a \circ \mu_a$,  $i_a: \I_a \subset \D_4$.
\end{lemma}

\subsubsection*{Доказательство Леммы $\ref{lemma12}$}

Рассмотрим двулистное накрытие $\bar \mu_a: \bar N_a^{n-2k} \to
K(\I_d,1)$ над отображением $\mu_a: N^{n-2k} \to K(\I_a,1)$,
индуцированное двулистным накрытием $K(\I_d,1) \to K(\I_a,1)$ над
пространством--образом отображения. Поскольку структурное
отображение $\eta$ является циклическим, многообразие $\bar
N^{n-2k}_a$ совпадает с каноническим двулистным накрывающим $\bar
N^{n-2k}$  над многообразием  самопересечения $N^{n-2k}$
погружения $f$. Доказательство леммы  вытекает из Леммы 7,
поскольку отображения $\bar \mu_a$ и $\bar \eta$ совпадают и
характеристическое число $h_{\mu_a,k}$ вычисляется как в левой
части равенства ($\ref{1.99}$).
\[  \]

\begin{lemma}\label{lemma13}
Для произвольного $\D_4$--оснащенного погружения $(g: N^{n-2k}
\looparrowright \R^n, \eta, \Psi)$ с многообразием самопересечения
$L^{n-4k}$, для которого классифицирующее отображение $\zeta$
нормального расслоения оказывается кватернионным, справедливо
равенство
$$h_k(g,\eta, \Psi) = h_{\lambda_a,k},$$
где характеристическое число в правой части вычислено для такого
отображения $\lambda_a$, что $\zeta = i_{a} \circ \lambda_a$, $i_{a}: \Q_a
\subset \Z/2^{[3]}$.
\end{lemma}

\subsubsection*{Доказательство Леммы $\ref{lemma13}$}

Рассмотрим 4-листное накрытие $\bar \lambda_a: \bar L_a^{n-4k} \to
K(\H_d,1)$ над отображением $\lambda_a: L^{n-4k} \to K(\Q_a,1)$,
индуцированное 4-листным накрытием $K(\Q_e,1) \to K(\Q_a,1)$ над
пространством--образом отображения. Поскольку структурное
отображение $\zeta$ является кватернионным, многообразие $\bar
L^{n-4k}_a$ совпадает с каноническим 4-листным накрывающим $\bar
L^{n-4k}_{sf}$ над многообразием точек самопересечения $L^{n-4k}$
погружения $g$. Доказательство леммы вытекает из Предложения 8,
поскольку отображения $\bar \lambda_a$ и $\bar \zeta^{sf}$
совпадают и характеристическое число $h_{\lambda_a,k}$ вычисляется
по отображению $\bar \lambda_a$ как в правой части равенства
($\ref{1.60}$).
\[  \]

\begin{definition}\label{def14}
Пусть $N^{n-2k}$ -- многообразие двукратных точек самопересечения
скошенно-оснащенного погружения   $(f,  \Xi, \kappa)$ в
предположении, что $f: M^{n-k} \looparrowright \R^n$ -- погружение
общего положения.
 Скажем, что это скошенно-оснащенное погружение
 допускает циклическую структуру, если существует
 отображение
$\mu_a : N^{n-2k} \to K(\I_a,1)$ (при этом, вообще говоря, классифицирующее отображение $\eta$
циклическим не предполагается) такое, что $h_{\mu_a,k}=h_k(f,\Xi,
\kappa)$.
\end{definition}

\begin{example}\label{def15}
Из Леммы 10 вытекает, что если классифицирующее отображение $\eta$
является циклическим, то  циклическую структуру можно задать
отображением  $\mu_a: N^{n-2k} \to K(\I_a,1)$, где $i_a \circ \mu_a =
\eta$, $i_a: K(\I_a,1) \subset K(\D_4,1)$.
\end{example}

\begin{definition}\label{def16}
Пусть $L^{n-4k}$--многообразие двукратных точек самопересечения
$\D_4$-оснащенного погружения $(g, \Psi, \eta)$ в предположении,
что $g: N^{n-2k} \looparrowright \R^n$--погружение общего
положения.
 Скажем, что это $\D_4$--оснащенное погружение
 допускает кватернионную структуру, если существует
 отображение
$\lambda_a : L^{n-4k} \to K(\Q_a,1)$ (при этом, вообще говоря, классифицирующее отображение $\zeta$
циклическим не предполагается) такое, что
$h_{\lambda_a,k}=h_k(g,\Psi,  \eta)$.
\end{definition}

\begin{example}\label{ex17}
Из Леммы 11 вытекает, что если структурное отображение $\zeta$
является кватернионным, то кватернионную структуру можно задать
отображением  $\lambda_a: L^{n-4k} \to K(\Q_a,1)$, где $i_{Q_a}
\circ \lambda_a = \zeta$, $i_{Q_a}: K(\Q_a,1) \subset K(\Z/2^{[3]},1)$.
\end{example}
\[  \]

Нам потребуется cформулировать понятие циклической и кватернионной
структур без предположения о том, что отображения $f: M^{n-k} \to
\R^n$, $g: N^{n-2k} \to \R^n$ являются погружениями. Мы
сформулируем необходимое определение в минимальной общности в
предположении $M^{n-k}=\RP^{n-k}$, $N^{n-2k} = S^{n-2k}/\i$. Итак,
пусть $d: \RP^{n-k} \to \R^n$ является гладким отображением общего
положения.


Рассмотрим конфигурационное пространство
\begin{eqnarray}\label{99}
(\RP^{n-k} \times
\RP^{n-k} \setminus \Delta_{\RP^{n-k}})/T',
\end{eqnarray}
которое также
называется "взрезанный квадрат" \- пространства $\RP^{n-k}$. Это
пространство получено путем факторизации прямого произведения без
диагонали по инволюции $T': \RP^{n-k} \times \RP^{n-k} \to
\RP^{n-k} \times \RP^{n-k}$, переставляющей координаты.
Построенное пространство является открытым многообразием.

Определим пространство $\bar
\Gamma_0$ как сферическое раздутие пространства $\RP^{n-k} \times \RP^{n-k}
\setminus \Sigma_{diag}$ в окрестности диагонали. Сферическим раздутием называется
многообразие с краем, которое определено в результате компактификации открытого многообразия
$\RP^{n-k} \times \RP^{n-k}
\setminus \Sigma_{diag}$ посредством послойного вклеивания слоев сферизации $ST \Sigma_{diag}$
касательного расслоения $T \Sigma_{diag}$ в окрестности нулей слоев нормального расслоения
диагонали $\Sigma_{diag} \subset \RP^{n-k} \times \RP^{n-k}$. Определены естественные включения:
$$\RP^{n-k} \times \RP^{n-k} \setminus \Sigma_{diag} \subset \bar \Gamma_0,$$
$$ST \Sigma_{diag} \subset \bar \Gamma_0.$$
На пространстве $\bar \Gamma_0$ определена свободная инволюция
$\bar T': \bar \Gamma_0 \to \bar \Gamma_0$, которая определена как продолжение инволюции $T'$.
Факторпространство
$\bar \Gamma_0 / \bar T'$ обозначается через $\Gamma_0$, а соответствующее двулистное накрытие через
$$
p_{\Gamma_0}:  \bar \Gamma_0 / \bar T' \to \Gamma_0.
$$
Пространство $\Gamma_0$ является многообразием с краем и оно называется пространством раздутия
конфигурационного пространства $(\ref{99})$. Определена проекция $p_{\partial \Gamma_0}: \partial \Gamma_0 \to \RP^{n-k}$, это отображение отображение называется разрешением диагонали.

Для произвольного отображения $d$ определено пространство $N(d)$ точек самопересечения отображения $d$ по формуле:

\begin{eqnarray}\label{Nd}
N(d) = Cl \{ ([x,y]) \in int(\Gamma_0) : y \ne x, d(y)=d(x) \}.
\end{eqnarray}

По теореме Портеуса [Por] в предположении о том, что отображение $d$ общего положения, пространство $N(d)$ является многообразием с краем
размерности $n-2k$. Это многообразие обозначим через $N^{n-2k}(d)$ и назовем многообразием точек самопересечения отображения $d$.
При этом формула $(\ref{Nd})$ определяет вложение многообразий с краем:
$$i_{N(d)}: (N^{n-2k}(d), \partial N^{n-2k}(d)) \subset (\Gamma_0,\partial \Gamma_0). $$

 Край $\partial N^{n-2k}(d)$ многообразия $N^{n-2k}(d)$ называется многообразием разрешений критических
 точек отображения $d$. Отображение $p_{\partial \Gamma_0} \circ  i_{\partial N(d)} \vert_{\partial N(d)}:
 \partial N^{n-2k}(d) \subset \partial \bar \Gamma_0  \to \RP^{n-k}$ называется отображением
 разрешения особенностей отображения $d$, обозначим это отображение через $res_{d}: \partial N(d) \to \RP^{n-k}$.
 Определено
каноническое двулистное накрытие
\begin{eqnarray}\label{pNd}
p_{N(d)}: \bar N(d)^{n-2k} \to N(d)^{n-2k},
\end{eqnarray}
разветвленное над краем $\partial N(d)^{n-2k}$ (над этим краем
накрытие является диффеоморфизмом). При этом следующая диаграмма
является коммутативной:

$$
\begin{array}{ccc}
i_{\bar N(d)}: (\bar N^{n-2k}(d), \partial N^{n-2k}(d)) & \subset & (\bar \Gamma_0, \partial \Gamma_0) \\
\\
\downarrow p_{N(d)}&   &  \downarrow  p_{\Gamma_0} \\
\\
i_{N(d)}: (N^{n-2k}(d), \partial N^{n-2k}(d)) & \subset & (\Gamma_0, \partial \Gamma_0).
\end{array}
$$

\subsubsection*{Определение структурного отображения $\eta_N: N^{n-k}(d) \to
K(\D_4,1)$}

Определим отображение $\eta_{\Gamma_0}: \Gamma_0 \to K(\D_4,1)$,
которое назовем структурным отображением "взрезанного квадрата".
Заметим, что включение $\bar \Gamma_0 \subset \RP^{n-k} \times
\RP^{n-k}$ индуцирует изоморфизм фундаментальных групп, т.к.
коразмерность диагонали $\Delta_{\RP^{n-k}} \subset \RP^{n-k}
\times \RP^{n-k}$ равна $n-k \ge 3$. Следовательно, справедливо равенство

\begin{eqnarray}\label{pi_1}
\pi_1(\bar \Gamma_0) = H_1(\bar \Gamma_0;\Z/2)= \Z/2 \oplus \Z/2.
\end{eqnarray}

Рассмотрим индуцированный автоморфизм $T'_{\ast}: H_1(\bar
\Gamma_0;\Z/2) \to H_1(\bar \Gamma_0;\Z/2)$. Заметим, что этот
автоморфизм не является тождественным. Зафиксируем изоморфизм
групп $H_1(\bar \Gamma_0;\Z/2)$ и $\I_c$, при котором образующая
первого (соответственно, второго) слагаемого группы $H_1(\bar \Gamma_0;\Z/2)$ (см.
(\ref{pi_1})) переходит в образующую $ab \in \I_c \subset \D_4$
(соответственно, $ba \in \I_c \subset \D_4$), которая в стандартном
представлении $\D_4$ определена преобразованием симметрии относительно второй
(соответственно, первой) координатной оси.

Легко проверить, что автоморфизм внешнего сопряжения подгруппы
$\I_c \subset \D_4$ на элемент $b \in \D_4 \setminus \I_c$ (в этой
формуле элемент $b$ можно выбрать произвольным), определенный
формулой $x \mapsto bxb^{-1}$ переходит при построенном
изоморфизме в автоморфизм $T'_{\ast}$. Фундаментальная группа $\pi_1(\Gamma_0)$
является квадратичным расширением группы $\pi_1(\bar \Gamma_0)$
посредством элемента $b$, и это расширение однозначно c точностью
до изоморфизма определяется автоморфизмом $T'_{\ast}$. Поэтому
 $\pi_1(\Gamma_0) \simeq \D_4$ и, значит, определено отображение
$\eta_{\Gamma_0}: \Gamma_0 \to K(\D_4,1)$.

Нетрудно проверить, что отображение $\eta_{\Gamma_0} \vert_
{\partial \Gamma_0}$ принимает значения в подпространстве
$K(\I_b,1) \subset K(\D_4,1)$. Отображение $\eta_{\Gamma_0}$
индуцирует отображение $\eta_N: (N(d),
\partial N(d)) \to (K(\D_4,1),K(\I_b,1))$, которое назовем структурным отображением.
Также нетрудно проверить, что гомотопический класс композиции
$\partial N(d) \stackrel{\eta}{\longrightarrow} K(\I_b,1)
\stackrel{p_b}{\longrightarrow} K(\I_d,1)$ совпадает с
характеристическим отображением $\kappa \circ res_{d}: \partial
N(d) \to \RP^{n-k} \to K(\I_d,1)$, которое является композиции
отображения разрешения особенностей $res_d: \partial N(d)^{n-2k}
\to  \RP^{n-k}$ и вложения остова $\RP^{n-k} \subset K(\I_d,1)$ в
классифицирующее пространство.

\subsubsection*{Определение циклической структуры для
отображения $d: \RP^{n-k} \to \R^n$ с особенностью}

Пусть $d: \RP^{n-k} \to \R^n$-- кусочно-линейное отображение
общего положения, вообще говоря, имеющее критические точки, $k
\equiv 0 \pmod{2}$. Пусть многообразие $N^{n-2k}(d)$, вообще
говоря, с непустым краем $\partial N^{n-2k}(d)$, служит
многообразием точек самопересечения отображения $d$. Отображение
$\mu_{a,N(d)} : (N^{n-2k}(d),\partial N^{n-2k}(d)) \to
(K(\I_a,1),K(\I_d,1))$ называется циклической структурой
отображения $d$, если $\mu_{a,N(d)}$ удовлетворяет краевому
условию:

\begin{eqnarray}
 \mu_{a,N(d)} \vert_{\partial N^{n-2k}(d)}= (i_a \circ \kappa)
\vert_{\partial N(d)}, \label{mu_a boundary}
\end{eqnarray}
где $\kappa: \RP^{n-k} \to K(\I_d,1)$ -- классифицирующее отображение для образующего класса
когомологий, $i_a: K(\I_d,1) \subset K(\I_a,1)$ -- отображение классифицирующих пространств, индуцированное включением подгрупп $i_{d,a}: \I_d \subset \I_a$, и, кроме
того, выполнено условие
\begin{eqnarray}
\langle \mu_{a,N(d)}^{\ast}(t);[N^{n-2k}(d),\partial N^{n-2k}(d)]
\rangle =1, \label{h mu_a}
\end{eqnarray}
 где  $t \in
H^{n-2k}(K(\I_a,1),K(\I_d,1);\Z/2)$ -- произвольный класс
когомологий, который переходит в образующий класс когомологий в
группе $H^{n-2k}(K(\I_a,1);\Z/2)$ при индуцированном гомоморфизме
$j^{\ast}: H^{n-2k}(K(\I_a,1),K(\I_d,1);\Z/2) \to
H^{n-2k}(K(\I_a,1);\Z/2)$, $[N^{n-2k}(d),\partial N^{n-2k}(d)]$ --
относительный фундаментальный цикл многообразия.
\[  \]

Пусть $c: S^{n-2k}/\i \to \R^n$ является
 гладким отображением общего положения.


Рассмотрим конфигурационное пространство
\begin{eqnarray}\label{99.1}
(S^{n-2k}/\i \times
S^{n-2k}/\i \setminus \Delta_{S^{n-2k}})/T',
\end{eqnarray}
которое также
называется "взрезанный квадрат" \- линзового пространства $S^{n-2k}/\i$. Это
пространство получено путем факторизации прямого произведения без
диагонали по инволюции $T': S^{n-2k}/\i \times S^{n-2k}/\i \to
S^{n-2k}/\i \times S^{n-2k}/\i$, переставляющей координаты.
Построенное пространство является открытым многообразием.

Определим пространство $\bar
\Gamma_1$ как сферическое раздутие пространства $S^{n-2k}/\i \times
S^{n-2k}/\i \setminus \Delta_{S^{n-2k}}$ в окрестности диагонали. Сферическим раздутием называется
многообразие с краем, которое определено в результате компактификации открытого многообразия
$S^{n-2k}/\i \times S^{n-2k}/\i
\setminus \Sigma_{diag}$ посредством послойного вклеивания слоев сферизации $ST \Sigma_{diag}$
касательного расслоения $T \Sigma_{diag}$ в окрестности нулей слоев нормального расслоения
диагонали $\Sigma_{diag} \subset S^{n-2k} \times S^{n-2k}$. Определены естественные включения:
$$S^{n-2k} \times S^{n-2k} \setminus \Sigma_{diag} \subset \bar \Gamma_1,$$
$$ST \Sigma_{diag} \subset \bar \Gamma_1.$$
На пространстве $\bar \Gamma_1$ определена свободная инволюция
$\bar T': \bar \Gamma_1 \to \bar \Gamma_1$, которая определена как продолжение инволюции $T'$.
Факторпространство
$\bar \Gamma_1 / \bar T'$ обозначается через $\Gamma_1$, а соответствующее двулистное накрытие через
$$
p_{\Gamma_1}:  \bar \Gamma_1 / \bar T' \to \Gamma_1.
$$
Пространство $\Gamma_1$ является многообразием с краем и оно называется пространством раздутия
конфигурационного пространства $(\ref{99.1})$. Определена проекция $p_{\partial \Gamma_1}: \partial \Gamma_1 \to S^{n-2k}/\i$, это отображение отображение называется разрешением диагонали.

Для произвольного отображения $c$ определено пространство $L(c)$ точек самопересечения отображения $c$ по формуле:

\begin{eqnarray}\label{Nc}
L(c) = Cl \{ ([x,y]) \in int(\Gamma_1) : y \ne x, c(y)=c(x) \}.
\end{eqnarray}

По теореме Портеуса [Por] в предположении о том, что отображение $c$ общего положения, пространство $L(c)$ является многообразием с краем
размерности $n-4k$. Это многообразие обозначим через $L^{n-4k}(c)$ и назовем многообразием точек самопересечения отображения $c$.
При этом формула $(\ref{Nc})$ определяет вложение многообразий с краем:
$$i_{L(c)}: (L^{n-4k}(c), \partial L^{n-4k}(c)) \subset (\Gamma_1,\partial \Gamma_1). $$

 Край $\partial L^{n-4k}(c)$ многообразия $L^{n-4k}(c)$ называется многообразием разрешений критических точек отображения $c$. Отображение $p_{\partial \Gamma_1} \circ  i_{\partial L(c)} \vert_{\partial L(c)}:
 \partial L^{n-4k}(c) \subset \partial \bar \Gamma_1  \to S^{n-2k}/\i$ называется отображением
 разрешения особенностей отображения $c$, обозначим это отображение через $res_{L(c)}: \partial L(c) \to S^{n-2k}/\i$.
 Определено
каноническое двулистное накрытие
\begin{eqnarray}\label{pLc}
p_{L(c)}: \bar L(c)^{n-4k} \to L(c)^{n-4k},
\end{eqnarray}
разветвленное над краем $\partial L(c)^{n-4k}$; над этим краем
накрытие является диффеоморфизмом. При этом следующая диаграмма
является коммутативной:

$$
\begin{array}{ccc}
i_{\bar L(c)}: (\bar L^{n-4k}(c), \partial L^{n-4k}(c)) & \subset & (\bar \Gamma_1, \partial \Gamma_1) \\
\\
\downarrow p_{L(c)}&   &  \downarrow  p_{\Gamma_1} \\
\\
i_{L(c)}: (L^{n-4k}(c), \partial L^{n-4k}(c)) & \subset & (\Gamma_1, \partial \Gamma_1).
\end{array}
$$

\subsubsection*{Определение подгруппы $\H \subset \Z/2^{[3]}$}
Рассмотрим пространство $\R^4$ c базисом $(\e_1,\e_2,\e_3,\e_4)$.
Вектора базиса удобно отождествить с базисными единичными
кватернионами $({\bf 1}, \i,\j,\k)$, что иногда будет
использоваться, т.к. позволит упростить формулы некоторых
преобразований. Определим подгруппу
\begin{eqnarray}\label{H}
\H \subset \Z/2^{[3]}
\end{eqnarray}
как подгруппу преобразований, преобразований следующих двух типов:

-- в каждой плоскости $(\e_1 ={\bf 1},\e_2=\i)$,
$(\e_3=\j,\e_4=\k)$ допускается (взаимно независимые)
преобразования, кратные умножению на кватернион $\i$. Подгруппу
всех таких преобразований обозначим через $\H_c$, это подгруппа
изоморфна $\Z/4 \oplus \Z/4$.

-- преобразование, переставляющее между собой одновременно пару
базисных векторов $\e_1 ={\bf 1}$ и $\e_3 =\j$ и
 пару базисных векторов  $\e_2 =\i$ и $\e_4 =\k$, сохраняя их направления. Обозначим это преобразование через
 $t \in \H_c$.

 Легко проверить, что сама группа $\H$ имеет порядок $32$ и является подгруппой $\H \subset \Z/2^{[3]}$
индекса 4.
\[  \]

\subsubsection*{Определение подгруппы $\H_b \subset \H$ и мономорфизма $i_{\I_a,\H}: \I_a \subset \H$.}
Определим гомоморфизм включения $i_{\I_a,\H}: \I_a \subset \H$,
который переводит образующую группы $\I_a$ в преобразование
умножения на кватернион $\i$, действующее одновременно в каждой
плоскости $(\e_1 ={\bf 1},\e_2=\i)$, $(\e_3=\j,\e_4=\k)$.
Определим подгруппу $\H_b \subset \H$ как сумму подгруппы
$i_{\I_a,\H}: \I_a \subset \H$ и подгруппы, порожденной образующей
$t \in \H$. Легко проверить, что сама группа $\H_b$ имеет порядок
$16$ и изоморфна группе $\Z/4 \oplus \Z/2$. Подгруппа $\H_b
\subset \H$ имеет индекс 2. Определены также гомоморфизм включения
$i_{\I_a,\H_b}: \I_a \subset \H_b$ подгрупп и гомоморфизм проекции
$p_{\H_b,\I_a}: \H_b \to \I_a$, такие, что гомоморфизм композиции
$\I_a \stackrel{i_{\I_a,\H_b}}{\longrightarrow} \H_b
\stackrel{p_{\H_a,\I_a}}{\longrightarrow} \I_a$ является
тождественным.
\[  \]

\subsubsection*{Определение структурного отображения $\zeta: L^{n-4k}(c) \to
K(\H,1)$}

Определим отображение $\zeta_{\Gamma_1}: \Gamma_1 \to K(\H,1)$,
которое назовем структурным отображением "взрезанного квадрата".
Заметим, что включение $\bar \Gamma_1 \subset S^{n-2k}/\i \times
S^{n-2k}/\i$ индуцирует изоморфизм фундаментальных групп, т.к.
коразмерность диагонали $\Delta_{S^{n-2k}/\i} \subset S^{n-2k}/\i
\times S^{n-2k}/\i$ равна $n-2k \ge 3$. Справедливо равенство:

\begin{eqnarray}\label{pi_1H}
\pi_1(\bar \Gamma_1) = H_1(\bar \Gamma_1;\Z/4)= \Z/4 \oplus \Z/4.
\end{eqnarray}

Рассмотрим индуцированный автоморфизм $T'_{\ast}: H_1(\bar
\Gamma_1;\Z/4) \to H_1(\bar \Gamma_1;\Z/4)$. Заметим, что этот
автоморфизм не является тождественным. Зафиксируем мономорфизм
групп $H_1(\bar \Gamma_1;\Z/4)$ в группу $\H_c$, при котором
образующая первого (соответственно, второго) слагаемого группы
$H_1(\bar \Gamma_1;\Z/2)$ (см. (\ref{pi_1H})) переходит в
образующую, которая в стандартном представлении подгруппы $\H_c
\subset \H \subset \Z/2^{[3]}$ определена умножением на кватернион
$\i$ в двумерной плоскости $({\bf 1},\i)$ (соответственно, в
плоскости $(\j, \k)$) и неподвижно в дополнительной плоскости.

Легко проверить, что автоморфизм внешнего сопряжения подгруппы
$\H_c \subset \H$ на произвольный элемент $b \in \H \setminus \H_c$, определенный
формулой $x \mapsto bxb^{-1}$, переходит при построенном
изоморфизме $\ref{pi_1H}$ в автоморфизм $T'_{\ast}$. Фундаментальная группа $\pi_1(\Gamma_1)$
является квадратичным расширением группы $\pi_1(\bar \Gamma_1)$
посредством элемента $b$, и это расширение однозначно c точностью
до изоморфизма определяется автоморфизмом $T'_{\ast}$. Поэтому
 $\pi_1(\Gamma_1) \simeq \H$ и, значит, определено отображение
$\zeta_{\Gamma_1}: \Gamma_1 \to K(\H,1)$.

Нетрудно проверить, что отображение $\zeta_{\Gamma_1} \vert_
{\partial \Gamma_1}$ принимает значения в подпространстве
$K(\H_b,1) \subset K(\H,1)$. Отображение $\zeta_{\Gamma_1}$
индуцирует отображение $\zeta: (L^{n-4k}(c),
\partial L^{n-4k}(c)) \to (K(\H,1),K(\H_b,1))$, которое назовем структурным отображением.
Также нетрудно
проверить, что гомотопический класс композиции $\partial L^{n-4k}(c)
\stackrel{\zeta}{\longrightarrow} K(\H_b,1)
\stackrel{p_{\H_b,\I_d}}{\longrightarrow} K(\I_a,1)$ совпадает с отображением
$\eta \circ res_c : \partial L^{n-4k}(c) \to S^{n-2k}/\i \to
K(\I_a,1)$, которое определяется как
композиция отображения разрешения особенностей $res_c: \partial L(c)^{n-2k} \to S^{n-2k}/\i$ и вложения остова $\eta: S^{n-2k}/i
\subset K(\I_a,1)$ в классифицирующее пространство.

\subsubsection*{Определение кватернионной структуры для
отображения $c: S^{n-2k}/\i \to \R^n$ с особенностью}

Пусть $c: S^{n-2k}/\i \to \R^n$-- кусочно-линейное отображение
общего положения, вообще говоря, имеющее критические точки, $k
\equiv 0 \pmod{2}$. Пусть многообразие $L^{n-4k}(c)$, вообще
говоря, с непустым краем $\partial L^{n-4k}(c)$ является
многообразием точек самопересечения отображения $c$. Отображение
$\lambda_{a,L(c)} : (L^{n-4k}(c),\partial L^{n-4k}(c)) \to
(K(\H,1),K(\H_b,1))$ называется циклической структурой отображения
$c$, если $\lambda_{a,L(c)}$ удовлетворяет краевому условию:

\begin{eqnarray}\label{lambda_a boundary}
 \lambda_{a,L(c)} \vert_{\partial L^{n-4k}(c)}= (i_{\H_b,\I_a} \circ \lambda_{a,L(c)})
\vert_{\partial L^{n-4k}(c)},
\end{eqnarray}
где $\eta: S^{n-2k} \to K(\I_a,1)$ -- характеристическое отображение образующего класса
когомологий (включение остова в классифицирующее пространство), $i_{\I_a,\H_b}: K(\I_a,1) \subset K(\H_b,1)$--
отображение классифицирующих пространств, индуцированное гомоморфизмом включения $i_{\I_a,\H_b}: \I_a \subset \H_b$
и, кроме
того, выполнено условие:
\begin{eqnarray}
\langle \lambda_{a,L(c)}^{\ast}(t);[L^{n-4k}(c),\partial
L^{n-4k}(c)] \rangle =1, \label{h lambda a}
\end{eqnarray}
 где  $t \in
H^{n-4k}(K(\Q_a,1),K(\I_a,1);\Z/2)$ -- произвольный класс
когомологий, который переходит в образующий класс когомологий в
группе $H^{n-4k}(K(\Q_a,1);\Z/2)$ при индуцированном гомоморфизме
$j^{\ast}: H^{n-4k}(K(\Q_a,1),K(\I_a,1);\Z/2) \to
H^{n-4k}(K(\Q_a,1);\Z/2)$, $[L^{n-4k}(c),\partial L^{n-4k}(c)]$ --
относительный фундаментальный цикл многообразия.
\[  \]

Следующие две леммы необходимы для проверки корректности
определений циклической и кватернионной структуры.

\begin{lemma}\label{lemma18}
Характеристическое число, определенное левой частью формулы
$(\ref{h mu_a})$, не зависит от выбора класса
когомологий $t \in H^{n-2k}(K(\I_a,1),K(\I_d,1);\Z/2)$, для
которого $j^{\ast}(t) \in H^{n-2k}(K(\I_a,1);\Z/2)$ является
образующим классом.
\end{lemma}

\begin{lemma}\label{lemma19}
Характеристическое число, определенное левой частью формулы
$(\ref{h lambda a})$, не зависит от выбора класса когомологий $t
\in H^{n-4k}(K(\Q_a,1),K(\I_a,1);\Z/2)$, для которого $j^{\ast}(t)
\in H^{n-4k}(K(\Q_a,1);\Z/2)$ является образующим классом.
\end{lemma}

\subsubsection*{Доказательство Леммы $\ref{lemma18}$}

Действительно, при другом выборе $t'$ получим, что $t-t' \in
\mathrm{Im} \delta^{\ast}: H^{n-2k-1}(K(\I_d,1);\Z/2) \to
H^{n-2k}(K(\I_a,1),K(\I_d,1);\Z/2))$. Следовательно, класс
$\mu_{a,N(d)}^{\ast}(t)-\mu_{a,N(d)}^{\ast}(t') \in
H^{n-2k}(N^{n-2k}(d),\partial N^{n-2k}(d);\Z/2))$ получен в
результате применения кограничного гомоморфизма $\delta^{\ast}:
H^{n-2k-1}(\partial N^{n-2k}(d);\Z/2) \to H^{n-2k}(N^{n-2k}(d),
\partial N^{n-2k}(d);\Z/2)$ к некоторому элементу $x = \mu_a^{\ast}(q)$,
где $q \in H^{n-2k-1}(K(\I_d,1);\Z/2)$ -- образующий класс
когомологий. Теперь заметим, что гомоморфизм
$$(\mu_{a,N(d)} \vert_{\partial N(d)})^{\ast}:
H^{n-2k-1}(K(\I_d,1);\Z/2) \to H^{n-2k-1}(\partial
N^{n-2k}(d);\Z/2)$$ является тривиальным. Для этого докажем, что
подмногообразие $\partial N^{n-2k}(d) \subset \RP^{n-k}$
представляет нулевой гомологический класс. Исходя из
геометрического смысла характеристических классов Штифеля-Уитни
(см. [M-S], параграф 12), класс $[\partial N^{n-2k}(d)] \in
H_{n-2k-1}(\RP^{n-k};\Z/2)$ двойственнен в смысле Пуанкаре
характеристическому классу $\bar w_{k+1}(\RP^{n-k})$ Штифеля-Уитни
нормального расслоения над $\RP^{n-k}$, поскольку является классом
гомологий, препятствующим к построению редукции стабильного
нормального расслоения этого многообразия до нормального
расслоения размерности $k$. Вспомним, что $n=2^l-1$, поэтому самый
старший характеристический класс $\bar w_{\ast}(\RP^{n-k})$
(нормального расслоения) имеет размерность $k$ и $\bar
w_{k+1}(\RP^{n-k})=0$ (Например, в простейшем случае $n=4$, $k=1$
получим $\bar w_2(\RP^2)=0$.) Лемма $\ref{lemma18}$ доказана.
\[  \]

Доказательство Леммы $\ref{lemma19}$ аналогично и опускается.
\[  \]

Пусть $d: \RP^{n-k} \to \R^n$ -- кусочно-линейное отображение
общего положения, $(N^{n-2k}(d), \partial N^{n-2k}(d))$ -- многообразие
с краем двойных точек самопересечения отображения $d$.
Пусть заданно
 отображение $\mu_{a,N(d)}: (N^{n-2k}(d),\partial N^{n-2k}(d)) \to
 (K(\I_a,1),K(\I_d,1))$, удовлетворяющее граничному условию ($\ref{mu_a boundary}$).
 Нам потребуется критерий проверки того, что отображение $\mu_{a,N(d)}$
  удовлетворяет условию ($\ref{h mu_a}$).

Рассмотрим двулистное накрытие
\begin{eqnarray}\label{p_a}
p_a: \bar N_a^{n-2k}(d) \to N^{n-2k}(d),
\end{eqnarray}
индуцированное из
универсального двулистного накрытия
\begin{eqnarray}\label{da}
K(\I_d,1) \to K(\I_a,1)
\end{eqnarray}
классифицирующих пространств отображением $\mu_{a,N(d)}:
N^{n-2k}(d) \to K(\I_a,1)$. Поскольку выполнено граничное условие
$(\ref{mu_a boundary})$, ограничение накрытия $(\ref{p_a})$ на
край $\partial N^{n-2k}(d)$ оказывается тривиальным двулистным
накрытием. Заметим, что в случае, когда структурное отображение
$\eta: N^{n-2k}(d) \to K(\D_4,1)$ является циклическим, накрытие
$(\ref{p_a})$ совпадает с каноническим двулистным накрытием,
определенным формулой $(\ref{pNd})$.

Рассмотрим отображение $\bar \mu_{a,N(d)}: \bar N^{n-2k}(d)_a \to
K(\I_d,1)$ накрывающих пространств, которое двулистно накрывает
отображение $\mu_{a,N(d)}: N^{n-2k}(d) \to K(\I_a,1)$ при
двулистных накрытиях $(\ref{p_a})$, $(\ref{da})$ над прообразом и
образом отображения $\mu_{a,N(d)}$.
 Определим замкнутое
многообразие $\bar R^{n-2k}$ как результат подклейки к
многообразию $\bar N^{n-2k}(d)_a$ цилиндра $C^{n-2k}$ двулистного
накрытия $\partial \bar N^{n-2k}(d)_a \to
\partial N^{n-2k}(d)$ вдоль общего края $\partial \bar N^{n-2k}(d)_a = \partial
C^{n-2k}$. Определим отображение $\mu_{\bar R}: \bar R^{n-2k} \to
K(\I_d,1)$ как результат склейки отображения $\bar \mu_a: \bar
N^{n-2k}(d)_a \to K(\I_d,1)$ c отображением
$$C^{n-2k}
\stackrel{p_C}{\longrightarrow}
\partial \bar N^{n-2k}(d)_a \stackrel{\mu_a \vert_{\partial \bar
N^{n-2k}(d)_a}}{\longrightarrow} K(\I_d,1),$$
где $p_C: C^{n-2k} \to
\partial N^{n-2k}(d)_a$ -- проекция цилиндра на основание.
Фундаментальный цикл многообразия $\bar R^{n-2k}$ обозначим через
$[\bar R]$.

\begin{lemma}\label{lemma20}

Для отображения $\mu_{a,N(d)}: (N^{n-2k}(d), \partial N^{n-2k}(d))
\to (K(\I_a,1),K(\I_d,1))$ равенство $(\ref{h mu_a})$ эквивалентно
тому, класс гомологий
$$ \mu_{\bar R,\ast}([\bar R^{n-2k}]) \in  H_{n-2k}(K(\I_d,1);\Z/2)$$ является образующим.
\end{lemma}

\subsubsection*{Доказательство Леммы $\ref{lemma20}$}
Из граничного условия ($\ref{mu_a boundary}$) вытекает, что класс
гомологий $\mu_{a,N(d);\ast}([\partial N^{n-2k}(d)])$ в группе
$H_{n-2k-1}(K(\I_d,1);\Z/2)$ равен нулю (см. эквивалентное
утверждение для классов когомологий в Лемме $\ref{lemma18}$).
Пусть $W^{n-2k}$ многообразие с краем $\partial W^{n-2k} =
\partial N^{n-2k}(d)$. Обозначим через
\begin{eqnarray}\label{varphi}
\varphi: W^{n-2k} \to
K(\I_d,1)
\end{eqnarray}
сингулярный относительный цикл, осуществляющий гомологию нулю
сингулярного цикла $\varphi \vert_{\partial W^{n-2k}} =
\mu_{a,N(d)} \vert_{\partial N^{n-2k}(d)}$ в подпространстве
$K(\I_d,1) \subset K(\I_a,1)$.

Образ фундаментального класса при отображении $(\mu_{a,N(d)}
\cup_{\partial} \varphi) : N^{n-2k}(d) \cup_{\partial} W^{n-2k}
\to K(\I_a,1)$ обозначим через $x \in H_{n-2k}(K(\I_a,1);\Z/2)$.
Очевидно, что при любом выборе относительного цикла
$(\ref{varphi})$ условие ($\ref{h mu_a}$), переформулированное для
соответствующих абсолютных классов гомологий, состоит в том, что
класс гомологий $x \in H_{n-2k}(K(\I_a,1);\Z/2)$ является
образующим.

Рассмотрим двулистное накрытие $\bar \mu_{a,N(d)} \cup_{\partial}
\bar \varphi: \bar N^{n-2k}(d)_a \cup_{\partial} \bar W^{n-2k} \to
K(\I_d,1)$ над отображением $\mu_{a,N(d)} \cup_{\partial}
\varphi$, индуцированное из накрытия $K(\I_d,1)  \to K(\I_a,1)$
над пространством-образом отображения $\mu_{a,N(d)}
\cup_{\partial} \varphi$. Образ фундаментального класса при этом
отображении обозначим через $\bar x \in H_{n-2k}(K(\I_d,1);\Z/2)$.
Элемент $\bar x$ является образующей, тогда и только тогда, когда
$x$ является образующей.

С другой стороны, класс гомологий $\bar x$ представлен образом
фундаментального класса при отображении $\mu_{\bar R}: \bar
R^{n-2k} \to K(\I_d,1)$. Лемма $\ref{lemma20}$ доказана.
\[  \]

Пусть $c: S^{n-2k}/\i \to \R^n$ -- кусочно-линейное отображение
общего положения, $(L^{n-4k}(c), \partial L^{n-4k}(c))$ --
многообразие с краем двойных точек самопересечения отображения
$c$. Пусть заданно
 отображение $\lambda_{a,L(c)}: (L^{n-4k}(c),\partial L^{n-4k}(c)) \to
 (K(\Q_a,1),K(\I_a,1))$, удовлетворяющее граничному условию $(\ref{lambda_a boundary})$.
 Нам потребуется критерий проверки того, что отображение $\lambda_{a,L(c)}$
  удовлетворяет уравнению ($\ref{h lambda a}$).

Рассмотрим двулистное накрытие
\begin{eqnarray}\label{p_aa}
p_a: \bar L_a^{n-4k}(c) \to L^{n-4k}(c),
\end{eqnarray}
индуцированное из
универсального двулистного накрытия
\begin{eqnarray}\label{daa}
K(\I_a,1) \to K(\Q_a,1)
\end{eqnarray}
классифицирующих пространств отображением $\lambda_{a,L(c)}:
L^{n-4k}(c) \to K(\Q_a,1)$. Поскольку выполнено граничное условие
$(\ref{mu_a boundary})$, ограничение накрытия $(\ref{p_aa})$ на
край $\partial L^{n-4k}(c)$ оказывается тривиальным двулистным
накрытием. Заметим, что в случае, когда структурное отображение
$\zeta: L^{n-4k}(c) \to K(\Z/2^{[3]},1)$ является циклическим,
накрытие $(\ref{p_aa})$ совпадает с каноническим двулистным
накрытием, определенным формулой $(\ref{pLc})$.

Рассмотрим отображение $\bar \lambda_{a,L(c)}: \bar L^{n-4k}(c)_a
\to K(\I_a,1)$ накрывающих пространств, которое двулистно
накрывает отображение $\lambda_{a,L(c)}: L^{n-4k}(c) \to
K(\Q_a,1)$ при двулистных накрытиях $(\ref{p_aa})$, $(\ref{daa})$
над прообразом и образом отображения $\lambda_{a,L(c)}$.
 Определим замкнутое
многообразие $\bar R^{n-4k}$ как результат подклейки к
многообразию $\bar L^{n-4k}(c)_a$ цилиндра $C^{n-4k}$ двулистного
накрытия $\partial \bar L^{n-4k}(c)_a \to
\partial L^{n-4k}(c)$ вдоль общего края $\partial \bar L^{n-4k}(c)_a = \partial
C^{n-4k}$. Определим отображение $\lambda_{\bar R}: \bar R^{n-4k}
\to K(\I_d,1)$ как результат склейки отображения $\bar
\lambda_{a,L(c)}: \bar L^{n-4k}(c)_a \to K(\I_a,1)$ c отображением
$$C^{n-4k}
\stackrel{p_C}{\longrightarrow}
\partial \bar L^{n-4k}(c)_a \stackrel{\lambda_a \vert_{\partial \bar
L^{n-4k}(c)_a}}{\longrightarrow} K(\I_a,1),$$
где $p_C: C^{n-4k} \to
\partial L^{n-4k}(c)_a$ -- проекция цилиндра на основание.
Фундаментальный цикл многообразия $\bar R^{n-4k}$ обозначим через
$[\bar R]$.

\begin{lemma}\label{lemma21}

Для отображения $\lambda_{a,L(c)}: (L^{n-4k}(c), \partial
L^{n-4k}(c)) \to (K(\Q_a,1),K(\I_a,1))$ равенство ($\ref{h lambda
a}$) эквивалентно тому, класс гомологий
$$ \lambda_{\bar R,\ast}([\bar R]) \in  H_{n-4k}(K(\Q_a,1);\Z/2)$$ является образующим.
\end{lemma}

Доказательство Леммы $\ref{lemma21}$ аналогично доказательству
Леммы $\ref{lemma20}$ и опускается.

\begin{proposition}\label{prop22}
При $k \equiv 0 \pmod{2}$, $k \ge 8$,  произвольный элемент группы
$Imm^{sf}(n-k,k)$ представлен скошенно--оснащенным погружением
$(f,\kappa, \Xi)$, допускающим циклическую структуру в смысле
Определения $\ref{defcycl}$.
\end{proposition}
\[  \]

\begin{proposition}\label{prop23}
При $k \equiv 0 \pmod{2}$, $10k-1>n$, $k \ge 8$,  произвольный
элемент группы $Imm^{\D_4}(n-2k,2k)$, лежащий в образе
гомоморфизма $\delta: Imm^{sf}(n-k,k) \to Imm^{\D_4}(n-2k,2k)$,
представлен $\D_4$--оснащенным погружением $(g,\eta, \Psi)$,
допускающим кватернионную структуру в смысле Определения
$\ref{defQ}$.
\end{proposition}
\[  \]

Доказательство Предложений $\ref{prop22}$, $\ref{prop23}$
основывается на следствии следующего принципа плотности
подпространства погружений в пространстве непрерывных отображений
с компактно-открытой топологией, см. [Gr, предложение 1.2.2 ].

\begin{proposition}\label{prop24}
Пусть $f_0: M \looparrowright N$ гладкое погружение между
замкнутыми многообразиями, причем многообразие $N$ снабжено
метрикой $b$ и $\dim(M) < \dim(N)$. Пусть $g: M \to N$ непрерывное
отображение, гомотопное погружению $f_0$. Тогда $\forall
\varepsilon > 0$ найдется погружение $f : M \looparrowright N $, регулярно
гомотопное погружению $f_0$, для которого $\dist(g;f)_{C^0} <
\varepsilon $ в метрике пространства отображений, индуцированной
метрикой $b$.
\end{proposition}

Нам потребуется следующее следствие Предложения 24.

\begin{corollary}\label{cor25}

Пусть $(f',\kappa,\Xi')$ -- произвольное скошенно-оснащенное
погружение, где $f': M^{n-k} \looparrowright \R^n$, представляющее
элемент $[(f',\kappa,\Xi')] \in Imm^{sf}(n-k,k)$ и пусть $f_1: M^{n-k}
\to \R^n$ -- произвольное непрерывное отображение. Тогда для
произвольного (сколь угодно малого) $\varepsilon
> 0$
существует скошенно-оснащенное погружение  $(f, \kappa, \Xi)$, где
$f : M^{n-k} \looparrowright \R^n$ -- погружение того же
многообразия, представляющая тот же элемент $[(f',\kappa, \Xi')]
\in Imm^{sf}(n-k,k)$, причем выполнено условие

\begin{eqnarray}
\dist(f_1;f) < \varepsilon. \label{dist}
\end{eqnarray}
\end{corollary}

\subsubsection*{Доказательство Следствия $\ref{cor25}$}

По предложению 24 существует погружение $f$, регулярно-гомотопное
$f'$ для которого выполнено условие ($\ref{dist}$). Регулярная
гомотопия скошенно-оснащенного погружения продолжается до
регулярной гомотопии в классе скошенно-оснащенных погружений.
Поэтому $f$ является скошенно-оснащенным и элементы
$[(f',\kappa,\Xi')]$, $[(f,\kappa,\Xi)]$ равны. Следствие
$\ref{cor25}$ доказано.
\[  \]

Нам также потребуется следующая лемма (Основная Лемма),
доказательство которой проводится в разделе 3.

\begin{lemma}\label{osnlemma1}
При $k \equiv 0 \pmod{2}$, $k\ge 8$,  существует кусочно-линейное
отображение $d: \RP^{n-k} \to \R^n$ (с особенностью) общего
положения,  допускающее циклическую структуру.
\end{lemma}

\subsubsection*{Вывод Предложения $\ref{prop22}$ из Леммы $\ref{osnlemma1}$}

Рассмотрим отображение $d: \RP^{n-k} \to \R^n$, построенное в
 Лемме 26. Рассмотрим произвольный элемент группы
 $Imm^{sf}(n-k,k)$, представленный погружением
 $f': M^{n-k} \looparrowright \R^n$ со скошенным оснащением $(\kappa,\Xi')$. Рассмотрим
  характеристический класс
 $\kappa: M^{n-k} \to \RP^{n-k}$ этого скошенного оснащения и рассмотрим отображение
 $f_1: M^{n-k} \to \R^n$, определенное как композиция $f_1 = d \circ
 \kappa$. Обозначим через $\varepsilon$ положительное число, много меньшее радиуса некоторой регулярной
 окрестности $U_{N(d)}$ полиэдра $N(d)$ (подмногообразия с особенностями) точек самопересечения
 отображения $d$. Обозначим через  $g_{N(d)}: N(d) \looparrowright \R^n$ параметризующее погружение полиэдра точек
 самопересечения.
 Воспользуемся Следствием $\ref{cor25}$ и определим
 скошенно-оснащенное погружение $f: M^{n-k} \looparrowright \R^n$,
 представляющее исходный элемент $[(f', \kappa, \Xi')]$ и при этом такое, что
 выполнено условие ($\ref{dist}$).

 Рассмотрим многообразие $N^{n-2k}$ точек самопересечения
 погружения $f$, определенное по формуле $(\ref{N})$. Определено погружение $g: N^{n-2k} \looparrowright \R^n$.
 Многообразие $N^{n-2k}$
 представим в виде
 объединения двух  подмногообразий $N^{n-2k}=N^{n-2k}_{cycl} \cup
 N^{n-2k}_{\I_b}$ по общей части границы $\partial N^{n-2k}_{cycl}= \partial
 N^{n-2k}_{\I_b}$. Все эти многообразия с краем определим далее.

 Для простоты обозначений предположим, что радиус регулярной погруженной окрестности
 погружения многообразия с краем $N^{n-2k}(d)$ точек самопересечения отображения $d$ (вдали от критических значений этого
 отображения т.е. края $\partial N^{n-2k}(d)$), радиус регулярной погруженной окрестности
 погружения, полученного ограничением отображения $d$ вне
 регулярной окрестности критических точек и радиус самой
 регулярной вложенной окрестности критических значений отображения $d$
  равны $\varepsilon$. Обозначим указанные окрестности через  $U_{N(d)}$ $U_d$, $U_{\partial
  N(d)}$ соответственно.

 Определим $N^{n-2k}_{cycl}$ как максимальное подмногообразие с краем, чтобы образ $g(N^{n-2k}_{cycl})$
 лежал в регулярной погруженной
 $\varepsilon$-окрестности точек самопересечения отображения $d$.
Это означает, что существует вложение $N^{n-2k}_{cycl} \subset
U_{N(d)}$, которое при композиции с параметризующим погружением
окрестности переходит в ограничение погружения многообразия
самопересечения на это подмногообразие. Легко определить
отображение
 $(N^{n-2k}_{cycl}, \partial  N^{n-2k}_{cycl}) \to (N^{n-2k}(d), \partial N^{n-2k}(d))$, при котором прообраз
 каждой точки состоит из пары точек, образы которых при отображении $\kappa$ $\varepsilon$-близки на $\RP^{n-k}$.

  Определим
 $N^{n-2k}_{\I_b}$ как оставшуюся часть многообразия $N^{n-2k}$, образ которой при отображении $g$ расположен в
 окрестности $U_d$ (погруженной) множества регулярных значений отображения
 $d$ (выше мы привели определение для аналогичного случая) или во вложенной окрестности $U_{\partial
  N(d)}$.

 Предположим, что погружение $g: N^{n-2k} \looparrowright \R^n$ трансверсально вдоль $\partial
 U_{N(d)}$.
 Тогда общая граница
 $\partial N^{n-2k}_{cycl} = \partial N^{n-2k}_{\I_b}$
 погружена в границу регулярной окрестности $U_{\partial N(d)}$.

Многообразие $N^{n-2k}$ снабжено характеристическим отображением
$\eta_{N}: N^{n-2k} \to K(\D_4,1)$.
 Легко проверить, что характеристическое отображение $\eta_{N}$,
 ограниченное на подмногообразие $N^{n-2k}_{\I_b} \subset N^{n-2k}$,
 гомотопно  отображению в подпространство $K(\I_b,1) \subset K(\D_4,1)$.

Структурное отображение   $\eta_{N(d)}: N^{n-2k}(d) \to K(\D_4,1)$
продолжается c $N^{n-2k}(d)$  до отображения
 $\eta_{U_{N(d)}}: U_{N(d} \to K(\D_4,1)$. При этом ограничение $\eta_{U_{\partial
 N(d)}}: U_{\partial N(d)} \to K(\D_4,1)$ гомотопно отображению в подпространство $K(\I_b,1) \subset K(\D_4,1)$.

Нетрудно проверить, что характеристическое отображение $\eta_{N}:
N^{n-2k} \to K(\D_4,1)$, ограниченное на подмногообразие
 $N^{n-2k}_{cycl}$ определено как ограничение структурного отображения $\eta_{U_N(d)}: U_{N(d)} \to
 K(\D_4,1)$, поскольку $g(N^{n-2k}_{cycl})$ лежит в $U_{N(d)}$.

По условию определена циклическая структура отображения $d$. Эта
структура  определена отображением $\mu_{a,N(d)}:
(N^{n-2k}(d),\partial N^{n-2k}(d)) \to (K(\I_a,1),K(\I_d,1))$.
Отображение $\mu_{a,N(d)}$ продолжается до отображения
$\mu_{a,U_{N(d)}}: (U_{N(d)},U_{\partial N(d)}) \to
(K(\I_a,1),K(\I_d,1))$. При этом из условия ($\ref{h mu_a}$)
вытекает, что отображения $\mu_{a,U_{N(d)}} \vert _{U_{\partial
N(d)}}$ и $\eta \vert_{U_{\partial N(d)}}$ гомотопны.
Характеристическое отображение $\eta_{N(f)} \vert_{N_{cycl}\subset
N(f)}: N_{cycl}^{n-2k} \to K(\D_4,1)$ индуцировано структурным
отображением $\eta_{U_{N(d)}}$ при включениии $N_{cycl}^{n-2k}
\subset U_{N(d)}$. Следовательно, отображения $\mu_{N_{cycl}}$ и
$\eta_{N_{\I_b}}$ можно склеить в одно отображение $\mu_{a}:
N^{n-2k} \to K(\I_a,1)$.

Докажем, что отображение $\mu_{a}$ определяет циклическую
структуру для $(f,\kappa,\Xi)$ в смысле Определения
$\ref{defcycl}$. Рассмотрим отображение многообразий с краем,
которое определяется композицией включения и проекции регулярной
окрестности на центральное подмногообразие: $F:(N_{cycl},\partial
N_{cycl}) \stackrel{g}{\subset} (U_{N(d)},U_{\partial N(d)})
\longrightarrow (N^{n-2k}(d),\partial N^{n-2k}(d))$. Из
элементарных геометрических соображений вытекает, что эта степень
этого отображения совпадает со степенью отображения $\kappa:
M^{n-k} \to \RP^{n-k}$ т.е.  с $h(f,\Xi, \kappa)$. С другой
стороны, характеристическое число $h_{\mu_a,k}$, посчитанное по
формуле $(\ref{mu_a})$, вычисляется как степень $deg(F)$
отображения $F$. Доказательство вытекает из Леммы $\ref{lemma20}$
и равенства $(\ref{h mu_a})$.

 Предложение $\ref{prop22}$ доказано.

\begin{lemma}\label{osnlemma2}
При $k \equiv 0 \pmod{2}$, $n<10k-1$, $k \ge 8$,  существует
кусочно-линейное отображение $с: S^{n-2k}/\i \to \R^n$ (с
особенностью) общего положения,  допускающее кватернионную
структуру.
\end{lemma}

Доказательство Леммы $(\ref{osnlemma2})$ аналогично доказательству
Леммы $(\ref{osnlemma1})$, оно также проводится в разделе 3.

\subsubsection*{Вывод Предложения $\ref{prop23}$ из Леммы
$\ref{osnlemma2}$}

Рассмотрим скошенно-оснащенное погружение $(f,\kappa,\Xi)$, такое,
что $\delta([(f,\kappa,\Xi)]) = [(g_1: N^{n-2k} \looparrowright
\R^n,\eta,\Psi)]$ представляет данный элемент из группы
$Imm^{\D_4}(n-2k,2k)$. По Предложению 22, не ограничивая общности,
можно считать, что погружение $f: M^{n-k} \looparrowright \R^n$
допускает циклическую структуру. Рассмотрим отображение $c:
S^{n-2k}/\i \to \R^n$, построенное в Лемме $(\ref{osnlemma2})$ и
рассмотрим погружение $g: N^{n-2k} \looparrowright \R^n$,
определенное в результате применения Предложения $\ref{prop22}$ к
отображению композиции $c \circ \mu_a: N^{n-2k} \to S^{n-2k}/\i
\to \R^n$.

Рассуждая по аналогии с доказательством Предложения
$\ref{prop22}$, заключаем, что определена кватернионная структура
$\D_4$--оснащенного погружения $[(g,\eta,\Psi)]$. Предложение
$\ref{prop23}$ доказано.

\subsubsection*{Доказательство Теоремы $\ref{th9}$}
Рассмотрим скошенно-оснащенное погружение $(f,\kappa,\Xi)$, такое,
что $\delta([(f,\kappa,\Xi)]) = [(g: N^{n-2k} \looparrowright
\R^n,\eta,\Psi)]$ представляет данный элемент из группы
$Imm^{\D_4}(n-2k,2k)$. Воспользуемся Предложением $\ref{prop23}$ и
представим элемент $[(g,\eta,\Psi)]$ $\D_4$--оснащенным
погружением, допускающим кватернионную структуру.


 Выберем натуральное
значение $k$ таким, что $n-4k=31$, $k \ge 8$, что возможно в
предположении $n \ge 63$. Обозначим через $L^{31}$ -- многообразие
самопересечения погружения $g$. Рассмотрим отображение $\lambda_a:
L^{31} \to S^{31}/\Q_a$, которое задает кватернионную структуру
для $(g, \eta, \Psi)$ и построено в Предложении $\ref{prop23}$.
Предположим, что $\lambda_a$ является трансверсальным вдоль
подмногообразия $S^7/\Q_a \subset S^{31}/\Q_a$. Рассмотрим в
многообразии $L^{31}$ подмногообразие $K^7 \subset L^{31}$,
определенное по формуле $K^7 = \lambda_a^{-1}(S^7/i)$. Определим
$\varphi = \lambda_a \vert_{K^7}$.

Докажем, что многообразие $K^7$ имеет нормальное расслоение $6
\varphi^{\ast}\chi_+$. Расслоение $\chi_+$ над $K(\Q_a,1)$ имеет
размерность $4$ и определено в начале раздела 4. Нормальное
расслоение $\nu_{L}$ изоморфно расслоению $k \zeta^{\ast}\omega$,
где $\omega: E(\omega) \to K(\Z/2^{[3]},1)$-- каноническое
4-мерное расслоение, определенное в разделе 4 перед Предложением
$\ref{prop33}$, $\zeta: L^{31} \to K(\Z/2^{[3]},1)$ --
характеристическое отображение, классифицирующее нормальное
расслоение многообразия $L^{31}$.

Легко проверить, что $k= \frac{n-1}{4}-7$ и $k \equiv 0 \pmod{8}$.
По Предложению $\ref{prop33}$ заключаем, что расслоение $\nu_{N}
\vert_{K}$ тривиально. Нормальное расслоение многообразия $K^7$
изоморфно сумме Уитни расслоений $6 \varphi^{\ast}\chi \oplus
\nu_{N} \vert_{K}$. Из-за тривиальности второго слагаемого получим
требуемый изоморфизм $\nu_K = 6 \varphi^{\ast}\chi$.

 По размерностным соображениям, не ограничивая общности, образ отображения $\lambda_a \vert_{K}$
 лежит в остове $S^7/\Q_a \subset S^{31}/\Q_a$. По Предложению $\ref{prop32}$ степень
 $\deg(\lambda_a \vert_K)$ четна. Поскольку $\deg(\varphi)=\deg(\lambda_a \vert_K)= \deg(\lambda_a)=h_a$ и
 $h_a = h(f,\kappa,\Xi)$,
 получим $h(f,\kappa,\Xi)=0$. Теорема $\ref{th9}$ доказана.


\section{Основной этап доказательства}

В этом разделе все отображения и многообразия предполагаются кусочно-линейными,
если специально не оговорено условие гладкости.

\subsubsection*{Определение  отображения    $d: \RP^{n-k} \to
\R^n$}

 Обозначим через $J_0$ стандартную $(n-k)$--мерную
сферу коразмерности $k$, $k \equiv 0 \pmod{8}$, $k>0$, которая
получена в результате джойна $\frac{n-k+1}{2}=r_0$ копий окружности
$S^1$.  Обозначим стандартное вложение $J_0$ в $\R^n$
через $i_{J_0}: J_0 \subset \R^n$.

Определено отображение $p': S^{n-k} \to J_0$, которое получается в
результате взятия джойна $r$ копий стандартных накрытий $S^1 \to
S^1/\i$. Стандартное действие $\Z/4 \times S^{n-k} \to S^{n-k}$
коммутирует с отображением $p'$. Тем самым, определено отображение
$\hat p: S^{n-k}/\i \to J_0$ и отображение $p: \RP^{n-k} \to J_0$
как композиция стандартного двулистного накрытия $\pi: \RP^{n-k}
\to S^{n-k}/\i$ с отображением $\hat p$.

Рассмотрим композицию $i_{J_0} \circ \hat p: S^{n-k}/i \to J_0 \to
\R^n$, которую мы обозначим через $\hat g$, и приведем это
отображение в отображение $\hat d$ общего положения малой
$\delta$-деформацией $(i_{J_0} \circ \hat p) \to \hat d$ (какую
именно деформацию и насколько малым следует выбирать положительное
число $\delta$, задащее калибр этой деформации, выяснится при
доказательстве Леммы $\ref{lemma30}$). Искомое отображение $d:
\RP^{n-k} \to \R^n$ определим  как гладкое отображение общего
положения, полученное в результате $C^0$ деформации калибра
$\delta'$, где $\delta' << \delta$, отображения, представленного
композицией стандартного двулистного накрытия $\pi$ и отображения
$\hat d$. В Лемме $\ref{osnlemma1}$ доказывается, что что для
отображения $d$ (с особенностью) существует циклическая структура.

\subsubsection*{Определение отображения $c: S^{n-2k}/\i \to \R^n$}

Обозначим через $J_1$ $(n-2k)$--мерный полиэдр, $k \equiv 0 \pmod
{8}$, $k>0$, который полученый в результате джойна
$\frac{n-2k+1}{4}=r_1$ копий кватернионного линзового пространства
$S^3/\Q_a$.
 Согласно теореме Масси [Ma] (см. также [Ме]) определено
вложение $S^3/\Q_a \subset \R^4$ коразмерности 1 и в предположении
$n \ge 4r-1$ определено стандартное вложение $J_1$ в $\R^n$,
которое обозначим через $i_{J_1}: J_1 \subset \R^n$.

Определено отображение $p': S^{n-2k} \to J_1$, которое получается
в результате взятия джойна $r$ копий стандартных накрытий $S^3 \to
S^3/\Q_a$. Стандартное действие $\Q_a \times S^{n-2k} \to
S^{n-2k}$, определенное прямой суммой $r$ экземпляров стандартного
действия на $S^3$, заданного формулами $(\ref{Q a1}), (\ref{Q
a2}), (\ref{Q a3})$, коммутирует с отображением $p'$. Тем самым,
определено отображение $\hat p: S^{n-2k}/\Q_a \to J_1$ и
отображение $p: S^{n-2k}/\i \to J_1$ как композиция стандартного
двулистного накрытия $\pi: S^{n-2k}/\i \to S^{n-2k}/\Q_a$ с
отображением $\hat p$.

Рассмотрим композицию $i_{J_1} \circ \hat p: S^{n-2k}/\Q_a \to J_1
\to \R^n$, которую обозначим через $\hat h$. Приведем это
отображение в отображение $\hat c$ общего положения малой
$\delta$-деформацией (какую именно следует выбрать деформацию и
насколько малым следует выбирать положительное число $\delta$,
задащее калибр этой деформации выяснится при доказательстве Леммы
$\ref{lemma31}$). Искомое отображение $c: S^{n-2k}/\i \to \R^n$
определим  как гладкое отображение общего положения, полученное в
результате $C^0$ деформации калибра $\delta'$, где $\delta' <<
\delta$, отображения, представленного композицией стандартного
двулистного накрытия $\pi$ и отображения $\hat c$.
 В Лемме
$\ref{osnlemma2}$ доказывается, что что для отображения $c$ (с
особенностью) существует кватернионная структура.

\subsubsection*{Подпространства и факторподпространства пополненного конфигурационного пространства для
$\RP^{n-k}$}

В разделе 2 было
определено пространство $\Gamma_0$, его двулистное накрытие $\bar \Gamma_0$
и структурное отображение $\eta_{\Gamma_0} : \Gamma_0 \to
K(\D_4,1)$. Пространство $\Gamma_0$ является многообразием с краем. Обозначим внутренность этого многообразия
через $\Gamma_{0\circ}$. Ограничение структурного отображения $\eta_{\Gamma_0}$ на $\Gamma_{0\circ}$
обозначим через $\eta_{\Gamma_{0\circ}} : \Gamma_{0\circ} \to
K(\D_4,1)$.

 Обозначим
через $\Sigma_{0\circ} \subset \Gamma_{0\circ}$ полиэдр двукратных особых точек
отображения $p: \RP^{n-k} \to J_0$, полученный раздутием полиэдра
$\{[(x, y)] \in \Gamma_{0\circ}, p(x) = p(y), x \ne y
\}$. Этот полиэдр снабжён структурным отображением
$\eta_{\Sigma_{0\circ}}: \Sigma_{0\circ} \to K(\D_4,1)$, которое индуцировано
ограничением структурного отображения $\eta_{\Gamma_{0\circ}}$ на
$\Sigma_{0\circ}$.

Определена свободная инволюция $T_{\RP}: \RP^{n-k} \to \RP^{n-k}$, переставляющая точки в
каждом слое стандартного двулистного накрытия $\RP^{n-k} \to S^{n-k}/\i$. На пространстве
$\bar \Gamma_{0\circ}$ действует инволюция $T_{\bar \Gamma_0}: \bar \Gamma_{0\circ} \to
\bar \Gamma_{0\circ}$, которая определена как ограничение инволюции
на $\RP^{n-k} \times \RP^{n-k}$, построенной по инволюции $T_{\RP}$ на каждом сомножителе,
на подпространство $\bar \Gamma_{0\circ} \subset
\RP^{n-k} \times \RP^{n-k}$. На факторпространстве $\Gamma_{0\circ}$ инволюции $T'$
определена факторинволюция $T_{\Gamma_{0\circ}}: \Gamma_{0\circ} \to \Gamma_{0\circ}$.

Обозначим через $\Sigma_{antidiag} \subset \Gamma_{0\circ}$
подпространство, называемое антидиагональю, которое образовано
всеми антиподальными парами  $\{[(x, y)] \in \Gamma_0 : x,y \in
\RP^{n-k}, x \ne y, T_{\RP}(x)=y \}$. Здесь и далее в обозначении
диагональных и антидиагональных пространств индекс $0$
опускается. Нетрудно проверить, что антидиагональ
$\Sigma_{antidiag} \subset \Gamma_0$ является множеством
неподвижных точек для инволюции $T_{\Gamma_0}$. Нормальное расслоение антидиагонали в $\Gamma_0$ изоморфно
касательному расслоению $T(\Sigma_{antidiag})$.

Определим пространство $\Gamma_{K_{0\circ}}$ на
котором инволюция $T_{\Gamma_{0\circ}}$ действует свободно. Определим
$\Gamma_{K_{0\circ}} = \Gamma_{0\circ} \setminus \Sigma_{antidiag}$.  Определена свободная инволюция
$T_{\Gamma_{K_{0\circ}}}: \Gamma_{K_{0\circ}} \to \Gamma_{K_{0\circ}}$.
Факторпространство $\Gamma_{K_{0\circ}}/T_{\Gamma_{K_{0\circ}}}$ обозначим через
$\hat \Gamma_{K_{0\circ}}$.

Подполиэдр $\Sigma_{0\circ} \subset \Gamma_{0\circ}$ кратных точек отображения
$p$ представляется в виде объединения $\Sigma_{0\circ} =
\Sigma_{antidiag} \cup K'_{0\circ}$, где   $K'_{0\circ}$--открытый подполиэдр, в
который входят все точки из $\Sigma_{0\circ}$, не попавшие на
антидиагональ.
Подполиэдр $K_{0\circ} \subset \Gamma_{K_{0\circ}}$ инвариантен при инволюции
$T_{\Gamma_{K_{0\circ}}}$. Обозначим $T_{\Gamma_{K_{0\circ}}} \vert_{K_{0\circ}}$ через
$T_{K_{0\circ}}$. Обозначим факторпространство
 $K_{0\circ}/ T_{K_{0\circ}}$ через  $\hat K_{0\circ}$.
Ограничение структурного отображение $\eta_{\Gamma_{0\circ}}: \Gamma_{0\circ}
\to K(\D_4,1)$ на $\Gamma_{K_{0\circ}}$ и $K_{0\circ}$ обозначим через
$\eta_{\Gamma_{0\circ}}$ и $\eta_{K_{0\circ}}$ соответственно.

Рассмотрим многообразие с особенностями, которое определено в результате замыкания открытого многообразия $\Gamma_{0\circ}$ и получено в результате присоединения
к этому открытому многообразию диагональной и антидиагональной компонент $\Sigma_{diag}$ и $\Sigma_{antidiag}$.
Обозначим замыкание $Cl(K_{0\circ})$ полиэдра $K_{0\circ}$ (соответственно замыкание
$Cl(\hat K_{0\circ})$ полиэдра $\hat K_{0\circ}$) в этом многообразии с особенностями через $K_{(0)}$ (соответственно через $\hat K_{(0)}$).
 Обозначим через
$Q_{antidiag} = \Sigma_{antidiag} \cap K_{(0)}$, $Q_{diag} =
\partial \Gamma_{diag} \cap K_{(0)}$, $Q_{diag}\subset K_{(0)}$,
$Q_{antidiag} \subset K_{(0)}$. Назовем эти подпространства компонентами границы полиэдра $K_{(0)}$.
Аналогично определим $\hat Q_{diag}$, $\hat Q_{antidiag}$ компоненты границы полиэдра $\hat K_{(0)}$.

Заметим, что структурное отображение $\eta_{K_{0\circ}}$ продолжается с $K_{0\circ}$ на компоненту границы
$Q_{antidiag}$. Обозначим это продолжение через $\eta_{Q_{antidiag}}: Q_{antidiag} \to K(\D_4,1)$.
Отображение $\eta_{Q_{antidiag}}$  представляется в виде композиции отображения
$\eta_{antidiag}: Q_{antidiag} \to K(\I_a,1)$ и отображения
включения $i_{\I_a,\D_4}: K(\I_a,1) \subset K(\D_4,1)$.

Cтруктурное
отображение $\eta_{K_{(0)}}$ на компоненту $Q_{diag}$ не продолжается.
Определено отображение $\eta_{diag}: Q_{antidiag} \to K(\I_d,1)$.
Обозначим через $U(Q_{diag}) \subset K_{0\circ}$ малую регулярную взрезанную окрестность
открытого конца, примыкающего к $Q_{diag}$. Определена проекция $proj_{diag}:U(Q_{diag}) \to Q_{diag}$
регулярной взрезанной окрестности на центальное подмногообразие.
Ограничение структурного отображения $\eta_{K_{0\circ}}$ на окрестность $U(Q_{diag})$ представлено композицией
отображения $\eta_{U(Q_{diag})}: U(Q_{diag}) \to K(\I_b,1)$ и отображения $i_{\I_b,\D_4}: K(\I_b,1) \to K(\D_4,1)$.
Гомотопические классы отображений $\eta_{K_{(0)}}$ и $\eta_{U(Q_{diag})}$ связаны равенством
$\eta_{diag} \circ proj_{diag} = p_{\I_b,\I_d} \circ \eta_{K_{(0)}}$.

Определим следующую коммутативную диаграмму мономорфизмов групп,
каждая из которых представлена в группе $O(4)$:
\begin{eqnarray}\label{HH}
\begin{array}{ccccccccc}
&&&& \I_b &&&&\\
&& \nearrow && \cap &&&&\\
\I_d& \subset &\I_a& \subset & \D_4 & \subset & \E & \subset &\H. \\
&& \searrow && \cap &&&&\\
&&&& \I_c &&&&\\
\end{array}
\end{eqnarray}

В этой диаграмме группа $\H$ имеет порядок $32$ и является
подгруппой в группе $O(4)$ и была ранее определена (см.
$(\ref{H})$) при помощи преобразований, записанных в стандартном
базисе $\{ \e_1, \e_2, \e_3, \e_4 \}$. Подгруппа $\E \subset \H$
имеет порядок $16$ и определяется как группа, порожденная
преобразованиями из $\H$, сохраняющими неупорядоченную пару
плоскостей $(\e_1,\e_3)$, $(\e_2,\e_4)$ (напимер, элемент из $\H$,
заданный преобразованием $\e_1 \to \e_2$, $\e_2 \to -\e_1$, $\e_3
\to \e_3$, $\e_4 \to \e_4$, не лежит в подгруппе $\E \subset \H$).
Включение $\D_4 \subset \E$ является центральным квадратичным
расширением группы $\D_4$ посредством элемента $c$ такого, что
$c^2$ совпадает с образующей подгруппы $\I_d \subset \D_4$.
Подгруппа $\D_4$ представлена диагональными преобразованиями $O(2)
\subset O(2) \oplus O(2) \subset O(4)$ стандартного вида в каждой
паре плоскостей $(\e_1, \e_3)$, $(\e_2, \e_4)$. Остальные
подгруппы в диаграмме $(\ref{HH})$ были определены ранее.

Определены абелевы подгруппы $\E_a,
\E_b, \E_c, \E_d$ в группе $\E$, которые являются квадратичными
расширениями соответствующих подгрупп $\I_a, \I_b, \I_c, \I_d$ в
группе $\D_4$. Нетрудно проверить, что определено естественное
отображение $\eta_{\hat K_{0\circ}}: \hat K_{0\circ} \to K(\E,1)$, которое индуцирует
отображение диаграмм 2-листных накрытий:

\begin{eqnarray}\label{14}
\begin{array}{ccccccc}
\bar K_{0\circ} & \stackrel{\bar r}{\longrightarrow} & \tilde K_{0\circ}& &K(\I_c,1)& \longrightarrow & K(\E_c,1)\\
\downarrow & & \downarrow & \longrightarrow & \downarrow & & \downarrow\\
K_{0\circ} & \stackrel{r}{\longrightarrow} & \hat K_{0\circ}& &K(\D_4,1)& \longrightarrow & K(\E,1).\\
\end{array}
\end{eqnarray}

Горизонтальные отображения между соответствующими пространствами
диаграмм мы переобозначим для краткости через $\bar \eta, \tilde \eta, \eta, \hat
\eta$ соответственно.

\[  \]
Включения подгрупп $\I_b \subset \E_b$, $\I_d \subset \E_d$, $\I_a
\subset \E_a$ входят в коммутативные диаграммы гомоморфизмов
групп:

\begin{eqnarray}\label{Ib}
\begin{array}{ccc}
\I_b & \longrightarrow & \E_b \\
p_{b,d} \downarrow & & \downarrow  \\
\I_d & \stackrel{i_{d,a}}{\longrightarrow} & \I_a
\end{array}
\end{eqnarray}

\begin{eqnarray}\label{Id}
\begin{array}{ccc}
\I_d & \longrightarrow & \E_d \\
i_{d,a} \downarrow & & \downarrow  \\
\I_a & \stackrel{id}{\longrightarrow} & \I_a
\end{array}
\end{eqnarray}

\begin{eqnarray}\label{Ia}
\begin{array}{ccc}
\I_a & \longrightarrow & \E_a \\
id \downarrow & & \downarrow  \\
\I_a & \stackrel{id}{\longrightarrow} & \I_a
\end{array}
\end{eqnarray}

На пространстве $\hat \Gamma_{K_{0\circ}}$ определена, в свою очередь,
свободная инволюция $T_{\hat \Gamma_{K_{0\circ}}} : \hat \Gamma_{K_{0\circ}} \to
\hat \Gamma_{K_{0\circ}}$ по формуле $T_{\hat \Gamma_{0\circ}} ([(x,y)]) =
[(T(x),y)]$. (Заметим, что $[(T(x),y)]=[(x,T(y)]$.) Пространства,
полученные из $\hat \Gamma_{K_{0\circ}}$, $\hat K_{0\circ}$, $\bar
\Gamma_{K_{0\circ}}$, $\bar K_{{0\circ}}$
после факторизации по инволюции
$T_{\hat \Gamma_{K_{0\circ}}}$ обозначим через $\hat
\Gamma_{K_{0\circ}}^{\downarrow}$, $\hat K_{0\circ}^{\downarrow}$ $\bar
\Gamma_{K_{0\circ}}^{\downarrow}$, $\bar K_{0\circ}^{\downarrow}$
соответственно. Определено 2-листное накрытие $\hat \Gamma_{K_{0\circ}}
\to \hat \Gamma_{K_{0\circ}}^{\downarrow}$, которое индуцирует двулистное
накрытие $\hat K_{0\circ} \to \hat K_{0\circ}^{\downarrow}$.
При этом определена коммутативная диаграмма:
\begin{eqnarray}\label{E-H}
\begin{array}{ccc}
\hat \Gamma_{K_{0\circ}} & \longrightarrow & K(\E,1) \\
\downarrow & & \downarrow  \\
\hat \Gamma_{K_{0\circ}}^{\downarrow} & \longrightarrow & K(\H,1).
\end{array}
\end{eqnarray}
\[  \]

\subsubsection*{Подпространства и факторподпространства пополненного конфигурационного пространства для
$S^{n-2k}/\i$}

В разделе 2 было
определено пространство $\Gamma_1$, его двулистное накрытие $\bar \Gamma_1$
и структурное отображение $\zeta_{\Gamma_1} : \Gamma_1 \to
K(\H,1)$. Пространство $\Gamma_1$ является многообразием с краем. Обозначим внутренность этого многообразия
через $\Gamma_{1\circ}$. Ограничение структурного отображения $\zeta_{\Gamma_1}$ на $\Gamma_{1\circ}$
обозначим через $\zeta_{\Gamma_{1\circ}} : \Gamma_{1\circ} \to
K(\H,1)$.

 Обозначим
через $\Sigma_{1\circ} \subset \Gamma_{1\circ}$ полиэдр двукратных особых точек
отображения $p: S^{n-2k}/\i \to J_1$, определенный формулой
$\{[(x, y)] \in \Gamma_{1\circ}, p(x) = p(y), x\ne y
\}$. Этот полиэдр снабжён структурным отображением
$\zeta_{\Sigma_{1\circ}}: \Sigma_{1\circ} \to K(\H,1)$, которое индуцировано
ограничением структурного отображения $\zeta_{\Gamma_{1\circ}}$ на
$\Sigma_{1\circ}$.

Определена свободная инволюция $T_{S^{n-2k}/\i}: S^{n-2k}/\i \to S^{n-2k}/\i$, переставляющая точки в
каждом слое стандартного двулистного накрытия $S^{n-2k}/\i \to S^{n-2k}/\i$. На пространстве
$\bar \Gamma_{1\circ}$ действует инволюция $T_{\bar \Gamma_{1\circ}}: \bar \Gamma_{1\circ} \to
\bar \Gamma_{1\circ}$, которая определена как ограничение инволюции
на $S^{n-2k} \times S^{n-2k}$, построенной по инволюции $T_{S^{n-2k}/\i}$ на каждом сомножителе,
на подпространство $\bar \Gamma_{1\circ} \subset
S^{n-2k}/\i \times S^{n-2k}/\i$. На факторпространстве $\Gamma_{1\circ}$ инволюции $T'$
определена факторинволюция $T_{\Gamma_{1\circ}}: \Gamma_{1\circ} \to \Gamma_{1\circ}$.

Обозначим через $\Sigma_{antidiag} \subset \Gamma_{1\circ}$
подпространство, называемое антидиагональю, которое образовано
всеми антиподальными парами  $\{[(x, y)] \in \Gamma_{1\circ} : x,y \in
S^{n-2k}/\i, x \ne y, T_{S^{n-2k}}(x)=y \}$. Здесь и далее в обозначении
диагональных и антидиагональных пространств индекс $1$
опускается. Нетрудно проверить, что антидиагональ
$\Sigma_{antidiag} \subset \Gamma_{1\circ}$ является множеством
неподвижных точек для инволюции $T_{\Gamma_{1\circ}}$.
Нормальное расслоение антидиагонали в $\Gamma_{1\circ}$ изоморфно
касательному расслоению $T(\Sigma_{antidiag})$.

Определим пространство $\Gamma_{K_{1\circ}}$ на
котором инволюция $T_{\Gamma_{1\circ}}$ действует свободно. Определим
$\Gamma_{K_{1\circ}} = \Gamma_{1\circ} \setminus \Sigma_{antidiag}$.  Определена свободная инволюция
$T_{\Gamma_{K_{1\circ}}}: \Gamma_{K_{1\circ}} \to \Gamma_{K_{1\circ}}$.
Факторпространство $\Gamma_{K_{1\circ}}/T_{\Gamma_{K_{1\circ}}}$ обозначим через
$\hat \Gamma_{K_{1\circ}}$.

Подполиэдр $\Sigma_{1\circ} \subset \Gamma_{1\circ}$ кратных точек отображения
$p$ представляется в виде объединения $\Sigma_{1\circ} =
\Sigma_{antidiag} \cup K'_{1\circ}$, где   $K'_{1\circ}$--открытый подполиэдр, в
который входят все точки из $\Sigma_{1\circ}$, не попавшие на
антидиагональ.
Подполиэдр $K_{1\circ} \subset \Gamma_{K_{1\circ}}$ инвариантен при инволюции
$T_{\Gamma_{K_{1\circ}}}$. Обозначим $T_{\Gamma_{K_{1\circ}}} \vert_{K_{1\circ}}$ через
$T_{K_{1\circ}}$. Обозначим факторпространство
 $K_{1\circ}/ T_{K_{1\circ}}$ через  $\hat K_{1\circ}$.
Ограничение структурного отображение $\zeta_{\Gamma_{1\circ}}: \Gamma_{1\circ}
\to K(\H,1)$ на $\Gamma_{K_{1\circ}}$ и $K_{1\circ}$ обозначим через
$\zeta_{\Gamma_{1\circ}}$ и $\zeta_{K_{1\circ}}$ соответственно.

Рассмотрим многообразие с особенностями, которое определено в результате замыкания открытого многообразия $\Gamma_{1\circ}$ и получено в результате присоединения
к этому открытому многообразию диагональной и антидиагональной компонент $\Sigma_{diag}$ и $\Sigma_{antidiag}$.
Обозначим замыкание $Cl(K_{1\circ})$ полиэдра $K_{1\circ}$ (соответственно замыкание
$Cl(\hat K_{1\circ})$ полиэдра $\hat K_{1\circ}$) в этом многообразии с особенностями через $K_{(1)}$ (соответственно через $\hat K_{(1)}$).
 Обозначим через
$Q_{antidiag} = \Sigma_{antidiag} \cap K_{(1)}$, $Q_{diag} =
\partial \Gamma_{diag} \cap K_{(1)}$, $Q_{diag}\subset K_{(1)}$,
$Q_{antidiag} \subset K_{(1)}$. Назовем эти подпространства компонентами границы полиэдра $K_{(1)}$.
Аналогично определим $\hat Q_{diag}$, $\hat Q_{antidiag}$ компоненты границы полиэдра $\hat K_{(1)}$.

Заметим, что структурное отображение $\zeta_{K_{1\circ}}$ продолжается с $K_{1\circ}$ на компоненту границы
$Q_{antidiag}$. Обозначим это продолжение через $\zeta_{Q_{antidiag}}: Q_{antidiag} \to K(\H,1)$.
Отображение $\zeta_{Q_{antidiag}}$  представляется в виде композиции отображения
$\zeta_{antidiag}: Q_{antidiag} \to K(\Q_a,1)$ и отображения
включения $i_{\Q_a,\H}: K(\Q_a,1) \subset K(\D_4,1)$.

Cтруктурное
отображение $\zeta_{K_{(1)}}$ на компоненту $Q_{diag}$ не продолжается.
Определено отображение $\zeta_{diag}: Q_{antidiag} \to K(\I_a,1)$.
Обозначим через $U(Q_{diag}) \subset K_{1\circ}$ малую регулярную взрезанную окрестность
открытого конца, примыкающего к $Q_{diag}$. Определена проекция $proj_{diag}:U(Q_{diag}) \to Q_{diag}$
регулярной взрезанной окрестности на центальное подмногообразие.
Ограничение структурного отображения $\zeta_{K_{1\circ}}$ на окрестность $U(Q_{diag})$ представлено композицией
отображения $\zeta_{U(Q_{diag})}: U(Q_{diag}) \to K(\I_a,1)$ и отображения $i_{\I_a,\H}: K(\I_a,1) \to K(\H,1)$.
Гомотопические классы отображений $\zeta_{K_{(1)}}$ и $\zeta_{U(Q_{diag})}$ связаны равенством
$\eta_{diag} \circ proj_{diag} = p_{\H_b,\I_a} \circ \eta_{K_{(1)}}$.

Определим следующую коммутативную диаграмму мономорфизмов групп,
каждая из которых представлена в группе $O(8)$:
\begin{eqnarray}\label{G}
\begin{array}{ccccccccccc}
&&&&&& \E_b &&&&\\
&&&&  && \cap &&&&\\
\I_d &\subset &\I_a & \subset &\Q_a& \subset & \H & \subset & \G & \subset &\G^{\downarrow}. \\
&&&& \searrow && \cap &&&&\\
&&&&&& \H_c &&&&\\
\end{array}
\end{eqnarray}

В этой диаграмме группа $\G^{\downarrow}$ определяется как подгруппа порядка $128$ в группе $O(8)$,
которая порождена  при помощи преобразований, записанных в стандартном базисе
$\{ \e_1, \e_2, \e_3, \e_4, \e_5, \e_6, \e_7, \e_8 \}$, следущего вида:

---1 преобразование в упорядоченной паре подпространств $( \e_1, \e_2, \e_3, \e_4 )$, $( \e_5, \e_6, \e_7, \e_8 )$
произвольной парой элементов из подгруппы $\Q_a \subset O(4)$.

---2  преобразование одновременно переставляющее пары векторов $(\e_1,\e_5)$, $(\e_2,\e_6)$, $(\e_3,\e_7)$,
$(\e_4,\e_8)$.

Подгруппа $\G \subset \G^{\downarrow}$ имеет порядок $64$ и
определяется как группа, порожденная преобразованиями из
$\G^{\downarrow}$, сохраняющими неупорядоченную пару
подпространств $(\e_1,\e_2,\e_5,\e_6)$, $(\e_3,\e_4,\e_7,\e_8)$
(напимер, элемент из $\G^{\downarrow}$, заданный преобразованием
$\j \in \Q_a$ в подпространстве $(\e_1,\e_2,\e_3,\e_4)$  и
тождественным преобразованием в подпространстве
$(\e_5,\e_6,\e_7,\e_8)$ не лежит в подгруппе $\G \subset
\G^{\downarrow}$). Включение $\H \subset \G$ является центральным
квадратичным расширением группы $\H$ посредством элемента
$\d_{\j}$ такого, что $\d_{\j}^2$ совпадает с образующей центра
$\I_d \subset \Q_a \subset \H$. Подгруппа $\Q_a$ представлена
диагональными преобразованиями $SO(4) \subset SO(4) \oplus SO(4)
\subset O(8)$ стандартного вида в каждой паре пространств $(\e_1,
\e_2, \e_5, \e_6)$, $(\e_3, \e_4, \e_7, \e_8)$. Подгруппа $\E_b
\subset \G$ порождена двумя коммутирующими преобразованиями. Первое преобразование имеент порядок $4$
и совпадает с преобразованием $\i \in \I_a \subset \Q_a$. Второе преобразование имеет порядок $2$, оно одновременно переставляюет пары базисных
векторов $(\e_1,\e_5)$, $(\e_2,\e_6)$, $(\e_3,\e_7)$,
$(\e_4,\e_8)$.

Подгруппа $\H_c \subset \G$ порождена собственными
преобразованиями в паре подпространств $(\e_1,\e_2,\e_3,\e_4)$,
$(\e_5,\e_6,\e_7,\e_8)$ на произвольную пару элементов из $\Q_a$.
Остальные подгруппы в диаграмме $(\ref{G})$ были определены ранее
в разделе 2.

Нетрудно проверить, что определено естественное отображение $\hat
\zeta: \hat K_{1\circ} \to K(\G,1)$, которое индуцирует отображение
диаграмм 2-листных накрытий:

\begin{eqnarray}\label{15.1}
\begin{array}{ccccccc}
\bar K_{1\circ} & \stackrel{\bar r}{\longrightarrow} & \tilde K_{1\circ}& &K(\H_c,1)& \longrightarrow & K(\G_c,1)\\
\downarrow & & \downarrow & \longrightarrow & \downarrow & & \downarrow\\
K_{1\circ} & \stackrel{r}{\longrightarrow} & \hat K_{1\circ}& &K(\H,1)& \longrightarrow & K(\G,1),\\
\end{array}
\end{eqnarray}
где подгруппа $\G_c \subset \G$ является  квадратичным расширением группы $\H_c$.

Горизонтальные отображения между соответствующими пространствами
диаграмм мы обозначим через $\bar \zeta, \tilde \zeta, \zeta, \hat
\zeta$ соответственно.

На пространстве $\hat \Gamma_{K_{1\circ}}$ определена, в свою очередь,
свободная инволюция $T_{\hat \Gamma_{K_{1\circ}}} : \hat \Gamma_{K_{1\circ}} \to
\hat \Gamma_{K_{1\circ}}$ по формуле $T_{\hat \Gamma_{K_{1\circ}}} ([(x,y)]) =
[(T_{S/i}(x),y)]$. (Заметим, что
$[(T_{S/i}(x),y)]=[(x,T_{S/i}(y)]$.) Пространства, полученные из
$\hat \Gamma_{K_{1\circ}}$, $\hat K_{1\circ}$, $\bar \Gamma_{K_{1\circ}}$, $\bar K_{1\circ}$
после факторизации по инволюции $T_{\hat \Gamma_{K_{1\circ}}}$
обозначим через $\hat \Gamma_{K_{1\circ}}^{\downarrow}$, $\hat
K_{1\circ}^{\downarrow}$ $\bar \Gamma_{K_{1\circ}}^{\downarrow}$, $\bar
K_{1\circ}^{\downarrow}$ соответственно. Определено 2-листное накрытие
$\hat \Gamma_{K_{1\circ}} \to \hat \Gamma_{K_{1\circ}}^{\downarrow}$, которое
индуцирует двулистное накрытие $\hat K_{1\circ} \to \hat
K_{1\circ}^{\downarrow}$. При этом определена коммутативная диаграмма:
\begin{eqnarray}\label{G-Gdownarrow}
\begin{array}{ccc}
\hat \Gamma_{K_{1\circ}} & \longrightarrow & K(\G,1) \\
\downarrow & & \downarrow  \\
\hat \Gamma_{K_{1\circ}}^{\downarrow} & \longrightarrow & K(\G^{\downarrow},1).
\end{array}
\end{eqnarray}
\[  \]


\subsubsection*{Полиэдры $K_{0}, \hat K_0$}
Введем обозначение:
\begin{eqnarray}\label{r0max}
r_{min,0}= \frac{n-7}{4}; \quad i_{max,0}=r_0-r_{min,0}.
\end{eqnarray}

Определим подпространство   $\hat K_{0\circ}^{(i)} \subset \hat
K_{0\circ}$ по формуле:
$$
\hat K_{0\circ}^{(i)} = p_{\hat K_{0\circ}}^{-1}(J_0^{(i)} \setminus
J_0^{(i+1)}),
$$
где $p_{\hat K_{0\circ}}: K_{0\circ} \to J_0$--естественная
проекция особенности на свой образ. Определим подпространство
$\hat K_{[0]\circ} \subset K_{0\circ}$ по формуле:
$$
\hat K_{[0]\circ} = p_{\hat K_{0\circ}}^{-1}(J_0 \setminus
J_0^{(i_{max,0})+1}).
$$
При этом справедлива формула:
$$\hat K_{[0]\circ} = \cup_{i=0}^{i=i_{max,0}} \hat K_{0\circ}^{(i)}. $$
Топология на пространстве $\hat K_{[0]\circ}$  индуцирована из
топологии локально-связного пространства в правой части формулы,
полученного при склейке дизъюнктного объединения пространств.
Определено двулистное накрывающее $K_{0\circ}^{(i)}$,
(соответственно $K_{[0]\circ}$) над пространством $\hat
K_{0\circ}^{(i)}$ (соответственно $\hat K_{[0]\circ}$).

Определим пространство $K_0$ (соответственно $\hat K_0$) как
пространство, полученное замыканием пространства $K_{[0]\circ}$
(соответственно $\hat K_{[0]\circ}$). Определим естественные
отображения  $\hat K_{0} \to \hat K_{(0)}$, $K_{0} \to K_{(0)}$.
Эти отображения будут $PL$--вложениями всюду, за исключением
прообразов диагонали и антидиагонали, где имеются особенности.

\subsubsection*{Разрешение особенностей полиэдра $K_{0}$}

 Для каждого $i$, $0 \le i \le i_{max,0}$, мы построим полиэдр $K^{(i)}_0$ и пространство
$RK_0^{(i)}$, которое назовем пространством, разрешающим
особенности полиэдра $K^{(i)}_0$. Параметр $i$ назовем глубиной.
Пространства $K^{(i)}_0$, $RK_0^{(i)}$ кроме того, являются
накрывающими пространствами при 2-листных накрытиях $r: K^{(i)}_0
\to \hat K^{(i)}_0$, $Rr: RK_0^{(i)} \to R\hat K_0^{(i)}$, эти
накрытия включены в следующую коммутативную диаграмму
отображений:

\begin{eqnarray}\label{16.1}
\begin{array}{ccccc}
 &  &  RK_0^{(i)} & \stackrel {pr}{\longrightarrow} & K_0^{(i)}  \\
&&&&\\
  & \phi_0 \swarrow  &  \downarrow Rr & & \downarrow r \\
&&&&\\
K(\I_a,1) &   \stackrel {\hat \phi_0}{\longleftarrow} & R\hat
K_0^{(i)}  &
\stackrel {p\hat r}{\longrightarrow} & \hat K_0^{(i)}.\\
\end{array}
\end{eqnarray}



Введем обозначения $R\hat Q_{diag}^{(i)} = (p \hat r)^{-1}(\hat
Q_{diag}^{(i)})$, $R\hat Q_{antidiag}^{(i)} = (p \hat r)^{-1}(\hat
Q_{antidiag}^{(i)})$. Аналогичные обозначения $RQ_{diag}^{(i)}$,
$RQ_{antidiag}^{(i)}$ введем для соответствующих двулистных
накрывающих. Рассматриваемые пространства включены в следующие
коммутативные диаграммы:


\begin{eqnarray}\label{18.1}
\begin{array}{ccc}
RQ_{antidiag}^{(i)} & \qquad \stackrel{pr}{\longrightarrow} \qquad & Q_{antidiag}^{(i)} \\
\phi_0 \searrow & & \swarrow \eta_{antidiag} \\
& K(\I_a,1), &
\end{array}
\end{eqnarray}

\begin{eqnarray}\label{18.2}
\begin{array}{ccc}
RQ_{diag}^{(i)} & \qquad \stackrel{pr}{\longrightarrow} \qquad & Q_{diag}^{(i)} \\
\phi_0 \searrow & & \swarrow \eta_{diag} \\
& K(\I_d,1), &
\end{array}
\end{eqnarray}

Для доказательства основного результата раздела нам потребуется
следующая лемма.

\begin{lemma}\label{lemma28}

--1. Для каждого $i$, $0 \le i \le i_{max,0}$, cуществует
пространство $\hat K_0^{(i)}$, его двулистное накрывающее
$K_0^{(i)}$ и база двулистного накрытия $\hat
K_0^{(i),\downarrow}$. Пространства $\hat K_0^{(i)}$, $K_0^{(i)}$,
$\hat K_0^{(i),\downarrow}$ являются подпространствами пространств
$\hat K_0$, $K_0$, $\hat K_0^{\downarrow}$ соответственно.
Существует пространство $R\hat K_0^{(i)}$ и его двулистное
накрывающее $RK_0^{(i)}$ и база 2-листного накрытия $R\hat
K_0^{(i),\downarrow}$, которые включаются в коммутативную
диаграмму $(\ref{16.1})$. На двулистном накрывающем выполняются
граничные условия, которые определены диаграммами $(\ref{18.1})$,
$(\ref{18.2})$.

--2. Cуществует пространство $R\hat K_0$ и его двулистное
накрывающее $RK_0$ и база 2-листного накрытия $R\hat
K_0^{\downarrow}$, которые включаются в нижеследующую
коммутативную диаграмму   $(\ref{16.2})$, аналогичную диаграмме
$(\ref{16.1})$. На двулистном накрывающем выполняются граничные
условия, которые определены диаграммами, аналогичными диаграммам
$(\ref{18.1})$, $(\ref{18.2})$.

\end{lemma}

\begin{eqnarray}\label{16.2}
\begin{array}{ccccc}
 &  &  RK_0 & \stackrel {pr}{\longrightarrow} & K_0  \\
&&&&\\
  & \phi_0 \swarrow  &  \downarrow Rr & & \downarrow r \\
&&&&\\
K(\I_a,1) &   \stackrel {\hat \phi_0}{\longleftarrow} & R\hat K_0
&
\stackrel {p\hat r}{\longrightarrow} & \hat K_0.\\
\end{array}
\end{eqnarray}

\subsubsection*{Полиэдры $K_{1}, \hat K_1$}
Введем обозначение:
\begin{eqnarray}\label{r0max}
r_{min,1}= \frac{n-15}{8}; \quad i_{max,1}=r_1-r_{min,1}.
\end{eqnarray}

Определим подпространство   $\hat K_{1\circ}^{(i)} \subset \hat
K_{1\circ}$ по формуле:
$$
\hat K_{1\circ}^{(i)} = p_{\hat K_{1\circ}}^{-1}(J_1^{(i)} \setminus
J_1^{(i+1)}),
$$
где $p_{\hat K_{1\circ}}: K_{1\circ} \to J_1$--естественная
проекция особенности на свой образ. Определим подпространство
$\hat K_{[1]\circ} \subset K_{1\circ}$ по формуле:
$$
\hat K_{[1]\circ} = p_{\hat K_{1\circ}}^{-1}(J_1 \setminus
J_1^{(i_{max,1})+1}).
$$
При этом справедлива формула:
$$\hat K_{[1]\circ} = \cup_{i=0}^{i=i_{max,1}} \hat K_{1\circ}^{(i)}. $$
Топология на пространстве $\hat K_{[1]\circ}$  индуцирована из
топологии локально-связного пространства в правой части формулы,
полученного при склейке дизъюнктного объединения пространств.
Определено двулистное накрывающее $K_{1\circ}^{(i)}$,
(соответственно $K_{[1]\circ}$) над пространством $\hat
K_{1\circ}^{(i)}$ (соответственно $\hat K_{[1]\circ}$).

Определим пространство $K_1$ (соответственно $\hat K_1$) как
пространство, полученное замыканием пространства $K_{[1]\circ}$
(соответственно $\hat K_{[1]\circ}$). Определим естественные
отображения  $\hat K_{1} \to \hat K_{(1)}$, $K_{1} \to K_{(1)}$.
Эти отображения будут $PL$--вложениями всюду, за исключением
прообразов диагонали и антидиагонали, где имеются особенности.

\subsubsection*{Разрешение особенностей полиэдра $K_{1}$}

 Для каждого $i$, $0 \le i \le i_{max,1}$, мы построим полиэдр $K^{(i)}_1$ и пространство
$RK_1^{(i)}$, которое назовем пространством, разрешающим
особенности полиэдра $K^{(i)}_1$. Параметр $i$ назовем глубиной.
Пространства $K^{(i)}_1$, $RK_1^{(i)}$ кроме того, являются
накрывающими пространствами при 2-листных накрытиях $r: K^{(i)}_1
\to \hat K^{(i)}_1$, $Rr: RK_1^{(i)} \to R\hat K_1^{(i)}$, эти
накрытия включены в следующую коммутативную диаграмму
отображений:

\begin{eqnarray}\label{20.1}
\begin{array}{ccccc}
 &  &  RK_1^{(i)} & \stackrel {pr}{\longrightarrow} & K_1^{(i)}  \\
&&&&\\
  & \phi_1 \swarrow  &  \downarrow Rr & & \downarrow r \\
&&&&\\
K(\Q_a,1) &   \stackrel {\hat \phi_1}{\longleftarrow} & R\hat
K_1^{(i)}  &
\stackrel {p\hat r}{\longrightarrow} & \hat K_1^{(i)}.\\
\end{array}
\end{eqnarray}



Введем обозначения $R\hat Q_{diag}^{(i)} = (p \hat r)^{-1}(\hat
Q_{diag}^{(i)})$, $R\hat Q_{antidiag}^{(i)} = (p \hat r)^{-1}(\hat
Q_{antidiag}^{(i)})$. Аналогичные обозначения $RQ_{diag}^{(i)}$,
$RQ_{antidiag}^{(i)}$ введем для соответствующих двулистных
накрывающих. Рассматриваемые пространства включены в следующие
коммутативные диаграммы:


\begin{eqnarray}\label{22.1}
\begin{array}{ccc}
RQ_{antidiag}^{(i)} & \qquad \stackrel{pr}{\longrightarrow} \qquad & Q_{antidiag}^{(i)} \\
\phi_1 \searrow & & \swarrow \eta_{antidiag} \\
& K(\Q_a,1), &
\end{array}
\end{eqnarray}

\begin{eqnarray}\label{23.1}
\begin{array}{ccc}
RQ_{diag}^{(i)} & \qquad \stackrel{pr}{\longrightarrow} \qquad & Q_{diag}^{(i)} \\
\phi_1 \searrow & & \swarrow \eta_{diag} \\
& K(\I_a,1), &
\end{array}
\end{eqnarray}

Для доказательства основного результата раздела нам потребуется
следующая лемма.

\begin{lemma}\label{lemma29}

--1. Для каждого $i$, $0 \le i \le i_{max,1}$, cуществует
пространство $\hat K_1^{(i)}$ и его двулистное накрывающее
$K_1^{(i)}$. Пространства $\hat K_1^{(i)}$, $K_1^{(i)}$ являются подпространствами пространств
$\hat K_1$, $K_1$ соответственно.
Существует пространство $R\hat K_1^{(i)}$ и его двулистное
накрывающее $RK_1^{(i)}$, которые включаются в коммутативную
диаграмму $(\ref{20.1})$. На двулистном накрывающем выполняются
граничные условия, которые определены диаграммами $(\ref{22.1})$,
$(\ref{23.1})$.

--2. Cуществует пространство $R\hat K_1$ и его двулистное
накрывающее $RK_1$,
которые включаются в нижеследующую коммутативную диаграмму
$(\ref{16.22})$, аналогичную диаграмме $(\ref{20.1})$. На
двулистном накрывающем выполняются граничные условия, которые
определены диаграммами, аналогичными диаграммам $(\ref{22.1})$,
$(\ref{23.1})$.

\end{lemma}

\begin{eqnarray}\label{16.22}
\begin{array}{ccccc}
 &  &  RK_1 & \stackrel {pr}{\longrightarrow} & K_1  \\
&&&&\\
  & \phi_1 \swarrow  &  \downarrow Rr & & \downarrow r \\
&&&&\\
K(\Q_a,1) &   \stackrel {\hat \phi_1}{\longleftarrow} & R\hat K_1
&
\stackrel {p\hat r}{\longrightarrow} & \hat K_1.\\
\end{array}
\end{eqnarray}


\subsubsection*{Полиэдры $K_{1}, \hat K_1$}
Введем обозначение:
\begin{eqnarray}\label{r1max}
r_{min,1}= \frac{n-15}{8}; \quad i_{max,1}=r_1 - r_{min,1}.
\end{eqnarray}

Определим пространство $\hat K_{(1)}^{(i)}$ по формуле:
$$
\hat K_{(1)}^{(i)} = p_{\hat K_{(1)}}^{-1}(J_1^{(i)} \setminus
J_1^{(i+1)}).
$$
Определим пространство $\hat K_{(1)}$ по формуле:
$$
\hat K_{(1)} = p_{\hat K_{(1)}}^{-1}(J_1 \setminus
J_1^{(i_{max,1})+1}) = \cup_{i=0}^{i=i_{max,1}}\hat K_{(1)}^{(i)}.
$$
Определено двулистное накрывающее $K_{(1)}^{(i)}$, (соответственно $K_{(1)}$) над пространством $\hat K_{(1)}^{(i)}$ (соответственно $\hat K_{(1)}$).

Для каждого $i$, $0 \le i \le i_{max,1}$, ниже определим
пространство $\hat K_{1}^{(i)}$ и пространство $\hat K_{1}$.
Определим также естественные отображения  $\hat K_{(1)}^{(i)} \to
\hat K_{1}^{(i)}$, $\hat K_{(1)} \to \hat K_{(1)}$. Эти
отображения будут $PL$--гомеоморфизмами всюду, за исключением
участков раздутой диагонали и антидиагонали, где имеется
вырождение. Аналогично определим двулистные накрывающие
$K_{1}^{(i)}$ и пространство $K_{1}$ и естественные отображения
$K_{(1)}^{(i)} \to K_{1}^{(i)}$, $K_{(1)} \to K_{1}$.

\subsubsection*{Начало доказательства Леммы $\ref{osnlemma1}$}

При определении отображения $d$ определялось также отображение
$\hat g: S^{n-k}/\i \to \R^{n}$, которое получено малой
$PL$-деформацией общего положения калибра $\delta$ из отображения
$i_{J_0} \circ p: S^{n-k}/\i \to J_0 \subset \R^n$. (Сама эта
деформация и её калибр $\delta$ будут явно указаны при
доказательстве Леммы $\ref{lemma30}$). Рассмотрим также
отображение $g
 = \pi \circ \hat d: \RP^{n-k} \to S^{n-k}/\i \to
\R^n$.

Рассмотрим полиэдр $ \hat N_{K_0}^{\downarrow}$ точек
самопересечения отображения $\hat g$ и двулистное накрывающее над
этим полиэдром $\hat N_{K_0}$, определенное в соответствии с
правым нижним пространством в диаграмме $(\ref{14})$.
 Чтобы подчеркнуть аналогию со случаем
гладких отображений, далее все рассматриваемые полиэдр (раздутие
особенностей на диагонали и антидиагонали в полиэдрах не
рассматривается) точек самопересечения будет называться
PL-многообразием с особенностями. Эти многообразия с особенностями
имеют особый край, состоящий из точек диагонали и антидиагонали.
Полиэдры без точек диагонали и антидиагонали  обозначаются $ \hat
N_{K_0\circ}^{\downarrow}$, $ \hat N_{K_0\circ}$ соответственно.
Каждый такой полиэдр назовем открытым многообразием с
особенностями. Открытое многообразие с особенностями  $ \hat
N_{K_0\circ}$ естественно вложено в открытое многообразие $\hat
\Gamma_{0\circ}$.

 Многообразие с особенностями $\hat
N_{K_0}$, в свою очередь, само служит базой двулистного накрытия
$r: N_{K_0} \to \hat N_{K_0}$. Накрытия определены в соответствии
с обозначениями в диаграмме ($\ref{14}$). Определено открытое
многообразие с особенностями $N_{K_{0\circ}}$ c краем, которое
естественно вкладывается во взрезанный квадрат $\Gamma_0$. При
этом вложении образ $N_{K_{0\circ}}$ лежит в регулярной
окрестности полиэдра $K_{(0)}$ и определена коммутативная
диаграмма отображений

\begin{eqnarray}\label{24}
\begin{array}{ccccc}
N_{K_{0\circ}} & \subset & U_{\varepsilon} & \longrightarrow & K(\D_4,1) \\
\downarrow & &\downarrow && \downarrow \\
\hat N_{K_{0\circ}} &\subset & \hat U_{\varepsilon} & \longrightarrow & K(\E,1), \\
\end{array}
\end{eqnarray}
в которой через $U_{\varepsilon}$, $\hat U_{\varepsilon}$
обозначены регулярные $\varepsilon$-окрестности подполиэдров
$K_{(0)} \subset \Gamma_0$, $\hat K_{(0)} \subset \hat \Gamma_0$
соответственно. Правые горизонтальные стрелки этой диаграммы
являются  продолжениями отображений диаграммы ($\ref{14}$) на
регулярную окрестность $U_{\varepsilon}$ и на факторпространство
$\hat U_{\varepsilon}$ этой окрестности в многообразиях с краем
$\Gamma_0$, $\hat \Gamma_0$ соответственно.

Рассмотрим многообразие с особенностями $N_{K_0}$ и обозначим его
компоненты края через $N_{Q,diag} \cup N_{Q,antidiag}$.
Переобозначим многообразие с особенностями $N_{K_0}$ через
$N_{K_0}^{(0)}$, а его компоненты края через $N^{(0)}_{Q,diag}
\cup N^{(0)}_{Q,antidiag}$. Ниже в доказательстве Леммы
$\ref{lemma30}$ определена стратификация:
\begin{eqnarray}\label{strat0}
N_{K_0}^{i_{max,0}} \subset \dots \subset N_{K_0}^{0},
\end{eqnarray}
(эта стратификация совпадает с каноническим двулистным накрытием
стратификации $(\ref{stra0})$, $(\ref{stra00})$).

Обозначим через
\begin{eqnarray}\label{strat11}
W_{K_0}^{(i)} =  N_{K_0}^{i}  \setminus N_{K_0}^{i+1}.
\end{eqnarray}
открытое многообразие с особенностями, содержащее компоненты края.

Обозначим компоненты его края $N^{i}_{Q,diag} \setminus
N^{i+1}_{Q,diag}$, $N^{i}_{Q,antidiag} \setminus
N^{i+1}_{Q,antidiag}$ через $W^{i}_{Q,diag}$ и
$W_{Q,antidiag}^{i}$ соответственно. Введем аналогичные
обозначения $\hat N_{K_0}^{i} \setminus \hat N^{i+1}_{K_0} = \hat
W_{K_0}^{(i)}$, $\hat N^{i}_{Q,diag} \setminus \hat
N^{i+1}_{Q,diag} = \hat W^{(i)}_{Q,diag}$, $\hat
N^{i}_{Q,antidiag} \setminus \hat N^{i+1}_{Q,antidiag} = \hat
W^{(i)}_{Q,antidiag}$. При вложении $W_{K_0}^{(i)} \subset
U_{\varepsilon} \subset \Gamma_0$ подмногообразие с осбенностями
$W_{Q,diag}^{(i)}$ вкладывается в компоненту края $\Gamma_{diag}$,
а подмногообразие c особенностями $W_{Q,antidiag}^{(i)}$
вкладывается в $\Gamma_{antidiag}$. Аналогичные утверждения
справедливы для компонент границы $\hat W_{\Gamma_0}^{(i)}$.

Для каждого $i$, $0 \le i \le i_{max,0}$, определена коммутативная
диаграмма ($\ref{25}$) отображений полиэдров с выписанными под
диаграммой граничными условиями. В пространства, стоящие в третьей
строке этой диаграммы, отображаются соответственно пространства
центральной строки коммутативной диаграммы ($\ref{26}$), которая
также определена для произвольного $i$, $0 \le i \le i_{max,0}$,
соответствующего номеру страта полиэдра точек самопересечения в
стратификации ($\ref{strat0}$).

\begin{eqnarray}\label{25}
\begin{array}{ccccccc}
K(\I_a,1) &&&& K(\I_a,1) &&\\
&&&&&&\\
\uparrow \hat \phi & \nwarrow \hat \phi &&& \uparrow \phi & \nwarrow \phi &\\
&&&&&&\\
 R \hat K_0^{(i)} &  \supset & R\hat Q_{diag}^{(i)} \cup R\hat Q_{antidiag}^{(i)} &
 & RK_0^{(i)}  & \longleftarrow & RQ_{diag}^{(i)} \cup RQ_{antidiag}^{(i)} \\
&&&&&&\\
 \downarrow p\hat r  &  & \downarrow& & \downarrow pr & & \downarrow  \\
&&&&&&\\
\hat K_0^{(i)} & \supset &   \hat Q_{diag}^{(i)} \cup \hat
Q_{antidiag}^{(i)} & & K_0^{(i)}
& \supset & Q_{diag}^{(i)} \cup Q_{antidiag}^{(i)} \\
&&&&&&\\
\cup &  & \cup && \cup  &  & \cup \\
&&&&&&\\
\hat K_{0\circ}^{(i)} & \supset &  \emptyset & & K_{0\circ}^{(0)}
& \supset & \emptyset \\
&&&&&&\\
\downarrow \hat \eta & &   & &   \downarrow  \eta & &    \\
&&&&&&\\
K(\E,1) &  &  &  &     K(\D_4,1).&& \\

\end{array}
\end{eqnarray}

Граничные условия в диаграмме ($\ref{25}$) на пространстве
$Q_{antidiag}^{(i)}$ записываются формулой $\phi  = i_a \circ
\eta_{Q_{antidiag}^{(i)}} : Q_{antidiag}^{(i)} \longrightarrow
K(\I_a,1) \stackrel{i_a}{\longrightarrow} K(\D_4,1)$, на
пространстве $Q_{diag}^{(i)}$ записываются формулой $\phi =
i_{d,a} \circ p_{b,d} \circ \eta_{Q_{diag}^{(i)}}: Q_{diag}^{(i)}
\longrightarrow K(\I_b,1) \longrightarrow K(\I_d,1)
\longrightarrow K(\I_a,1)$.

\begin{eqnarray}\label{26}
\begin{array}{ccccccc}
K(\I_a,1) &  &   & & K(\I_a,1)   & &  \\
&&&&&&\\
\uparrow \hat \mu_a  &  & \nwarrow \hat \mu_a & & \uparrow \mu_a & & \nwarrow  \mu_a \\
&&&&&&\\
\hat W_{K_0}^{(i)} & \supset & \hat W_{Q,diag}^{(i)} \cup \hat
W_{Q,antidiag}^{(i)} &
&W_{K_0}^{(i)} & \supset & W_{Q,diag}^{(i)} \cup W_{Q,antidiag}^{(i)} \\
&&&&&&\\
\cup  &  & \cup & & \cup & & \cup  \\
&&&&&&\\
\hat W_{K_{0\circ}}^{(i)} & \supset & \emptyset &
&W_{K_{0\circ}}^{(i)} & \supset & \emptyset \\
&&&&&&\\
\downarrow \hat \eta & &  & &   \downarrow \eta & &  \\
&&&&&&\\
K(\E,1) &  &  &  &     K(\D_4,1).&& \\
\end{array}
\end{eqnarray}

Граничные условия в диаграмме (\ref{26}) на пространстве
$W_{Q,antidiag}^{(i)}$ записываются формулой
 $\phi  = i_a \circ
\eta_{W_{Q,antidiag}^{(i)}} : W_{Q,antidiag}^{(i)} \longrightarrow
K(\I_a,1) \stackrel{i_a}{\longrightarrow} K(\D_4,1)$, на
пространстве $W_{Q,diag}^{(i)}$ записываются формулой $\phi = i_a
\circ p_b \circ \eta_{W_{Q,diag}^{(i)}}: W_{Q,diag}^{(i)}
\longrightarrow K(\I_b,1) \longrightarrow K(\I_d,1)
\longrightarrow K(\I_a,1)$.

Воспользуемся следующей леммой, доказательство которой проводится в конце раздела.

\begin{lemma}\label{lemma30}
Существует малая $PL$-деформация общего положения $(i_{J_0} \circ
\hat p) \mapsto \hat d$ такая, что для каждого кусочно-линейного
страта $W_{K_0}^{(i)}$ полиэдра $N(d)$ в средней строке
коммутативной диаграммы $(\ref{26})$ определены отображения $\hat
t^{(i)}: W^{(i)}_{\hat K_0} \to R\hat K_0^{(i)}$, $ t^{(i)}:
W^{(i)}_{K_0} \to RK_0^{(i)}$, называемые отображениями разрешения
особенностей, в соответствующие пространства второй строки
диаграммы $(\ref{25})$.
 При этом определены отображения
$\hat \mu_a^{(i)}: W_{\hat K_0}^{(i)} \to K(\I_a,1)$,
$\mu_a^{(i)}: W_{K_0}^{(i)} \to K(\I_a,1)$ в диаграмме
$(\ref{26})$ по формулам $\hat \mu_a^{(i)} = \hat \phi \circ \hat
t^{(i)}$, $\mu_a^{(i)} = \phi \circ t^{(i)}$.

Построенные семейства отображений при различных $i$ в
$(\ref{strat0})$ склеиваются до отображений
$\hat t: \hat N_{K_0} \to R\hat K_0$ (соответственно $t: N_{K_0} \to RK_0$)
и определены индуцированные отображения
$\hat \mu_a = \hat \phi \circ \hat t: \hat
N_{K_0} \to K(\I_a,1)$ (соответственно $\mu_a = \phi \circ t: N_{K_0} \to K(\I_a,1)$).
Более того, выполнено соотношение $r \circ \hat
\mu_a = \mu_a: N_{K_0} \to \hat N_{K_0} \to K(\I_a,1)$ и выполнены
граничные условия на компонетах $N_{Q,diag}$, $N_{Q,antidiag}$
границы стратифицированного многообразия с особенностями
$N_{K_0}$, которые аналогичны выписаным под формулой $(\ref{26})$
граничным условиям для компонент $W_{Q,diag}^{(i)}$,
$W_{Q,antidiag}^{(i)}$ соответственно.
\end{lemma}

\subsubsection*{Окончание доказательства Леммы $\ref{osnlemma1}$}
Докажем, что отображение $\mu_a$, определенное из Леммы
$\ref{lemma30}$ на многообразии с особенностями с краем $N_{K_0}$,
продолжается до отображения на всем $N(d)$ и указанное продолжение
определяет циклическую структуру для отображения $d$.

Напомню, что $N(d) = N_{antidiag} \cup N_{K_0}$ по общему краю
$N_{Q,antidiag}$. Полиэдр $N_{K_0}$ является базой 2-листного
накрытия $r: N_{K_0} \to \hat N_{K_0}$ и корректно определена
диаграмма $(\ref{24})$.  На компоненте $N_{K_0}$ циклическая
структура определена в результате продолжения отображения $\mu_a$
с $W_{K_0}$. В силу выполнения граничных условий для отображения
$\phi$ на $RQ_{antidiag}$, $RQ_{diag}$ следует, что отображение
$\mu_a$ имеет правильные граничные условия на $N_{Q,diag}$ и на
$N_{Q,antidiag}$. В частности, на $N_{Q,antidiag}$ отображение
$\mu_a$ совпадает со структурным отображением, которое является
циклическим.

Проверим условие ($\ref{h mu_a}$) в определении циклической
структуры для отображений с особенностью. Воспользуемся Леммой 18,
согласно которой достаточно рассмотреть абсолютный цикл $\mu_{\bar
R}: \bar R^{n-k} \to K(\I_d,1)$ (напомню, что этот цикл получен
при заклейке отображения $\mu_a: N(d) \to K(\I_a,1)$ вдоль границы
 двумя
одинаковыми копиями относительного цикла, которые выбираются произвольными) и проверить, что цикл
$\mu_{\bar R}$ определяет образующую группы гомологий
$H_{n-k}(K(\I_d,1);\Z/2)$. Обозначим эту образующую через $x \in
H_{n-k}(K(\I_d,1);\Z/2)$.

Рассмотрим еще один гомологический класс $p_{c,\ast} \circ \bar
\eta_{\ast}([\bar N(d)]) \in H_{n-k}(K(\I_d,1);\Z/2)$, который
определен циклом, полученным при помощи композиции канонического
накрытия $\bar \eta: \bar N(d) \to K(\I_c,1)$ над структурным
отображением $\eta: N(d) \to K(\D_4,1)$ с отображением $p_c:
K(\I_c,1) \to K(\I_d,1)$. Многообразие с особенностями $\bar N(d)$
рассматривается как замкнутое, т.е. каноническое накрытие берется
разветвленным над краем $\partial N(d)$. Обозначим построенный
класс гомологий через $y \in H_{n-k}(K(\I_d,1);\Z/2)$.

Докажем, что $y$ является образующим классом гомологий. Далее
проверим, что $x=y$. Действительно, класс когомологий $y^{op}$,
двойственный по Пуанкаре классу гомологий $y$, вычисляется
как нормальный характеристический класс $\bar w_{k} \in
H^{k}(\RP^{n-k};\Z/2)=H^{k}(K(\I_d,1);\Z/2)$. (При $n=2^l-1$
когомологический класс $\bar w_{k}$ не равен нулю.)

Вспомним, что $N(d) = N_{antidiag} \cup_{N_{Q,antidiag}} N_{K_0}$.
Поэтому двулистное накрывающее $\bar N(d)$ также представляется в
виде объединения двулистных накрывающих по формуле $\bar N(d) =
\bar N_{antidiag} \cup_{\bar N_{Q,antidiag}} \bar N_{K_0}$.
Справедлива аналогичная формула $\bar R = \bar R_{Q,antidiag} \cup
\bar R_{K_0}$. При этом $\bar R_{Q,antidiag}$ и $\bar
N_{Q,antidiag}$ PL-гомеоморфны как многообразия c особенностями с
краем, а отображения $\bar \mu_a$, $p_c \circ \bar \eta$ на
указанной общей части совпадают.

Нетрудно проверить, что многообразие с особенностями $\bar
R_{K_0}$ является двулистным накрывающим над некоторым
многообразием с особенностями $\tilde R_{K_0}$. Действительно,
многообразие с особенностями $N_{K_0}$ служит двулистным
накрывающим для $\hat N_{K_0}$. Здесь и далее для многообразий с
особенностями используются обозначения, соответствующие
обозначениям групп в диаграмме ($\ref{14}$). При этом отображение
$\bar \mu_{a} \vert_{R_{K_0}} : R_{K_0} \to K(\I_d,1)$
пропускается через двулистное накрытие $\bar R_{K_0} \to \tilde
R_{K_0}$, поскольку само отображение $\mu_a \vert _{N_{K_0}}$
накрывает отображение $\hat \mu_a \vert_{\hat N_{K_0}}$. Также
легко проверить, что многообразие $\bar N_{\Gamma_0}$ является
двулистным накрывающим над $\bar N_{K_0}^{\downarrow}$, при этом
отображение $p_c \circ \bar \eta \vert_{\bar N_{K_0}}$
пропускается через двулистное накрытие $\bar N_{K_0} \to \bar
N_{K_0}^{\downarrow}$. Поэтому каждый из циклов $x$, $y$
гомологичен циклу с носителем на общей части $\bar N_{antidiag}$.

 Лемма $\ref{osnlemma1}$ доказана.

\subsubsection*{Начало доказательства Леммы $\ref{osnlemma2}$}

При определении отображения $c$ определялось также отображение
$\hat h: S^{n-2k}/\Q_a \to \R^{n}$, которое получено малой
$PL$-деформацией общего положения калибра $\delta$ из отображения
$i_{J_1} \circ p: S^{n-2k}/\Q_a \to J_1 \subset \R^n$. (Сама эта
деформация и её калибр $\delta$ будут явно указаны при
доказательстве Леммы $\ref{lemma31}$ которое полностью аналогично
доказательству Леммы $\ref{lemma30}$.) Рассмотрим также
отображение $h
 = \pi \circ \hat c: S^{n-2k}/\i \to S^{n-2k}/\Q_a \to
\R^n$.

Рассмотрим полиэдр $ \hat L_{K_1}^{\downarrow}$ точек
самопересечения отображения $\hat h$ и двулистное накрывающее над
этим полиэдром $\hat L_{K_1}$, определенное в соответствии с
правым нижним пространством в диаграмме $(\ref{15.1})$.
 Чтобы подчеркнуть аналогию со случаем
гладких отображений, далее все рассматриваемые полиэдр (раздутие
особенностей на диагонали и антидиагонали в полиэдрах не
рассматривается) точек самопересечения будет называться
PL-многообразием с особенностями. Эти многообразия с особенностями
имеют особый край, состоящий из точек диагонали и антидиагонали.
Полиэдры без точек диагонали и антидиагонали  обозначаются $ \hat
L_{K_1\circ}^{\downarrow}$, $ \hat L_{K_1\circ}$ соответственно.
Каждый такой полиэдр назовем открытым многообразием с
особенностями. Открытое многообразие с особенностями  $ \hat
L_{K_1\circ}$ естественно вложено в открытое многообразие с краем $\hat
\Gamma_{1\circ}$.

Многообразие с особенностями $\hat L_{K_1}$, в свою очередь, само
служит базой двулистного накрытия $r: L_{K_1} \to \hat L_{K_1}$.
Накрытия определены в соответствии с обозначениями в диаграмме
($\ref{15.1}$). Определено открытое многообразие с особенностями
$L_{K_{1\circ}}$ c краем, которое естественно вкладывается во
взрезанный квадрат $\Gamma_1$. При этом вложении образ
$L_{K_{1\circ}}$ лежит в регулярной окрестности полиэдра $K_{(1)}$
и определена коммутативная диаграмма отображений:

\begin{eqnarray}\label{Q8}
\begin{array}{ccccc}
L_{K_{1\circ}} & \subset & U_{\varepsilon} & \longrightarrow & K(\H,1) \\
\downarrow & &\downarrow && \downarrow \\
\hat L_{K_{1\circ}} &\subset & \hat U_{\varepsilon} & \longrightarrow & K(\G,1), \\
\end{array}
\end{eqnarray}
в которой через $U_{\varepsilon}$, $\hat U_{\varepsilon}$
обозначены регулярные $\varepsilon$-окрестности подполиэдров
$K_{1\circ} \subset \Gamma_1$, $\hat K_{1\circ} \subset \hat
\Gamma_1$ соответственно. Правые горизонтальные стрелки этой
диаграммы являются продолжениями отображений диаграммы
($\ref{Q8}$) на регулярную окрестность $U_{\varepsilon}$ и на
факторпространство $\hat U_{\varepsilon}$ этой окрестности в
многообразиях $\Gamma_1$, $\hat \Gamma_1$ соответственно.

Рассмотрим многообразие с особенностями $L_{K_1}$ и обозначим его
компоненты края через $L_{Q,diag} \cup L_{Q,antidiag}$.
Переобозначим многообразие $L_{K_1}$ через $L_{K_1}^{(0)}$, а его
компоненты края через $L^{(1)}_{Q, diag} \cup L^{(1)}_{Q,
antidiag}$. Ниже в доказательстве Леммы $\ref{lemma30}$ по
аналогии со стратификацией  $(\ref{strat0})$  определена
стратификация
\begin{eqnarray}\label{strat1}
L_{K_1}^{i_{max,1}} \subset \dots \subset L_{K_1}^{0}.
\end{eqnarray}

Обозначим через
\begin{eqnarray}\label{strat11}
W_{K_1}^{(i)} =  L_{K_1}^{i}  \setminus L_{K_1}^{i+1}.
\end{eqnarray}
открытое многообразие с особенностями, содержащее компоненты края.

Обозначим компоненты его края $L^{i}_{Q,diag} \setminus
L^{i+1}_{Q,diag}$, $L^{i}_{Q,antidiag} \setminus
L^{i+1}_{Q,antidiag}$ через $W^{(i)}_{Q,diag}$ и
$W_{Q,antidiag}^{(i)}$ соответственно. Введем аналогичные
обозначения $\hat L_{K_1}^{i} \setminus \hat L^{i+1}_{K_1} = \hat
W_{K_1}^{(i)}$, $\hat L^{i}_{Q,diag} \setminus \hat
L^{i+1}_{Q,diag} = \hat W^{(i)}_{Q,diag}$, $\hat
L^{i}_{Q,antidiag} \setminus \hat L^{i+1}_{Q,antidiag} = \hat
W^{(i)}_{Q,antidiag}$. При вложении $W_{K_1}^{(i)} \subset
U_{\varepsilon} \subset \Gamma_1$ подмногообразие с осбенностями
$W_{Q,diag}^{(i)}$ вкладывается в компоненту края $\Gamma_{diag}$,
а подмногообразие c особенностями $W_{Q,antidiag}^{(i)}$
вкладывается в $\Gamma_{antidiag}$. Аналогичные утверждения
справедливы для компонент границы $\hat W_{\Gamma_1}^{(i)}$.

Обозначим открытое многообразие с особенностями с краем
$L_{K_1}^{(i)} \setminus L^{i+1}_{K_1}$ через $W_{K_1}^{(i)}$.
Обозначим эти компоненты его края через $W^{(i)}_{Q,diag} \cup
W_{Q,antidiag}^{(i)}$. При вложении $W_{K_1}^{(i)} \subset
U_{\varepsilon} \subset \Gamma_1$ подмногообразие
$W_{Q,diag}^{(i)}$ вкладывается в компоненту края $\Gamma_{diag}$,
а подмногообразие $W_{Q,antidiag}^{(i)}$ вкладывается в
$\Gamma_{antidiag}$. Аналогичные утверждения справедливы для
компонент границы $\hat W_{\Gamma_1}^{(i)}$.

Для каждого $i$, $1 \le i \le i_{max,1}$, определена коммутативная
диаграмма ($\ref{25Q}$) отображений полиэдров с выписанными под
диаграммой граничными условиями. В пространства, стоящие в третьей
строке этой диаграммы, отображаются соответственно пространства
центральной строки коммутативной диаграммы ($\ref{26Q}$), которая
также определена для произвольного $i$, $1 \le i \le i_{max,1}$,
соответствующего номеру страта полиэдра точек самопересечения в
стратификации ($\ref{strat1}$):

\begin{eqnarray}\label{25Q}
\begin{array}{ccccccc}
K(\Q_a,1) &&&& K(\Q_a,1) &&\\
&&&&&&\\
\uparrow \hat \phi & \nwarrow \hat \phi &&& \uparrow \phi & \nwarrow \phi &\\
&&&&&&\\
 R \hat K_1^{(i)} &  \supset & R\hat Q_{diag}^{(i)} \cup R\hat Q_{antidiag}^{(i)} &
 & RK_1^{(i)}  & \longleftarrow & RQ_{diag}^{(i)} \cup RQ_{antidiag}^{(i)} \\
&&&&&&\\
 \downarrow p\hat r  &  & \downarrow& & \downarrow pr & & \downarrow  \\
&&&&&&\\
\hat K_1^{(i)} & \supset &   \hat Q_{diag}^{(i)} \cup \hat
Q_{antidiag}^{(i)} & & K_1^{(i)}
& \supset & Q_{diag}^{(i)} \cup Q_{antidiag}^{(i)} \\
&&&&&&\\
\cup &  & \cup && \cup  &  & \cup \\
&&&&&&\\
\hat K_{1\circ}^{(i)} & \supset &  \emptyset & & K_{1\circ}^{(0)}
& \supset & \emptyset \\
&&&&&&\\
\downarrow \hat \zeta & &   & &   \downarrow  \zeta & &   \\
&&&&&&\\
K(\G,1) &  &  &  &     K(\H,1).&& \\
\end{array}
\end{eqnarray}

Граничные условия в диаграмме ($\ref{25Q}$) на пространстве
$Q_{antidiag}^{(i)}$ записываются формулой $\phi  = i_a \circ
\zeta_{Q_{antidiag}^{(i)}} : Q_{antidiag}^{(i)} \longrightarrow
K(\Q_a,1) \stackrel{i_a}{\longrightarrow} K(\H,1)$. Граничные
условия в диаграмме ($\ref{25Q}$) на пространстве $Q_{diag}^{(i)}$
записываются формулой $\phi = i_a \circ p_b \circ
\zeta_{Q_{diag}^{(i)}}: Q_{diag}^{(i)} \longrightarrow K(\H_b,1)
\longrightarrow K(\H_d,1) \longrightarrow K(\Q_a,1)$.

\begin{eqnarray}\label{26Q}
\begin{array}{ccccccc}
K(\Q_a,1) &  &   & & K(\Q_a,1)   & &  \\
&&&&&&\\
\uparrow \hat \lambda_a  &  & \nwarrow \hat \lambda_a & & \uparrow \lambda_a & & \nwarrow \lambda_a \\
&&&&&&\\
\hat W_{K_1}^{(i)} & \supset & \hat W_{Q,diag}^{(i)} \cup \hat
W_{Q,antidiag}^{(i)} &

&W_{K_1}^{(i)} & \supset & W_{Q,diag}^{(i)} \cup W_{Q,antidiag}^{(i)} \\
&&&&&&\\
\cup  &  & \cup & & \cup & & \cup  \\
&&&&&&\\
\hat W_{K_{1\circ}}^{(i)} & \supset & \emptyset &
&W_{K_{1\circ}}^{(i)} & \supset & \emptyset \\
\downarrow \hat \zeta & & & &   \downarrow \zeta & &  \\
&&&&&&\\
K(\G,1) &  &  &  &     K(\H,1)&& \\
\end{array}
\end{eqnarray}

Граничные условия в диаграмме (\ref{26Q}) на пространстве
$W_{Q,antidiag}^{(i)}$ записываются формулой
 $\phi  = i_a \circ
\zeta_{W_{Q,antidiag}^{(i)}} : W_{Q,antidiag}^{(i)}
\longrightarrow K(\Q_a,1) \stackrel{i_a}{\longrightarrow}
K(\H,1)$. Граничные условия в диаграмме (\ref{26Q}) на
пространстве $W_{Q,diag}^{(i)}$ записываются формулой $\phi =
i_{d,a} \circ p_{b,d} \circ \zeta_{W_{Q,diag}^{(i)}}:
W_{Q,diag}^{(i)} \longrightarrow K(\H_b,1) \longrightarrow
K(\H_d,1) \longrightarrow K(\Q_a,1)$.

Воспользуемся следующей леммой, доказательство которой проводится
в конце раздела.

\begin{lemma}\label{lemma31}
Существует малая $PL$-деформация общего положения $(i_{J_1} \circ
\hat p) \mapsto \hat c$ такая, что для каждого кусочно-линейного
страта $W_{K_1}^{(i)}$ полиэдра $L(d)$ в средней строке
коммутативной диаграммы $(\ref{26Q})$ определены отображения $\hat
t^{(i)}: W^{(i)}_{\hat K_1} \to R\hat K_1^{(i)}$, $ t^{(i)}:
W^{(i)}_{K_1} \to RK_1^{(i)}$, называемые отображениями разрешения
особенностей, в соответствующие пространства второй строки
диаграммы $(\ref{25Q})$.
 При этом определены отображения
$\hat \lambda_a^{(i)}: W_{\hat K_1}^{(i)} \to K(\Q_a,1)$,
$\lambda_a^{(i)}: W_{K_1}^{(i)} \to K(\Q_a,1)$ в диаграмме
$(\ref{26Q})$ по формулам $\hat \lambda_a^{(i)} = \hat \phi \circ
\hat t^{(i)}$, $\lambda_a^{(i)} = \phi \circ t^{(i)}$.

Построенные семейства отображений при различных $i$ в
$(\ref{strat1})$ склеиваются до отображений $\hat t: \hat L_{K_1}
\to R\hat K_1$ (соответственно $t: N_{K_1} \to RK_1$) и определены
индуцированные отображения $\hat \lambda_a = \hat \phi \circ \hat
t: \hat L_{K_1} \to K(\Q_a,1)$ (соответственно $\lambda_a = \phi
\circ t: L_{K_1} \to K(\Q_a,1)$). Более того, выполнено
соотношение $r \circ \hat \lambda_a = \lambda_a: L_{K_1} \to \hat
L_{K_1} \to K(\Q_a,1)$ и выполнены граничные условия на компонетах
$L_{Q,diag}$, $L_{Q,antidiag}$ границы стратифицированного
многообразия с особенностями $L_{K_1}$, которые аналогичны
выписаным под формулой $(\ref{26Q})$ граничным условиям для
компонент $W_{Q,diag}^{(i)}$, $W_{Q,antidiag}^{(i)}$
соответственно.
\end{lemma}

\subsubsection*{Окончание доказательства Леммы $\ref{osnlemma2}$}

Доказательство аналогично рассуждениям из окончания доказательства
Леммы $\ref{osnlemma1}$. Лемма $\ref{osnlemma2}$ доказана.


\subsubsection*{Предварительный этап в доказательстве Леммы $\ref{lemma28}$}

Изложим план доказательства. Мы начнем с того, что явно опишем
полиэдры $K_0$, $\hat K_0$ и структурные отображения $\eta$, $\hat
\eta$ на этих полиэдрах при помощи координат. Затем для каждого
$i$, $0 \le i \le i_{max,0}$, построим пространства $RK_0^{(i)}$,
$R\hat K_0^{(i)}$, снабженные отображениями $pr: RK_0^{(i)} \to
K_0^{(i)}$, $p\hat r: R\hat  K_0^{(i)} \to \hat K_0^{(i)}$ и
отображения $\hat \phi: R\hat K_0^{(i)} \to K(\I_a,1)$, $\phi :
RK_0^{(i)} \to K(\I_a,1)$, которые удовлетворяют требуемым
граничным условиям. Затем построим пространства $RK_0$, $R\hat
K_0$, снабженные отображениями $pr: RK_0 \to K_0$, $p\hat r: R\hat
K_0 \to \hat K_0$ и отображения $\hat \phi: R\hat K_0 \to
K(\I_a,1)$, $\phi : RK_0 \to K(\I_a,1)$, которые удовлетворяют
требуемым граничным условиям.

\subsubsection*{Предварительный этап в доказательстве Леммы $\ref{lemma29}$}

Изложим план доказательства. Мы начнем с того, что явно опишем
полиэдры $K_1$, $\hat K_1$ и структурные отображения $\zeta$,
$\hat \zeta$ на этих полиэдрах при помощи координат. Затем для
каждого $i$, $0 \le i \le i_{max,1}$, построим пространства
$RK_1^{(i)}$, $R\hat K_1^{(i)}$, снабженные отображениями $pr:
RK_1^{(i)} \to K_1^{(i)}$, $p\hat r: R\hat  K_1^{(i)} \to \hat
K_1^{(i)}$ и отображения $\hat \phi: R\hat K_1^{(i)} \to
K(\Q_a,1)$, $\phi : RK_1^{(i)} \to K(\Q_a,1)$, которые
удовлетворяют требуемым граничным условиям. Затем построим
пространства $RK_1$, $R\hat K_1$, снабженные отображениями $pr:
RK_1 \to K_1$, $p\hat r: R\hat K_1 \to \hat K_1$ и отображения
$\hat \phi: R\hat K_1 \to K(\Q_a,1)$, $\phi : RK_1 \to K(\Q_a,1)$,
которые удовлетворяют требуемым граничным условиям.

\subsubsection*{Построение стратификаций полиэдров $J_0$, $K_0$, $\hat K_0$
и отображений $K_{(0)} \to K_0$, $\hat K_{(0)} \to \hat K_0$.}

 Упорядочим линзовые пространства, образующие джойн, числами от 1
до $r_0$ и обозначим через $J_0(k_1, \dots, k_s) \subset J_0$
подджойн, образованный выбранным набором окружностей (одномерных
линзовых пространств $S^1/\i$ с номерами $1 \le k_1 < \dots < k_s
\le r$, $0 \ge s \ge r_0$. Построенная стратификация индуцирована из стандартной
стратификации открытых граней стандартного $r_0$-мерного симплекса
$\delta^r$ при естественной проекции $J_0 \to \delta^r$.
Прообразами вершин симплекса служат линзовые пространства $J_0(j)
\subset J_0$, $J_0(j) \approx S^1/\i$, $1 \le j \le r$,
порождающие джойн.

Определим пространство $J_0^{s}$ как подпространство в $J_0$, полученное в результате объединения
всех подпространств $J_0(k_1, \dots, k_s) \subset J_0$. Определим пространство $J_0^{(i)}$, $i_{max,0} \ge i \ge r_0$
по формуле:
$$
J_0^{(i)} = J_0^{i} \setminus J_0^{r_{min,0}+1}.
$$

Обозначим максимальную открытую клетку пространства $\hat
p^{-1}(J_0(k_1, \dots, k_s))$ через $\hat U(k_1, \dots, k_s)
\subset S^{n-k}/\i$. Такую открытую клетку назовем элементарным
стратом глубины $r-s$. Точка на элементарном страте $\hat U(k_1,
\dots, k_s) \subset S^{n-k}/\i$ определяется набором координат
$(\check x_{k_1}, \dots, \check x_{k_s}, l)$, где $\check x_{k_i}$
-- координата на 1-сфере (окружности), накрывающей линзовое
пространство с номером $k_i$, $l$-- координата на соответствующем
$(s-1)$-мерном симплексе джойна. При этом если два набора координат
отождествляются между собой при
преобразовании трансляции циклического $\I_a$-накрытия на
образующую, которое является общим для всего набора координат,
то эти наборы определяют одну и ту же точку на $S^{n-k}/\i$.
Точки на элементарном страте $\hat U(k_1, \dots, k_s)$ лежат в
объединении симплексов c вершинами на линзовых подпространствах
джойна с соответствующими координатами. Каждый элементарный страт
$\hat U(k_1, \dots, k_s)$ служит базой двулистного накрытия
$U(k_1, \dots, k_s) \to \hat U(k_1, \dots, k_s)$, которое
индуцировано из двулистного накрытия $\RP^{n-k} \to S^{n-k}/\i$
при включении $\hat U(k_1, \dots, k_s) \subset S^{n-k}/\i$.

Полиэдр $\hat K_{(0)} \setminus (\hat Q_{0,diag} \cup \hat Q_{0,antydiag})$
разбивается в объединение
открытых подмножеств (элементарных стратов) $K_0(k_1, \dots,
k_s)$, $1 \le s \le r$ в соответствии со стратификацией
пространства $J_0$. Для рассматриваемого страта число $r-s$
недостающих координат до полного набора назовем глубиной страта.
Определено двулистное накрытие $K_0(k_1, \dots,
k_s) \to \hat K_0(k_1, \dots,
k_s)$,
согласованое с построенной стратификацией базы.

Опишем элементарный страт $K_0(k_1, \dots, k_s)$ при помощи системы координат.
Для простоты обозначений разберем случай $s=r_0$
Пусть для пары точек $(x_1$, $x_2)$, определяющих точку на $K(1, \dots, r_0)$,
 фиксирована пара точек $(\check x_1, \check x_2)$ на накрывающей сфере
 $S^{n-k}$, которая переходит в рассматриваемую пару $(x_1$,
 $x_2)$ при проекции $S^{n-k} \to \RP^{n-k}$.
 В соответствии с проделанным выше построением, обозначим через
 $(\check x_{1,i}, \check x_{2,i})$, $i=1, \dots, r_0$ набор сферических координат
 каждой точки. Каждая такая координата определяет точку на
 1-мерной сфере
 (окружности)
$S^1_i$ с тем же номером $i$,
 которая накрывает  соответствующую окружность
$J_0(i) \subset J_0$ в джойне. Заметим, что пара
координат с общим номером необходимо определяет пару точек в одном
слое стандартного циклического $\I_a$-накрытия $S^1 \to S^1/\i$.

Наборы координат  $(\check x_{1,i}, \check x_{2,i})$
рассматриваются с точностью до независимых замен на антиподальные.
Кроме того, точки в паре $(x_1,x_2)$ не допускают естественного
упорядочивания и поднятие точки с $K_0$ до пары точек $(\bar x_1,
\bar x_2)$, лежащих на сфере $S^{n-k}$, определяется с точностью
до 8 различных возможностей. (Порядок группы $\D_4$ равен 8.)

Аналогичное построение справедливо и для точек более глубоких
элементарных стратов $K_0(k_1, \dots, k_s)$, $1 \le s \le r$. При
этом, если наборы координат элементарного страта получаются друг
из друга действием некоторого элемента $\I_a$--накрытия, то
меньший страт лежит целиком на раздутой диагонали (если действие
задается элементом подгруппы $\I_d \subset \I_a$) или
антидиагонали (если действие задается на элемент смежного класса
$\I_a \setminus \I_d$).

Определим полиэдр $\hat K^{(i)}_0$, $0 \le i \le i_{max,0}$ как
дизъюнктное объединение всех элементарных стратов глубины $i$. При
$i \ge 1$ рассматриваются страты, попавшие на диагональ или
антидиагональ. (При $i=0$ диагональный и антидиагональный страт
отсутствует.) Определим полиэдр $\hat K_0$ в результате склейки
объединения всех полиэдров $\hat K^{(i)}_0$ при $0 \le i \le
i_{max,0}$, причем примыкания граней, определяющие склейку,
согласованны с построенными системами координат. По построению
полиэдр $\hat K^{(0)}_0$ совпадает с определенным ранее полиэдром
$\hat K^{(0)}_{(0)}$. При $1 \le i \le i_{max,0}$ полиэдр $\hat
K^{(i)}_0$ отличается от полиэдра  $\hat K^{(i)}_{(0)}$ лишь в
диагональных и антидиагональных подполиэдрах.

Определены естественные отображения
$$
\hat K_{(0)}^{(i)} \to \hat K_{0}^{(i)}, \qquad 1 \le i \le
i_{max,0},
$$

$$
\hat K_{(0)} \to \hat K_{0},
$$
которые являются $PL$--гомеоморфизмом всюду, за исключением,
диагональных и антидиагональных подполиэдров, на которых эти
отображения совпадают с естественными проекциями части раздутой
диагонали (соответственно, антидиагонали) на диагональ
(соответственно антидиагональ), понижающими размерность. Обозначим
через $\hat K_{0 \circ}^{(i)}$, $\hat K_{0 \circ}$  полиэдры,
полученные из $\hat K_{0}^{(i)}$, $\hat K_{0}$ в результате
удаления точек на диагонали и антидиагонали.

Стандартным способом определяются двулистные накрывающие $K_{0}^{(i)}$, $K_{(0)}$, $K_{0 \circ}^{(i)}$, $K_{0 \circ}$
над построенными полиэдрами и отображения
$$
K_{(0)}^{(i)} \to K_{0}^{(i)}, \qquad 1 \le i \le i_{max,0},
$$

$$
K_{(0)} \to K_{0}.
$$

\subsubsection*{Построение стратификаций полиэдров $J_1$, $K_1$, $\hat K_1$
и построение отображений $K_{(1)} \to K_1$, $\hat K_{(1)} \to \hat K_1$.}

 Упорядочим линзовые пространства, образующие джойн, числами от 1
до $r_1$ и обозначим через $J_1(k_1, \dots, k_s) \subset J_1$
подджойн, образованный выбранным набором кватернионных
линзовых пространств $S^3/\Q_a$ с номерами $1 \le k_1 < \dots < k_s
\le r_1$, $0 \ge s \ge r_1$. Построенная стратификация индуцирована из стандартной
стратификации открытых граней стандартного $r_1$-мерного симплекса
$\delta^{r_1}$ при естественной проекции $J_1 \to \delta^{r_1}$.
Прообразами вершин симплекса служат линзовые пространства $J_1(j)
\subset J_1$, $J_1(j) \approx S^3/\Q_a$, $1 \le j \le r_1$,
порождающие джойн.

Определим пространство $J_1^{s}$ как подпространство в $J_1$, полученное в результате объединения
всех подпространств $J_1(k_1, \dots, k_s) \subset J_1$. Определим пространство $J_1^{(i)}$, $i_{max,1} \ge i \ge r_1$
по формуле:
$$
J_1^{(i)} = J_1^{i} \setminus J_1^{r_{min,1}+1}.
$$

Обозначим максимальную открытую клетку пространства $\hat
p^{-1}(J_1(k_1, \dots, k_s))$ через $\hat U(k_1, \dots, k_s)
\subset S^{n-2k}/\Q_a$. Такую открытую клетку назовем элементарным
стратом глубины $r_1-s$. Точка на элементарном страте $\hat U(k_1,
\dots, k_s) \subset S^{n-2k}/\Q_a$ определяется набором координат
$(\check x_{k_1}, \dots, \check x_{k_s}, l)$, где $\check x_{k_i}$
-- координата на 1-сфере (окружности), накрывающей линзовое
пространство с номером $k_i$, $l$-- координата на соответствующем
$(s-1)$-мерном симплексе джойна. При этом если два набора координат
отождествляются между собой при
преобразовании трансляции циклического $\Q_a$-накрытия на
образующую, которое является общим для всего набора координат,
то эти наборы определяют одну и ту же точку на $S^{n-2k}/\Q_a$.
Точки на элементарном страте $\hat U(k_1, \dots, k_s)$ лежат в
объединении симплексов c вершинами на линзовых подпространствах
джойна с соответствующими координатами. Каждый элементарный страт
$\hat U(k_1, \dots, k_s)$ служит базой двулистного накрытия
$U(k_1, \dots, k_s) \to \hat U(k_1, \dots, k_s)$, которое
индуцировано из двулистного накрытия $S^{n-2k}/\i \to S^{n-2k}/\Q_a$
при включении $\hat U(k_1, \dots, k_s) \subset S^{n-2k}/\Q_a$.

Полиэдр $\hat K_{(1)} \setminus (\hat Q_{1,diag} \cup \hat Q_{1,antydiag})$
разбивается в объединение
открытых подмножеств (элементарных стратов) $K_1(k_1, \dots,
k_s)$, $1 \le s \le r$ в соответствии со стратификацией
пространства $J_1$. Для рассматриваемого страта число $r-s$
недостающих координат до полного набора назовем глубиной страта.
Определено двулистное накрытие $K_1(k_1, \dots,
k_s) \to \hat K_1(k_1, \dots,
k_s)$,
согласованое с построенной стратификацией базы.

Опишем элементарный страт $K_1(k_1, \dots, k_s)$ при помощи системы координат.
Для простоты обозначений разберем случай $s=r_1$
Пусть для пары точек $(x_1$, $x_2)$, определяющих точку на $K(1, \dots, r_1)$,
 фиксирована пара точек $(\check x_1, \check x_2)$ на накрывающей сфере
 $S^{n-2k}$, которая переходит в рассматриваемую пару $(x_1$,
 $x_2)$ при проекции $S^{n-2k} \to S^{n-2k}/\i$.
 В соответствии с проделанным выше построением, обозначим через
 $(\check x_{1,i}, \check x_{2,i})$, $i=1, \dots, r_1$ набор сферических координат
 каждой точки. Каждая такая координата определяет точку на
 3-мерной сфере
$S^3_i$ с тем же номером $i$,
 которая накрывает  соответствующее кватернионное линзовое пространство
$J_1(i) \subset J_1$ в джойне. Заметим, что пара
координат с общим номером необходимо определяет пару точек в одном
слое стандартного циклического $\Q_a$-накрытия $S^3 \to S^3/\Q_a$.

Наборы координат  $(\check x_{1,i}, \check x_{2,i})$
рассматриваются с точностью до независимых замен на антиподальные.
Кроме того, точки в паре $(x_1,x_2)$ не допускают естественного
упорядочивания и поднятие точки с $K_1$ до пары точек $(\bar x_1,
\bar x_2)$, лежащих на сфере $S^{n-2k}$, определяется с точностью
до 32 различных возможностей. (Порядок группы $\H$ равен $32$.)

Аналогичное построение справедливо и для точек более глубоких
элементарных стратов $K_0(k_1, \dots, k_s)$, $1 \le s \le r_1$. При
этом, если наборы координат элементарного страта получаются друг
из друга действием некоторого элемента $\Q_a$--накрытия, то
меньший страт лежит целиком на раздутой диагонали (если действие
задается элементом подгруппы $\I_a \subset \Q_a$) или
антидиагонали (если действие задается на элемент смежного класса
$\Q_a \setminus \I_a$).

Определим полиэдр $\hat K^{(i)}_1$, $0 \le i \le i_{max,1}$ как
дизъюнктное объединение всех стратов глубины $i$. При $i \ge 1$
рассматриваются, в том числе, страты, попавшие на диагональ или
антидиагональ. (При $i=0$ диагональный и антидиагональный страт
отсутствует.) Определим полиэдр $\hat K_1$ в результате склейки
объединения всех полиэдров $\hat K^{(i)}_1$ при $0 \le i \le
i_{max,1}$, причем примыкания граней, определяющие склейку,
согласованны с построенными системами координат. По построению
полиэдр $\hat K^{(0)}_1$ совпадает с определенным ранее полиэдром
$\hat K^{(0)}_{(1)}$. При $1 \le i \le i_{max,1}$ полиэдр $\hat
K^{(i)}_1$ отличается от полиэдра  $\hat K^{(i)}_{(1)}$ лишь на
компонентах, лежащих на диагональных и антидиагональных
подполиэдрах.

Определены естественные отображения
$$
\hat K_{(1)}^{(i)} \to \hat K_{1}^{(i)}, \qquad 1 \le i \le
i_{max,1},
$$

$$
\hat K_{(1)} \to \hat K_{1},
$$
которые являются $PL$--гомеоморфизмом всюду, за исключением,
диагональных и антидиагональных подполиэдров, на которых эти
отображения совпадают с естественными проекциями части раздутой
диагонали (соответственно, антидиагонали) на диагональ
(соответственно антидиагональ), понижающими размерность. Обозначим
через $\hat K_{1,\circ}^{(i)}$, $\hat K_{1,\circ}$  полиэдры,
полученные из $\hat K_1^{(i)}$, $\hat K_1$ в результате удаления
точек на диагонали и антидиагонали.

Стандартным способом определяются двулистные накрывающие $K_{1}^{(i)}$, $K_{(1)}$, $K_{1,\circ}^{(i)}$, $K_{1,\circ}$
над построенными полиэдрами и отображения
$$
K_{(1)}^{(i)} \to K_{1}^{(i)},  \qquad 1 \le i \le i_{max,1},
$$

$$
K_{(1)} \to K_{0}.
$$

\subsubsection*{Координатное описание элементарных стратов пространств
$K_{0,\circ}$, $\hat K_{0,\circ}$}

Пусть $x \in K_{0}(1, \dots, r_0)$ -- точка, лежащая на
максимальном элементарном страте. Рассмотрим наборы сферических
координат $\check x_{1,i}$ и $\check x_{2,i}$, $1 \le i \le r_0$
точки $x$. При каждом $i$ возможны следующие случаи: пара $i$-ых
координат совпадает; антиподальна; вторая координата получается из
первой в результате преобразования на образующую (или минус
образующую) циклического накрытия. Сопоставим упорядоченной паре
координат $\check x_{1,i}, \check x_{2,i}$ вычет $v_i$ со
значением $+1$, $-1$, $+\i$ или $-\i$ соответственно.

При изменении набора координат одной из точек на антиподальный
набор, скажем, набора координат точки $x_2$,  набор значений
вычетов при новом выборе прообраза $\bar x_2$ на сферическом
накрытии получается из первоначального набора изменением знаков.
При перенумерации точек набор вычетов изменится на
комплексно-сопряженный. Очевидно, что набор вычетов не изменится,
если выбрать другую точку на том же элементарном страте
пространства $K_{0 \circ}$. Для пространства $\hat K_{0 \circ}$ справедливы
аналогичные построения, т.к. обе точки слоя накрытия $K_{0 \circ} \to \hat
K_{0 \circ}$ лежат на одном элементарном страте.

Элементарные страты пространства $K_0(1, \dots, r)$ в соответствии
с определенным набором вычетов делятся на 3 типа:
$\I_a,\I_b,\I_d$. Если среди вычетов набора встречаются лишь
вычеты $\{+\i,-\i\}$ ($\{+1,-1\}$, будем говорить об элементарном страте типа
$\I_a$ ($\I_b$). Если среди вычетов встречаются как вычеты из
набора $\{+\i,-\i\}$, так и из набора $\{+1,-1\}$, будем говорить об
элементарном страте типа $\I_d$. Легко проверить, что ограничение
характеристического отображения $\eta : K_{0 \circ} \to K(\D_4,1)$ на
элементарный страт типа $\I_a,\I_b,\I_d$ представлено композицией некоторого отображения в
пространство $K(\I_a,1)$, $K(\I_b,1)$, $K(\I_d,1)$ с отображением
$i_a: K(\I_a,1) \to K(\D_4,1)$, $i_b: K(\I_b,1) \to K(\D_4,1)$
$i_d: K(\I_d,1) \to K(\D_4,1)$ соответственно. Для стратов первых
двух типов редукция структурного отображения (с точностью до
гомотопии) неоднозначна, но определена лишь с точностью до
отображения, соответствующему гомоморфизму внешнего сопряжения в подгруппах
$\I_a$, $\I_b$.

Каждый элементарный страт пространства $K_{0 \circ}(1, \dots, r_0)$
двулистно накрывает  элементарный страт пространства $\hat K_{0,\circ}(1,
\dots, r_0)$. При этом элементарные страты пространства $\hat K_{0,\circ}(1,
\dots, r_0)$ делятся на типы $\E_a$, $\E_b$, $\E_c$ в соответствии с
типом $\I_a, \I_b, \I_d$ накрывающей компоненты. Эта классификация
соответствует редукции структурного отображения в пространство
$K(\E,1)$ к отображению в подпространство $K(\E_a,1)$,
$K(\E_b,1)$, $K(\E_d,1)$ соответственно. Для элементарных стратов
большей глубины построения аналогичны.

\subsubsection*{Координатное описание элементарных стратов пространств $K_{1,\circ}$, $\hat K_{1,\circ}$}

Для каждой точки $x \in K_{1}(1, \dots, r_1)$, лежащей на
максимальном элементарном страте рассмотрим наборы сферических
координат $\check x_{1,i}$ и $\check x_{2,i}$, $1 \le i \le r_1$.
При каждом $i$ возможны восемь случаев, в зависимости от того,
какой из элементов $v_i = (\check x_{2,i})^{-1} \check x_{1,i}$
структурной группы $\Q_a$ переводит при умножении справа
сферическую координату $\check x_{1,i}$ в координау $\check
x_{2,i}$. Сопоставим упорядоченной паре координат $\check x_{1,i},
\check x_{2,i}$ вычет $v_i$ со значением $+1$, $-1$, $+\i$, $-\i$,
$+\j$, $-\j$, $+\k$, $-\k$ соответственно.

При изменении набора координат точки $x_1$ (соответственно
$x_2$) на преобразование трансляции (правое умножение) на элемент $t$ из группы
$\Q_a$ набор вычетов при новом выборе прообраза $\bar x_2$ на
сферическом накрытии получается из первоначального набора вычетов
умножением справа на $t$ (соответственно умножением слева на $t^{-1}$).
При перенумерации точек набор вычетов
изменится на обратный (т.е. вычет $v_i$ изменится на $v_i^{-1}$). Набор вычетов
не изменится при изменении точки $x$ при преобразовании трансляции в накрывающем при накрытии
$\hat K_{1}(1, \dots, r_1) \to K_{1}(1, \dots, r_1)$.
Поэтому для элементарных стратов пространства $\hat K_{\circ,1}$ справедливы
аналогичные построения.

Элементарные страты пространства $K^{(0)}_1(1, \dots, r_1)$ в
соответствии с построенным выше набором вычетов делятся на 3 типа:
$\Q_a,\E_b,\I_a$. Если среди вычетов набора встречаются лишь
вычеты $\{+\j,-\j,+\k,-\k\}$ из класса смежности $\Q_a \setminus
\I_a$ (соответственно $\{+\i, -\i, +1, -1 \}$) из подгруппы
$\I_a$), причем только вычеты $+\j$ $-\j$, либо только вычеты
$+\k$ $-\k$, (соответственно  только вычеты $+1$, $-1$, либо
только вычеты $+\i$, $-\i$) будем говорить о страте типа $\Q_a$
(соответвтвенно $\H_b$). В другом случае, например, если среди
вычетов встречаются как вычеты из подгруппы $\I_a$, так и из
класса смежности $\Q_a \setminus \I_a$, будем говорить о
компоненте типа $\I_a$.

Легко проверить, что ограничение характеристического отображения
$\zeta : K_1 \to K(\H,1)$ на элементарный страт типа $\Q_a$,
($\E_b$)  представлено композицией некоторого отображения в
пространство $K(\Q_a,1)$ ($K(\E_b,1)$)  с отображением включения
$i_a: K(\Q_a,1) \to K(\H,1)$ ($i_b: K(\E_b,1) \to K(\H,1)$). Для
стратов этих двух подтипов редукция структурного отображения (с
точностью до гомотопии) неоднозначна, но лишь с точностью до
отображения, соответствующего внешнему сопряжению в подгруппах
$\Q_a$ ($\E_b$) некоторым произвольным элементом из класса
смежности $\H \setminus \Q_a$ ($\H \setminus \E_b$). Для
элементарных  типа $\I_a$ ограничение отображение $\zeta :
K_{1,\circ} \to K(\H,1)$
 представлено
композицией некоторого отображения в пространство $K(\I_a,1)$ с
отображением включения $i_{\I_a}: K(\I_a,1) \to K(\H,1)$.

Каждый элементарный страт пространства $K_{1 \circ}(1, \dots, r)$
двулистно накрывает  элементарный страт пространства $\hat K_{1 \circ}(1,
\dots, r)$. При этом элементарные страты пространства $\hat K_{1 \circ}(1,
\dots, r)$ делятся на типы $\G_a$, $\G_b$, $\G_d$ в
соответствии с группой симметрии компоненты (явное описание этих структурных групп не потребуется). Рассматриваемая классификация соответствует редукции структурного
отображения в пространство $K(\G,1)$ к отображению в
подпространство $K(\G_a,1)$, $K(\G_b,1)$, $K(\G_d,1)$
соответственно. Для элементарных стратов большей глубины все
построения аналогичны.

\subsubsection*{Построение структурных отображений $\eta: K_{0\circ} \to K(\D_4,1)$,
$\hat \eta: \hat K_{0 \circ} \to K(\E,1)$ при помощи координат}

Пусть $x = [(x_1, x_2)]$ -- отмеченная точка на $K_{0,\circ}$,
лежащая на максимальном элементарном страте. Рассмотрим замкнутый
путь $\lambda: S^1 \to K_0$, начинающийся и заканчивающийся в этой
отмеченной точке, пересекающий особые страты глубины 1 в общем
положении по конечному множеству точек. Пусть $(\check x_1, \check
x_2)$-- выбор сферических прообразов точки $x$. Определим другую
пару $(\check x'_1,\check x'_2)$ сферических прообразов точки $x$,
которая будет называться системой координат, полученной в
результате естественного преобразования системы координат $(\check
x_1, \check x_2)$ вдоль пути $\lambda$.

В регулярных точках пути $\lambda$  семейство пар сферических
прообразов в однопараметрическом семействе меняется непрерывно,
что однозначно определяет прообразы в конечной точке пути по
начальным значениям. При пересечении пути со стратом глубины 1
соответствующая пара сферических координат c номером $l$
претерпевает разрыв. Поскольку все остальные координаты остаются
регулярными, их продолжение вдоль пути в критический момент
времени однозначно определено. Но для заданной точки $x$ на
элементарном страте глубины 0 пространства $K_0$ выбор хотябы
одной пары сферических координат с некотором номером однозначно
определяет выбор сферических координат со всеми остальными
номерами. Следовательно, продолжение сферических координат вдоль
пути однозначно определено в окрестности особой точки пути.

Преобразование упорядоченной пары $(\check x_1, \check x_2)$ в
упорядоченную пару $(\check x'_1,\check x'_2)$ определяет элемент
в группе $\D_4$. Этот элемент не зависит от выбора пути $l$ в классе эквивалентных путей, определенных
отношением гомотопности в группе $\pi_1(K_{0 \circ},x)$. Тем самым, определен гомоморфизм
$\pi_1(K_{0 \circ},x) \to \D_4$, причем индуцированное отображение
\begin{eqnarray}\label{eta}
\eta: K_{0 \circ} \to K(\D_4,1)
\end{eqnarray}
совпадает со структурным отображением, которое было определено
ранее. Структурное отображение
\begin{eqnarray}\label{hateta}
\hat \eta: K_{0 \circ} \to K(\E,1)
\end{eqnarray}
определено аналогичной конструкцией. Нетрудно проверить, что
ограничение структурного отображения $\eta$ ($\hat \eta$) на
компоненты связности одного элементарного страта $K_{0 \circ}(1,
\dots, r)$ ($\hat K_{0 \circ}(1, \dots, r)$) гомотопно отображению
в подпространства $K(\I_a,1)$, $K(\I_b,1)$, $K(\I_d,1)$
($K(\E_a,1)$, $K(\E_b,1)$, $K(\E_d,1)$), что соответствует типу и
подтипу элементарного страта.

Рассмотрим элементарный страт $ K_{0}(k_1, \dots, k_s) \subset
K^{(r_0-s)}_{0 \circ}$ (соответственно $\hat K_0(k_1, \dots, k_s)
\subset \hat K^{(r_0-s)}_{0 \circ}$) глубины $r_0-s$. Обозначим
через
\begin{eqnarray}\label{pi}
\pi: K_{0}(k_1, \dots, k_s) \to K(\Z/2,1)
\end{eqnarray}
(соответственно через
\begin{eqnarray}\label{hatpi}
\hat \pi: \hat K_{0}(k_1, \dots, k_s) \to K(\Z/2,1))
\end{eqnarray}
 -- классифицирующее отображение, отвечающее за перестановку
точек пары при обходе по замкнутому пути на этом элементарном страте. Отображение $\pi$ (соответственно $\hat \pi$) совпадает с композицией
$$ K_{0}(k_1, \dots, k_s) \stackrel{\eta}{\longrightarrow} K(\D_4,1) \stackrel{p}{\longrightarrow} K(\Z/2,1), $$
(соответственно с композицией
$$ \hat K_{0}(k_1, \dots, k_s) \stackrel{\hat \eta}{\longrightarrow} K(\E,1)
\stackrel{\hat p}{\longrightarrow} K(\Z/2,1),) $$ где $K(\D_4,1)
\stackrel{p}{\longrightarrow} K(\Z/2,1)$ (соответственно $K(\E,1)
\stackrel{\hat p}{\longrightarrow} K(\Z/2,1)$) -- отображение
классифицирующих пространств, которое индуцировано эпиморфизмом
$\D_4 \to \Z/2$ с ядром $\I_c \subset \D_4$ (cоответственно
эпиморфизмом $\E \to \Z/2$ с ядром $\E_c \subset \E$).

\begin{lemma}\label{lemma32}
Ограничение отображения $(\ref{pi})$ $($ соответственно
$(\ref{hatpi})$$)$ на каждый элементарный страт гомотопно
композиции
\begin{eqnarray}\label{pi1}
\pi_1: K_{0}(k_1, \dots, k_s) \to S^1 \subset K(\Z/2,1),
\end{eqnarray}
$($ соответственно композиции
\begin{eqnarray}\label{hatpi1}
\hat \pi_1: \hat K_{0}(k_1, \dots, k_s) \to S^1 \subset
K(\Z/2,1),)
\end{eqnarray}
где $S^1 \subset K(\Z/2,1)$-- вложение 1-мерного остова
классифицирующего пространства, $\pi_1$ $($ cоответственно $\hat
\pi_1$ $)$--некоторое отображение.

\end{lemma}

\subsubsection*{Доказательство Леммы $\ref{lemma32}$}
Структурное отображение $\eta$ определяется преобразованием
координат на стратах типов $\I_a$, $\I_b$, $\I_d$. Отображение
$(\ref{pi})$ определяется преобразованием выбранной (произвольной)
пары координат с вычетами $\i$, $-\i$ для страта $\I_a$ и с
вычетами $+1$, $-1$ для страта $\I_b$. Для страта типа $\I_d$
отображение $(\ref{pi})$ нулевое. Преобразование выбранной
координаты на страте данного типа определяется отображением вида
$(\ref{pi1})$, что доказывает требуемое в случае отображения
$(\ref{pi})$. Для отображения $(\ref{hatpi1})$ доказательство
аналогично. Лемма $\ref{lemma32}$ доказана.
\[  \]

Нам также потребуется лемма  о 3-точечном конфигурационном
пространстве, ассоциированном с накрытием $p: U(k_1, \dots, k_s)
\to J_0(k_1, \dots, k_s)$ над элементарным стратом полиэдра $J_0$.
Сформулируем требуемые определения (см. ниже Лемма
$\ref{lemma33}$).

Для накрытия $p: U(k_1, \dots, k_s) \to J_0(k_1, \dots, k_s)$, $1
\le k_1 < \dots < k_s \le r_0$, $r_{0,min} \le s$ над произвольным элементарным стратом
пространства $J_0$ определим конфигурационное пространство $\bar K3_0(k_1, \dots,
k_s)$ по формуле:
$$
\bar K3_0(k_1, \dots, k_s) =
$$
\begin{eqnarray}\label{K3}
\{(x,y,z) \in U^3(k_1, \dots k_s): p(x)=p(y)=p(z), x \ne y, y \ne
z, z \ne x\}.
\end{eqnarray}
На пространстве $U^3(k_1, \dots k_s)$ действует группа $\Sigma_3$:
\begin{eqnarray}\label{sigma}
\sigma: \Sigma_3 \times U^3(k_1, \dots, k_s) \to U^3(k_1, \dots,
k_s).
\end{eqnarray}
Ограничение действия $(\ref{sigma})$ на подпространство
$(\ref{K3})$  является свободным действием, которое мы обозначим
через
\begin{eqnarray}\label{sigma1}
\sigma_{K3_0}: \Sigma_3 \times \bar K3_0(k_1, \dots, k_s) \to \bar
K3_0(k_1, \dots, k_s).
\end{eqnarray}
Определено факторпространство по действию $(\ref{sigma1})$,
которое мы обозначим через $K3_0(k_1, \dots, k_s)$.

На пространстве $\bar K3_0(k_1, \dots, k_s)$ определена свободная
инволюция, индуцированная инволюцией  в накрывающем накрытия
$U(k_1, \dots, k_s) \to \hat U(k_1, \dots, k_s)$, обозначаемая
через $\sim$.

Определено свободное действие на факторпространствах
соответствующих инволюций
\begin{eqnarray}\label{sigma2}
\hat \sigma_{K3_0}: \Sigma_3 \times \bar K3_0(k_1, \dots, k_s)/\sim
\to \bar K3_0(k_1, \dots, k_s)/ \sim.
\end{eqnarray}
Определено пространство $\hat K3_0(k_1, \dots k_s)$, которое
является базой 2-листного накрытия $K3_0(k_1, \dots, k_s) \to \hat
K3_0(k_1, \dots k_s)$, определенное как факторпространство по
действию $(\ref{sigma2})$.

Определена коммутативная диаграмма:
\begin{eqnarray}\label{K33}
\begin{array}{ccc}
K3_0(k_1,\dots,k_s)& \stackrel{\sigma_0}{\longrightarrow} & K(\Sigma_3,1)\\
\downarrow & & \| \\
\hat K3_0(k_1, \dots, k_s) & \stackrel{\hat
\sigma_0}{\longrightarrow} & K(\Sigma_3,1),
\end{array}
\end{eqnarray}
в которой через $\sigma_0$ (соответственно $\hat \sigma_0$)
обозначены характеристические отображения накрытий
$(\ref{sigma1})$ (соответственно $(\ref{sigma2})$).

\begin{lemma}\label{lemma33}
Отображения $\sigma_0$, $\hat \sigma_0$ в диаграмме $(\ref{K33})$
тривиальны.
\end{lemma}

\subsubsection*{Доказательство Леммы $\ref{lemma33}$}
Поскольку группа накрытия $U(k_1, \dots, k_s) \to J_0(k_1, \dots,
k_s)$
не содержит элементов нечетного порядка, то группа накрытия
$(\ref{sigma1})$
может содержать лишь нетривиальные элементы порядка степеней $2$. Пусть $l:
S^1 \to \bar K3_0(k_1, \dots, k_s)$-- замкнутый путь, вдоль
которого точка с координатами $(\check{x},\check{y},\check{z}) \in
\bar K3_0(k_1, \dots, k_s)$ преобразуется в точку с координатами
$(\check{x}_1,\check{y}_1, \check{z}') \in \bar K3_0(k_1, \dots,
k_s)$ в том же слое накрытия $(\ref{sigma1})$. Тогда по условию
координаты $\check{z}$ и $\check{z}'$ либо совпадают, либо связаны
преобразованием на образующую из $\I_d$. Но тогда и пары координат
($\check{x},\check{x}_1$), ($\check{y},\check{y}_1$) также связаны
таким же преобразованием, поскольку при переносе вдоль пути каждая из трех соответствующих
координат точек слоя накрытия $(\ref{sigma1})$ меняется
одинаково. Следовательно, три пары точек с соответствующими
координатами совпадают. Доказано, что отображение $\sigma_0$
тривиально. Доказательство для отображения $\hat \sigma_0$
аналогично. Лемма $\ref{lemma33}$ доказана.
\[  \]

\subsubsection*{Построение структурных отображений $\zeta: K_{1\circ} \to K(\H,1)$,
$\hat \eta: \hat K_{1\circ} \to K(\G,1)$ при помощи координат}

Пусть $x = [(x_1, x_2)]$ -- отмеченная точка на $K_{1 \circ}$,
лежащая на максимальном элементарном страте. Рассмотрим замкнутый
путь $\lambda: S^1 \to K_{1 \circ}$, начинающийся и
заканчивающийся в этой отмеченной точке, пересекающий особые
страты глубины 1 в общем положении по конечному множеству точек.
Пусть $(\check x_1, \check x_2)$-- выбор сферических прообразов
точки $x$. Определим другую пару $(\check x'_1,\check x'_2)$
сферических прообразов точки $x$, которая будет называться
системой координат, полученной в результате естественного
преобразования системы координат $(\check x_1, \check x_2)$ вдоль
пути $\lambda$.

В регулярных точках пути $\lambda$  семейство пар сферических
прообразов в однопараметрическом семействе меняется непрерывно,
что однозначно определяет прообразы в конечной точке пути по
начальным значениям. При пересечении пути со стратом глубины 1
соответствующая пара сферических координат претерпевает разрыв.
Поскольку все остальные координаты остаются регулярными, их
продолжение вдоль пути в критический момент времени однозначно
определено. Но для заданной точки $x$ на элементарном страте
глубины 0 пространства $K_{1 \circ}$ выбор хотябы одной пары сферических
координат с некотором номером однозначно определяет выбор
сферических координат со всеми остальными номерами. Следовательно,
продолжение сферических координат вдоль пути однозначно определено
в окрестности особой точки пути.

Преобразование упорядоченной пары $(\check x_1, \check x_2)$ в
упорядоченную пару $(\check x'_1,\check x'_2)$ определяет элемент
в группе $\H$. Этот элемент не зависит от выбора пути $\lambda$ в
группе $\pi_1(K_{1 \circ},x)$. Тем самым, определен гомоморфизм
$\pi_1(K_{1 \circ},x) \to \H$. Индуцированное отображение
\begin{eqnarray}\label{zeta}
\zeta: K_{1 \circ} \to K(\H,1)
\end{eqnarray}
совпадает со структурным отображением. Структурное отображение
\begin{eqnarray}\label{hatzeta}
\hat \zeta: K_{1 \circ} \to K(\G,1)
\end{eqnarray}
 определено аналогичной
конструкцией. Нетрудно проверить, что ограничение структурного
отображения $\zeta$ (соответственно $\hat \zeta$) на один
элементарный страт $K_{1}(1, \dots, r)$ (соответственно $\hat
K_{1}(1, \dots, r)$) гомотопно отображению в подпространства
$K(\Q_a,1)$, $K(\E_b,1)$, $K(\I_a,1)$ (в классифицирующие
пространства, соответствующие квадратичным расширениям групп
$\Q_a$, $\E_b$, $\I_a$).

Рассмотрим элементарный страт $K_1(k_1, \dots, k_s)$
(соответственно $\hat K_1(k_1, \dots, k_s)$) глубины $r_1-s$
пространства $K_{1 \circ}$ (соответственно $\hat K_{1 \circ}$). Обозначим через
\begin{eqnarray}\label{piQ}
\pi: K_{1}(k_1, \dots, k_s) \to K(\Z/2,1)
\end{eqnarray}
(соответственно через
\begin{eqnarray}\label{hatpiQ}
\hat \pi: \hat K_1 (k_1, \dots, k_s) \to K(\Z/2,1))
\end{eqnarray}
 -- классифицирующее отображение, отвечающее за перестановку
точек пары при обходе по замкнутому пути на элементарном страте. Отображение $\pi$
(соответственно $\hat \pi$) совпадает с композицией
$$ K_1(k_1, \dots, k_s) \stackrel{\zeta}{\longrightarrow} K(\H,1) \stackrel{p}{\longrightarrow} K(\Z/2,1), $$
(соответственно с
$$ \hat K_1(k_1, \dots, k_s) \stackrel{\hat \zeta}{\longrightarrow} K(\G,1) \stackrel{\hat p}
{\longrightarrow} K(\Z/2,1),$$ где $K(\H,1)
\stackrel{p}{\longrightarrow} K(\Z/2,1)$ (соответственно $K(\G,1)
\stackrel{\hat p}{\longrightarrow} K(\Z/2,1)$) -- отображение
классифицирующих пространств, которое индуцировано эпиморфизмом
$\H \to \Z/2$ с ядром $\H_c \subset \H$ (cоответственно
эпиморфизмом $\G \to \Z/2$ с ядром $\G_c \subset \G$).

\begin{lemma}\label{lemma34}
Отображение $(\ref{piQ})$ $($ соответственно $(\ref{hatpiQ})$ $)$
гомотопно композиции
\begin{eqnarray}\label{pi1Q}
\pi_1: K_1(k_1, \dots, k_s) \to \RP^2 \subset K(\Z/2,1),
\end{eqnarray}
$($ соответственно
\begin{eqnarray}\label{hatpi1Q}
\hat \pi_1: \hat K_1(k_1, \dots, k_s) \to \RP^2 \subset
K(\Z/2,1),)
\end{eqnarray}
где $\RP^2 \subset K(\Z/2,1)$-- вложение 2-мерного остова
классифицирующего пространства, $\pi_1$ $($ cоответственно $\hat
\pi_1$ $)$--некоторое отображение.
\end{lemma}

\subsubsection*{Доказательство Леммы $\ref{lemma34}$}
Структурное отображение $\zeta$ определяется преобразованием
координат на стратах типов $\Q_a$, $\H_b$, $\I_a$. Отображение
$(\ref{pi})$ определяется преобразованием выбранной (произвольной)
пары координат с вычетами $\i$, $-\i$ для страта $\Q_a$ и с
вычетами $+1$, $-1$ для страта $\H_b$. Для страта типа $\I_a$
отображение $(\ref{piQ})$ нулевое. Преобразование выбранной
координаты на страте данного типа определяется отображением вида
$(\ref{pi1Q})$, поскольку как отображение $S^3/\Q_a \to
K(\Z/2,1)$, соответствующее эпиморфизму $\Q_a \to Z/2$ с ядром
$\I_a \subset \Q_a$, так и отображение $S^3/\i \to K(\Z/2,1)$,
соответствующее эпиморфизму $\I_a \to Z/2$ с ядром $\I_d \subset
\I_a$, гомотопны отображениям в 2-мерный остов классифицирующего
пространства. Это доказывает требуемое в случае отображения
$(\ref{piQ})$. Для отображения $\ref{hatpi1Q}$ доказательство
аналогично. Лемма $\ref{lemma34}$ доказана.
\[  \]

Нам также потребуется лемма, аналогичная Лемме $\ref{lemma33}$.
Сформулируем требуемые определения.

Для накрытия $p: U(k_1, \dots, k_s) \to J_1(k_1, \dots, k_s)$, $1
\le k_1 < \dots < k_s \le r_1$, $r_{1,min} \le s$ над произвольным элементарным стратом
определим конфигурационное пространство $\bar K3_1(k_1, \dots,
k_s)$ по формуле:
$$
\bar K3_1(k_1, \dots, k_s) =
$$
\begin{eqnarray}\label{K3Q}
\{(x,y,z) \in U^3(k_1, \dots k_s): p(x)=p(y)=p(z), x \ne y, y \ne
z, z \ne x\}.
\end{eqnarray}
На пространстве $U^3(k_1, \dots k_s)$ действует группа $\Sigma_3$,
переставляя координаты:
\begin{eqnarray}\label{sigmaQ}
\sigma: \Sigma_3 \times U^3(k_1, \dots, k_s) \to U^3(k_1, \dots,
k_s).
\end{eqnarray}
Ограничение действия $(\ref{sigmaQ})$ на подпространство
$(\ref{K3Q})$  является свободным действием, которое мы обозначим
через
\begin{eqnarray}\label{sigma1Q}
\sigma_{K3_1}: \Sigma_3 \times \bar K3_1(k_1, \dots, k_s) \to \bar
K3_1(k_1, \dots, k_s).
\end{eqnarray}
Определено факторпространство по действию $(\ref{sigma1Q})$,
которое мы обозначим через $K3_1(k_1, \dots, k_s)$.

На пространстве $\bar K3_1(k_1, \dots, k_s)$ определена свободная
инволюция, индуцированная инволюцией  в накрывающем накрытия
$U(k_1, \dots, k_s) \to \hat U(k_1, \dots, k_s)$, обозначаемая ниже снова
через $\sim$. Определено свободное действие на факторпространствах
соответствующих инволюций
\begin{eqnarray}\label{sigm2Q}
\hat \sigma_{K3_1}: \Sigma_3 \times \bar K3_1(k_1, \dots, k_s)/\sim
\to \bar K3_1(k_1, \dots, k_s)/\sim.
\end{eqnarray}
Определено пространство $\hat K3_1(k_1, \dots k_s)$, которое
является базой 2-листного накрытия $K3_1(k_1, \dots, k_s) \to \hat
K3_1(k_1, \dots k_s)$, определенное как факторпространство по
действию $(\ref{sigm2Q})$.

Определена коммутативная диаграмма:
\begin{eqnarray}\label{K33Q}
\begin{array}{ccc}
K3_1(k_1,\dots,k_s)& \stackrel{\sigma_1}{\longrightarrow} & K(\Sigma_3,1)\\
\downarrow & & \| \\
\hat K3_1(k_1, \dots, k_s) & \stackrel{\hat
\sigma_1}{\longrightarrow} & K(\Sigma_3,1),
\end{array}
\end{eqnarray}
в которой через $\sigma_1$ (соответственно $\hat \sigma_1$)
обозначены характеристические отображения накрытий
$(\ref{sigma1Q})$ (соответственно $(\ref{sigm2Q})$).

\begin{lemma}\label{lemma33Q}
Отображения $\sigma_1$, $\hat \sigma_1$ в диаграмме $(\ref{K33Q})$
тривиальны.
\end{lemma}

\subsubsection*{Доказательство Леммы $\ref{lemma33Q}$}
Доказательство аналогично доказательству Леммы $\ref{lemma33}$.
\[  \]

\subsubsection*{Подпространства $K3_{0 \circ} \subset K3_{0} $, $\hat K3_{0 \circ} \subset \hat K3_{0}$}

Определим пространство $K3_0$ по формуле
$$ [(x,y,z)], (x,y,z) \in (\RP^{n-k})^3 \vert p_{J_0}(x) = p_{J_0}(y) = p_{J_0}(z) \in J_0^{(0)} \setminus J_0^{(r_{min})}, $$
(где $[\quad ]$ означает класс эквивалентности троек при изменении
порядка точек) причем если $x=y=z$, то $p_{J_0}(x)$ не принадлежит
страту пространства $J_0$ максимальной размерности. Определим
подпрстранство $K3_{0 \circ} \subset K3_0$ как простроанство,
полученное из $K3_0$ в результате удаления неупорядоченных троек
точек, попавших на двойную или тройную диагональ, т.е. троек
$[(x,y,z)]$, для которых выполнено хотябы одно из равенств $x=y$,
$y=z$, $z=x$. Подпространство $K3_{0 \circ} \subset K3_0$ является
двулистно накрывающим пространством над подпространством $\hat
K3_{0 \circ} \subset K3_0$. Определения этих накрытий стандартно и
опускается.

Пространства $K3_{0 \circ}$, $\hat K3_{0 \circ}$ можно определить по-другому при помощи склеек элементарных стратов различной
глубины $s$, $0 \le s \le i_{max,0}$. Пространство $K3_{0 \circ}$ (соответственно $\hat K3_{0 \circ}$) склеивается из элементарных стратов  $K3_1(k_1, \dots, k_s)$
(соответственно $\hat K3_1(k_1, \dots, k_s)$) глубины $i$,
$0 \le i \le r_{max,0}$. При этом используются элементарные страты, состоящие из неупорядоченных троек попарно различных точек.

\subsubsection*{Подпространства $K3_{1 \circ} \subset K3_1$, $\hat K3_{1 \circ} \subset \hat K3_1$}
Построение пространств $K3_{1 \circ}$, $\hat K3_{1 \circ}$
полностью аналогично построению пространств $K3_{0 \circ}$, $\hat
K3_{0 \circ}$.

\subsubsection*{Построение пространств $RK_{0\circ}$, $R\hat K_{0\circ}$, $RK_{0}$, $R\hat K_{0}$
для разрешения  особенностей и построение
 3-листных накрытий $p3_{\circ}: RK_{0\circ} \to
K3_{0\circ}$, $\hat p3_{\circ}: R\hat K_{0\circ} \to \hat
K3_{0\circ}$;  $p3: RK_{0} \to K3_{0}$, $\hat p3: R\hat K_{0} \to
\hat K3_{0}$.}

Рассмотрим подполиэдр $K_{0\circ}^{(s)} \subset K_0$, определенный
в результате объединения всех особых стратов глубины $s$,
$0 \le s \le r_{min,0}$, не лежащих целиком на диагонали
или антидиагонали.

Рассмотрим диаграмму ($\ref{14}$) отображений классифицирующих
пространств. Эта диаграмма индуцирует следующую
коммутативную диаграмму:

\begin{eqnarray}\label{pK3s}
\begin{array}{ccc}
RK^{(s)}_{0\circ} & \longrightarrow & R\hat K^{(s)}_{0\circ} \\
&&\\
\downarrow p3^{(s)}_{\circ}& & \downarrow \hat p3^{(s)}_{\circ}\\
&&\\
K3^{(s)}_{0\circ} & \longrightarrow & \hat K3^{(s)}_{0\circ}.\\
\end{array}
\end{eqnarray}

Пространство $RK^{(s)}_{0\circ}$ (соответственно $R\hat
K^{(s)}_{0\circ}$) в этой диаграмме определим как пространство
упорядоченных троек неупорядоченных пар $\{[x,y],[y,z],[z,x]\}$,
$x,y,z \in p^{-1}(J_0^{(s)}) \subset \RP^{n-k}$, $x \ne y, y \ne
z, z \ne x$, $p(x)=p(y)=p(z)$ (соответственно как
факторпространство указанных троек при умножении на $\i$).
Пространство $RK^{(s)}_{0\circ}$ (соответственно $R\hat
K^{(s)}_{0\circ}$) является накрывающим пространством 3-листного
накрытия $p3^{(s)}_{\circ}$ (соответственно $\hat
p3^{(s)}_{\circ}$) над пространством $K3^{(s)}_{0\circ}$
(соответственно $\hat K3^{(s)}_{\circ}$). Накрытие
$p3^{(s)}_{\circ}$ (соответственно $\hat p3^{(s)}_{\circ}$)  над
каждым элементарным стратом $K3_{0}(k_1, \dots k_s)$
(соответственно $K3_{0}(k_1, \dots k_s)$) определено по формуле
$\{[x,y],[y,z],[z,x]\} \to [x,y,z]$.  Пространства
$RK^{(s)}_{0\circ}$ $R\hat K^{(s)}_{0\circ}$ включены в
коммутативную диаграмму отображений (накрытий):

\begin{eqnarray}\label{piK3s}
\begin{array}{ccc}
RK^{(s)}_{0\circ} & \longrightarrow & R\hat K^{(s)}_{0\circ} \\
&&\\
\downarrow \pi3^{(s)}_{\circ}& & \downarrow \hat \pi3^{(s)}_{\circ}\\
&&\\
K^{(s)}_{0\circ} & \longrightarrow & \hat K^{(s)}_{0\circ}.\\
\end{array}
\end{eqnarray}

В этой диаграмме отображение (накрытие) $\pi3^{(s)}_{\circ}$
(соответственно $\hat\pi3^{(s)}_{\circ}$) является тавтологическим накрытием,
переводящим первую неупорядоченную пару точек $[x,y]$,
$x,y \in p^{-1}(J_{0}^{(s)}) \subset \RP^{n-k}$, $p(x)=p(y)$
в ту же самую неупорядоченную пару точек (соответственно в класс
эквивалентности той же пары при действии $\i$), которая рассматривается
как точка на $K^{(s)}_{0\circ}$ (соответственно на  $\hat K^{(s)}_{\circ}$).

\begin{lemma}\label{etaK30}

Композиция $\eta^{(s)} \circ \pi3^{(s)}_{\circ}: RK^{(s)}_{0\circ}
\longrightarrow K^{(s)}_{0\circ} \longrightarrow K(\D_4,1)$, где
$\eta^{(s)}$--ограничение структурного отображения $\eta$ $($ см.
$(\ref{eta})$ $)$ на соответствующий элементарный страт глубины
$s$, гомотопна отображению $ \i_{\I_d,\D_4} \circ
\eta_{RK^{(s)}_{0\circ}}: RK^{(s)}_{0\circ} \longrightarrow
K(\I_d, 1)  \longrightarrow K(\D_4,1)$.

Композиция $\hat \eta^{(s)} \circ \pi3^{(s)}_{\circ}: R \hat
K^{(s)}_{0\circ} \longrightarrow \hat K^{(s)}_{0\circ}
\longrightarrow K(\E,1)$, где $\hat \eta^{(s)}$--ограничение
структурного отображения $\hat \eta$ $($ см. $(\ref{hateta})$ $)$
на соответствующий элементарный страт глубины $s$, гомотопна
отображению $ \i_{\I_a,\E} \circ \eta_{R \hat K^{(s)}_{0\circ}}: R
\hat K^{(s)}_{0\circ} \longrightarrow K(\I_a, 1) \longrightarrow
K(\E,1)$.

\end{lemma}

\subsubsection*{Доказательство Леммы $\ref{etaK30}$}

 Ограничение
структурного отображения $\eta$ (соответственно $\hat \eta$) на
элементарный страт гомотопно одному из отображений в
подпространства $K(\I_b,1) \subset K(\D_4,1)$, $K(\I_a,1) \subset
K(\D_4,1)$, $K(\I_d,1) \subset K(\D_4,1)$. Доказательство вытекает
из Леммы $\ref{lemma33}$, поскольку отображение
$\pi3^{(s)}_{\circ}$ естественно по отношению к каноническому
накрытию в образе и прообразе. Доказательство для отображения
$\hat \eta$ аналогично. Лемма $\ref{etaK30}$ доказана.
\[  \]

Определим коммутативную диаграмму 3-листных накрытий:
\begin{eqnarray}\label{pK3}
\begin{array}{ccc}
RK_{0\circ} & \longrightarrow & R\hat K_{0\circ} \\
&&\\
\downarrow p3_{\circ}& & \downarrow \hat p3_{\circ}\\
&&\\
K3_{0\circ} & \longrightarrow & \hat K3_{0\circ},\\
\end{array}
\end{eqnarray}

и коммутативную диаграмму отображений:

\begin{eqnarray}\label{piK3}
\begin{array}{ccc}
RK_{0\circ} & \longrightarrow & R\hat K_{0\circ}^1 \\
&&\\
\downarrow \pi3_{\circ}& & \downarrow \hat \pi3_{\circ}\\
&&\\
K_{0\circ} & \longrightarrow & \hat K_{0\circ}.\\
\end{array}
\end{eqnarray}

Диаграмма $(\ref{pK3})$ (соответственно $(\ref{piK3})$) определена аналогичными формулами как $(\ref{pK3s})$
(соответственно $(\ref{piK3s})$), вместо элементарного страта в базе накрытия
рассматривается все пространство $K3_{0\circ} \to \hat K3_{0\circ}$ (соответственно $K_{0\circ} \to \hat K_{0\circ}$). Диаграмма $\ref{piK3}$
также определена аналогично.
Диаграммы $(\ref{pK3})$, $(\ref{piK3})$ можно определить по-другому как результат склейки семейства диаграмм $(\ref{pK3s})$, $(\ref{piK3s})$
над всеми элементарными стратами глубины $s$, $0 \le s \le i_{max,0}$.

Определим пространства $RK_0$, $R\hat K_0$ и диаграмму накрытий
(вертикальные стрелки):
\begin{eqnarray}\label{ppK3}
\begin{array}{ccc}
RK_{0} & \longrightarrow & R\hat K_{0} \\
&&\\
\downarrow p3^1& & \downarrow \hat p3^1\\
&&\\
K3_{0} & \longrightarrow & \hat K3_{0},\\
\end{array}
\end{eqnarray}

\begin{eqnarray}\label{ppiK3}
\begin{array}{ccc}
RK_{0} & \longrightarrow & R\hat K_{0} \\
&&\\
\downarrow \pi3& & \downarrow \hat \pi3\\
&&\\
K_{0} & \longrightarrow & \hat K_{0}.\\
\end{array}
\end{eqnarray}

Пространство $RK_{0}$ (соответственно $R\hat K_{0}$) определим как
замыкание пространства $RK_{0\circ}$ (соответственно $R\hat
K_{0\circ}$) т.е. как результат присоединения всех элементарных
диагональных и антидиагональных стратов глубины $i$, $1 \le i \le
i_{max,0}$, лежащих на границе. Присоединение диагональных и
антидиагональных стратов происходит так, что каждый присоединенный
диагональный или антидиагональный элементарный глубины $i$ страт
лежит в границе ровно одного элементарного страта глубины $i-1$.

Слои над граничными элементарными стратами $K3_{0}\setminus K3_{0\circ}$
(соответственно $\hat K3_{0}\setminus \hat K3_{0\circ}$)
делятся на следующие типы. Для элементарного страта первого типа пространства
$K3_0 \setminus K3_{0\circ}$ (соответственно $\hat K3_{0}\setminus \hat K3_{0\circ}$))
прообраз при накрытии  $p3$ (соответственно
$\hat p3^{(s)}$ ) состоит из дизъюнктного объединения элементарного страта границы $Q_{diag}$ или
$Q_{antidiag}$ (соответственно $\hat Q_{diag}$ или
$\hat Q_{antidiag}$) и пары совпадающих регулярных элементарных стратов.
Для элементарного страта  второго типа  пространства $K3_0 \setminus K3_{0\circ}$
(соответственно $\hat K3_{0}\setminus \hat K3_{0\circ}$), прообраз состоит из трех
экземпляров элементарного страта границы $Q_{diag}$ или
$Q_{antidiag}$ (соответственно $\hat Q_{diag}$ или
$\hat Q_{antidiag}$).

\subsubsection*{Построение пространств $RK_{1\circ}$, $R\hat K_{1\circ}$, $RK_{1}$, $R\hat K_{1}$
для разрешения  особенностей и построение
 3-листных накрытий $p3_{\circ}: RK_{1\circ} \to
K3_{1\circ}$, $\hat p3_{\circ}: R\hat K_{1\circ} \to \hat
K3_{1\circ}$;  $p3: RK_{1} \to K3_{1}$, $\hat p3: R\hat K_{1} \to
\hat K3_{1}$.}

Построения аналогичны предыдущим. Сформулируем лемму, аналогичную
лемме $\ref{etaK30}$.

\begin{lemma}\label{etaK31}

Композиция $\zeta^{(s)} \circ \pi3^{(s)}_{\circ}:
RK^{(s)}_{1\circ} \longrightarrow K^{(s)}_{1\circ} \longrightarrow
K(\H,1)$ гомотопна отображению $\i_{\I_a,\H} \circ
\zeta_{RK^{(s)}_{1\circ}} : RK^{(s)}_{1\circ} \longrightarrow
K(\I_a, 1) \longrightarrow K(\H,1)$.

Композиция $\hat \zeta^{(s)} \circ \pi3^{(s)}_{\circ}: R \hat
K^{(s)}_{1\circ} \longrightarrow \hat K^{(s)}_{1\circ}
\longrightarrow K(\G,1)$ гомотопна отображению $ \i_{\I_a,\G}
\circ \zeta_{R \hat K^{(s)}_{1\circ}}:  R \hat K^{(s)}_{1\circ}
\longrightarrow K(\I_a, 1)  \longrightarrow K(\G,1)$.

\end{lemma}

\subsubsection*{Доказательство Леммы $\ref{etaK31}$}

Доказательство аналогично доказательству Леммы $\ref{etaK30}$ и
вытекает из Леммы $\ref{lemma33Q}$.
\[  \]

\subsubsection*{Пространства $RK^{(s)}_{0}$, $R\hat K^{(s)}_{0}$
разрешающие особенности включены в диаграграмму ($\ref{16.2}$);
пространства $RK^{(s)}_{1}$, $R\hat K^{(s)}_{1}$ разрешающие
особенности включены в диаграмму ($\ref{16.22}$)}

Приступим к построению отображений в диаграммах разрешающих
особенности.

\subsubsection*{Отображения $\phi^{(s)}_0$, $\hat \phi^{(s)}_0$ в диаграмме ($\ref{16.1}$)}

Определим  отображение $\phi^{(s)}_0$ как  структурные
отображения. Структурное отображение  $\eta_{R
K^{(s+1)}_{0\circ}}: RK^{(s+1)}_{0\circ} \to K(\I_c,1)$ на
пространстве-прообразе по Лемме $\ref{etaK30}$ гомотопно
отображению в подпространство $K(\I_d,1) \subset K(\I_a,1)$.

Отображение $\phi^{(s)}_{0\circ}$ продолжается до отображения
$\phi^{(s)}_{0}: RK^{(s)}_{0} \to  K(\I_a,1)$. Для построенного
отображения выполнены граничные условия $(\ref{22.1})$,
$(\ref{23.1})$.

Построение отображения $\hat \phi^{(s)}_0: R\hat K^{(s)}_{0} \to K(\I_a,1)$ совершенно аналогично.

\subsubsection*{Отображения $\phi^{(s)}_1$, $\hat \phi^{(s)}_1$ в диаграмме ($\ref{20.1}$)}

Это построение аналогично предыдущему и опускается.
\[  \]

Для определения отображений $\phi_0$ $\hat \phi_0$ в диаграмме
($\ref{16.2}$) нам потребуются вспомогательные построения.

\subsubsection*{Построение вспомогательных пространств $X^s_{0\circ}$, $X^s_{0}$,  $Y^s_{0\circ}$, $Y^s_{0}$,
$Z^s_{0\circ}$, $Z^s_{0}$ }

Обозначим пространство $K_{0\circ}^{s} \cup K_{0\circ}^{s+1}$
(соответственно $\hat K_{0\circ}^{s} \cup  \hat K_{0\circ}^{s+1}$)
через $Y^{s}_{0\circ}$ (сооответственно через $\hat
Y^{s}_{0\circ}$). Обозначим пространство $RK_{0\circ}^{s} \cup
RK_{0\circ}^{s+1}$ (соответственно $R\hat K_{0\circ}^{s} \cup
R\hat K_{0\circ}^{s+1}$) через $X^{s}_{0\circ}$ (сооответственно
через $\hat X^{s}_{0\circ}$).

Обозначим пространство $RK_{0\circ}^{s-1} \cup RK_{0\circ}^{s}
\cup RK_{0\circ}^{s+1}$ (соответственно $\hat RK_{0\circ}^{s-1}
\cup \hat RK_{0\circ}^{s} \cup \hat RK_{0\circ}^{s+1}$) через
$Z^{s-1}_{0\circ}$ (сооответственно через $\hat
Z^{s-1}_{0\circ}$).

Определены естественные включения
\begin{eqnarray}\label{59}
j_{X^{s}_{0\circ}}: X^{s}_{0\circ} \subset Z^{s-1}_{0\circ},
\end{eqnarray}

\begin{eqnarray}\label{59.1}
j_{\hat X^{s}_{0\circ}}: \hat X^{s}_{0\circ} \subset \hat Z^{s-1}_{0\circ},
\end{eqnarray}

\begin{eqnarray}\label{59.2}
j_{Y^{s}_{0\circ}}: Y^{s}_{0\circ} \subset K^{s}_{0\circ},
\end{eqnarray}

\begin{eqnarray}\label{59.3}
j_{\hat Y^{s}_{0\circ}}: \hat Y^{s}_{0\circ} \subset \hat K^{s}_{0\circ}.
\end{eqnarray}

Определены естественные проекции
\begin{eqnarray}\label{60}
p_{X^{s}_{0\circ}}: X^{s}_{0\circ} \to Y^s_{0\circ},
\end{eqnarray}

\begin{eqnarray}\label{61}
p_{\hat X^{s}_{0\circ}}: \hat X^{s}_{0\circ} \to \hat
Y^s_{0\circ}.
\end{eqnarray}

Определены пространства $X^s_0$, $Y^s_0$, $Z^{s-1}_0$ в результате
пополнения пространств $X^s_{0\circ}$, $Y^s_{0\circ}$,
$Z^{s-1}_{0\circ}$, т.е. добавления диагональных и
антидиагональных стратов. Определены естественные включения
$i_{X^s_{0\circ}}: X^s_{0\circ} \subset X^s_0$, $i_{Y^s_{0\circ}}:
Y^s_{0\circ} \subset Y^s_0$, $i_{Z^s_{0\circ}}: Z^s_{0\circ}
\subset Z^s_0$, при этом определена коммутативная диаграмма,
горизонтальные стрелки которой являются естественными вложениями:
\begin{eqnarray}\label{62}
 X^{s}_{0\circ} & \stackrel{i_{X^s_{0\circ}}}{\longrightarrow} & X^s_0\\
 \downarrow p_{X^{s}_{0\circ}}&  &\downarrow p_{X^{s}_0}\\
 Y^s_{0\circ} & \stackrel{i_{Y^s_{0\circ}}}{\longrightarrow} & Y^s_0.\\
\end{eqnarray}

Определено структурные отображения $\eta_{Y^s_{0\circ}}: Y^s_{0\circ} \to K(\D_4,1)$,
$\eta_{\hat Y^{s}_{0\circ}}: \hat Y^{s}_{0\circ} \to K(\E,1)$, которые получается в
результате ограничения структурных отображений
$\eta_{K_0}$, $\eta_{\hat K_0}$ при включениях $j_{Y^{s}_{0\circ}}$, $j_{\hat Y^{s}_{0\circ}}$.

Обозначим композицию $p_{X^s_{0\circ}} \circ \eta_{Y^s_{0\circ}}$ через
$\eta_{\hat X^s_{0\circ}}: X^s_{0\circ} \to K(\D_4,1)$, композицию
$p_{\hat X^s_{0\circ}} \circ \eta_{\hat Y^s_{0\circ}}$ через
$\eta_{\hat X^{s}_{0\circ}}: \hat X^{s}_{0\circ} \to K(\E,1)$.
Эти отображения также назовем структурными.

Определено пространство $\bar X^{s}_{0\circ}$, которое называется каноническим
2-листное накрывающим над $X^s_{0\circ}$, соответствующее 2-листное накрытие индуцированно из накрытия
$K(\I_b,1) \to K(\D_4,1)$ отображением $\eta_{X^s_{0\circ}}$. Обозначим это каноническое 2-листное накрытие через
$ p_{\bar X^{s}_{0\circ}}: \bar X^{s}_{0\circ} \to \bar X^s_{0\circ}$. Аналогично определяется пространство
$\tilde X^{s}_{0\circ}$ и 2-накрытие $ p_{\tilde X^{s}_{0\circ}}: \tilde X^{s}_{0\circ} \to \hat X^s_{0\circ}$.

\begin{lemma}\label{barX}
В предположении $s = -1 \pmod{4}$ cуществует отображение
$\phi_{X^s_{0\circ}}: X^s_{0\circ} \to K(\I_a,1)$ $($
соответственно $\phi_{\hat X^s_{0\circ}}: \hat X^s_{0\circ} \to
K(\I_a,1)$ $)$, удовлетворяющее следующему условию:

---  Отображение $\phi_{X^s_{0\circ}}: X^s_{0\circ} \to K(\I_a,1)$
продолжается до отображения $\phi_{Z^{s-1}_0}: Z^{s-1}_0 \to
K(\I_a,1)$ $($отображение $\phi_{\hat Z^{s-1}_{0\circ}}: \hat
Z^{s-1}_{0\circ} \to K(\I_a,1)$ продолжается до отображения
$\phi_{\hat Z^{s-1}_0}: \hat Z^{s-1}_0 \to K(\I_a,1)$$)$, причем
для продолженного отображения выполняются граничные условия:
$$ \eta_{RQ^s_0} = \phi_{X^s_{0\circ}} \vert_{RQ^s_0}: RQ^s_0 \to K(\I_a,1), $$
$$ \eta_{R \hat Q^s_0} = \phi_{\hat X^s_{0\circ}} \vert_{R \hat Q^s_0}: R \hat Q^s_0 \to K(\I_a,1). $$

\end{lemma}

\subsubsection*{Доказательство Леммы $(\ref{barX})$}

Начнем с построения отображения $\phi_{X^s_{0\circ}}$. Обозначим
через $X^s_{\I_b\circ}$, (соответственно $X^s_{\I_a,\circ}$)
подполиэдры в $X^s_{0\circ}$, которые определим как объединение
всех элементарных стратов глубины $s$, для которых произведение
вычетов $\prod_{i=1}^{r_0} v_{1,2;i}$  выделенных координат
$(\check{x_1}, \check{x_2})$ из тройки координат
$(\check{x_1},\check{x_1},\check{x_3})$  равно $\pm 1$
(соответственно $\pm \i$) и всех  элементарных стратов глубины
$s+1$, которые лежат на границе указанных выше элементарных
стратов глубины $s$.

Построим отображения $\phi_{X^s_{\I_b\circ}}: X^s_{\I_b\circ} \to
K(\I_a,1)$, $\phi_{X^s_{\I_a\circ}}: X^s_{\I_a\circ} \to
K(\I_a,1)$. Далее проверим, что построенные отображения
склеиваются в одно отображение $\phi_{X^s_{0\circ}}: X^s_{0\circ}
\to K(\I_a,1)$.

Начнем с определения гомоморфизма $\phi_{X^s_{\I_b\circ},\ast}:
\pi_1(X^s_{\I_b\circ}) \to \I_a$.
 Пусть $l: S^1 \subset X^s_{\I_b
\circ}$ -- произвольный кусочно-линейный путь, трансверсально
пересекающий подпространство особых стратов $PK^{(s-1)}_{0,\circ}
\subset X^s_{\I_b \circ}$ в конечном числе точек, которые
обозначим через $\{ t_1, \dots, t_j\}$.

На  элементарном страте глубины $s$ полиэдра $\bar X^s_{\I_b
\circ}$, содержащем точку $pt=l(pt_0)$, $pt_0 \in S^1$--отмеченная
точка, выберем систему координат
$(\check{x_1}(pt),\check{x_2}(pt),\check{x_3}(pt))$ (замечу, что
первые две координаты произвольной точки на  элементарном страте
пространства $RK^{(s)}_{0\circ}$ неупорядочены). Обозначим
$\prod_{i=1}^r v_{1,2;i}=s(1,2)$, $\prod_{i=1}^r v_{3,1;i}=s(1,3)$
$\prod_{i=1}^r v_{2,3;i}=s(2,3)$. В частности, поскольку $pt \in
X^s_{\I_b \circ}$, можно выбрать систему координат
$(\check{x_1}(pt),\check{x_2}(pt),\check{x_3}(pt))$ так, что
произведение вычетов
\begin{eqnarray}\label{s12}
s(1,2)=-1.
\end{eqnarray}

Поскольку справедливо равенство $s(1,2)s(2,3)s(3,1)=1$, то
\begin{eqnarray}\label{s13}
s(1,3)=s(2,3)^{-1}.
\end{eqnarray}

(Если дополнительно предположить, что $s(1,3) = 1$, то
$s(2,3)=-1$.) Указанное соглашение о произведении вычетов
облегчает проверку (по сравнению с условием $s(1,2)=1$) граничного
условия на диагональной компоненте.

Рассмотрим произвольную точку $t \in \{t_1, \dots, t_j\}$.
Продолжим естественным способом систему координат
$(\check{x_1}(pt),\check{x_2}(pt))$ из точки $pt$ вдоль пути, т.е.
согласовывая регулярные координаты в окрестности каждой из точек
$\{ t_1, \dots, t_j\}$ (по поводу естественного продолжения
системы координат см. координатное построение структурного
отображения $\eta$). Обозначим через
$(\check{x_1}(t-\varepsilon),\check{x_2}(t-\varepsilon))$,
$(\check{x_1}(t+\varepsilon),\check{x_2}(t+\varepsilon))$ системы
координат, полученные продолжением исходной системы координат в
точки $t-\varepsilon$, $t+\varepsilon$ соответственно.

Определим понятие преобразования параллельного выделенной третьей
координате вдоль пути $l$. Пусть $t \in \{ t_1, \dots, t_j\}$.

Обозначим через $\alpha_1(t)$, $\alpha_2(t)$ -- два максимальных
элементарных страта, которые примыкают друг к другу по
элементарному страту $\beta(t)$ глубины $s+1$, содержащему точку
$l(t)$. На страте $\alpha_1(t)$ выбрана система координат
$(\check{x_1}(t-\varepsilon),
\check{x_2}(t-\varepsilon),\check{x_3}(t-\varepsilon))$. На страте
$\alpha_2(t)$ выбрана система координат
$(\check{x_1}(t+\varepsilon),\check{x_2}(t+\varepsilon),\check{x_3}(t-\varepsilon))$.
В гомотопическом классе $l$ выберем  путь, что при каждом $t \in
\{ t_1, \dots, t_j\}$ каждая регулярная координата набора
$\check{x}_{i;j}$, $i \ne i_0$ или $j \ne j_0$, в окрестности $t$
сохраняется, т.е
$\check{x}_{i;j}(t-\varepsilon)=\check{x}_{i,j}(t+\varepsilon)$.
При этом критическая координата $(\check{x}_{i_0;j_0}$ меняется
умножением на элемент из $\I_a$. Нетрудно видеть, что либо одна из
координат $\check{x}_{1;j_0})$, $\check{x}_{2;j_0})$
 меняется умножением на $-1$, либо координата $\check{x}_{2;j_0})$ меняется умножением на один из элементов
 $\{-1, +\i, -\i\}$ из $\I_a$. Такой путь назовем стандартным.


Помимо системы координат
$(\check{x}_1(2\pi),\check{x}_2(2\pi),\check{x}_3(2\pi))$,
определенной в результате естественного продолжения системы
координат $(\check{x}_1(0),\check{x}_2(0),\check{x}_3(0))$ вдоль
$l$, определим еще одну систему координат, которую обозначим через
$(\check{y}_1(2\pi),\check{y}_2(2\pi),\check{y}_3(2\pi))$. Эта
система координат будет называться преобразованной системой
координат вдоль стандартизованного $l$ параллельно третьей
координате.

Разобъем множество критических точек пути $l$ на два подмножества
$\{ t_1, \dots, t_j\} = U \cup V$. Отнесем к подмножеству
$U$--критические точки, в которых первая или вторая координаты
меняются на элемент $\pm 1$. Отнесем к подмножеству
$V$--критические точки, в которых третья координата меняется на
элемент из $\I_a$.

Обозначим критические точки подмножества $U$ через $U=\{\tau_1,
\dots, \tau_k\}$. На каждом из отрезков $[\tau_i \tau_{i+1}]$ пути
$l$ (возможно, что начальная или конечная точка такого пути
совпадет с отмеченной) продолжим систему координат следующим
способом.

Для каждого отрезка $[\tau_i \tau_{i+1}]$  определим
$\check{y}_1(\tau_{i+1}) = \theta(1,1)\check{x}_1(\tau_{i+1})$,
$\check{y}_2(\tau_{i+1}) = \theta(2,2)\check{x}_2(\tau_{i+1})$,
$\check{y}_3(\tau_{i+1}) = \check{x}_3(\tau_{i+1})$.

Значения $\theta(1,1), \theta(2,2) \in \I_d$ выбраны так, что
выполнены равенства:

$$\prod_{i=1}^{r_0} w_{1,2;i}=\prod_{i=1}^{r_0} v_{1,2;i},$$
$$\prod_{i=1}^{r_0} w_{2,3;i}=\prod_{i=1}^{r_0} v_{2,3;i},$$
$$\prod_{i=1}^{r_0} w_{3,1;i}=\prod_{i=1}^{r_0} v_{3,1;i},$$
 где
 $v_{1,2;i}, v_{2,3;i}, v_{3,1;i}$-- вычеты пар соответствующих
координат $(\check{y}_1(\tau_i),\check{y}_2(\tau_i),
\check{y}_2(\tau_i+0))$,  $w_{1,2;i}, w_{2,3;i}, w_{3,1;i}$--
вычеты пар соответствующих координат
$(\check{y}_1(\tau_{i+1}),\check{y}_2(\tau_{i+1}))$,
$(\check{y}_3(\tau_{i+1}-0)$, $(\check{y}_3(\tau_{i+1}-0)$
получена естественным продолжением $(\check{y}_3(\tau_{i}+0)$
вдоль $[\tau_i \tau_{i+1}] \subset l$.

Иными словами, определим систему координат
$(\check{y}_1(\tau_{i+1}),\check{y}_2(\tau_{i+1}),\check{y}_3(\tau_{i+1}-0))$
так, что третья координата этой новой системы определена
естественным продолжением, а попарные произведения вычетов
координат в правом конце отрезка такие же как и в левом.

Геометрически можно представлять себе, что система координат
$(\check{y}_1(\tau_{i+1}),\check{y}_2(\tau_{i+1}),\check{y}_3(\tau_{i+1}-0))$
определена следующим способом. Сначала переносим систему координат
$(\check{x}_1(\tau_i),\check{x}_2(\tau_i),\check{x}_3(\tau_i+0))$
вдоль $l$, причем в каждый критический момент времени один  из
наборов критических координат с номером $1$ или $2$ изменяется.
Пусть для определенности набор с номером $1$. Подправим весь набор
координат c номером $1$ так, чтобы в критический момент времени
поднабор регулярных координат с номером 1 изменяется на обратный,
а особая координата с номером 1 в критической точке не изменяется.

В точках $U$, где меняется одна координата из набора третьих
координат на элемент $x \in \I_a$, определим преобразование
наборов первых и вторых координат при котором  каждая регулярная
третья координата преобразуется на элемент $x \in \I_a$.
 В результате получим
преобразование третьих координат, при котором произведение вычетов
по всем координатам в начальной и преобразованной системы не
изменяется (здесь используется равенство $s=-1 \pmod 4$).

Преобразование исходной системы координат
$(\check{x}_1(0),\check{x}_2(0),\check{x}_3(0))$ в систему
координат
$(\check{y}_1(2\pi),\check{y}_2(2\pi),\check{y}_3(2\pi))$
определим как последовательное преобразование вдоль системы
отрезков, определенных по $U$. Это преобразование представлено
прямой суммой двух преобразований $\aleph_{1,2} \in \I_b$,
$\aleph_{3} \in \I_d$ по первым двум и третьей координате
соответственно. По построению $\aleph_{3}$ совпадает с
естественным преобразованием третьей координаты.

Определим значение гомоморфизма $\phi_{X^s_{\I_b\circ}}$ на
гомотопическом классе пути $[l]$ по формуле
$\phi_{X^s_{\I_b\circ},\ast}([l])=i_{\I_d,\I_a} \circ \bar
p_{\I_b,\I_d}(\aleph_{1,2})$, где $\aleph_{1,2} \in \I_b$, $\bar
p_{\I_b, \I_d}: \I_b \to \I_d$-- проекция с ядром $a^2 b \in
\I_b$. (Если вместо условия $(\ref{s12})$ использовать условие
$s(1,2)=1$, то проекцию $\bar p_{\I_b,\I_d}$ следует заменить на
проекцию $p_{\I_b,\I_d}$ с ядром $b$.)

Докажем, что гомоморфизм $\phi_{X^s_{\I_b\circ},\ast}:
\pi_1(X^s_{\I_b\circ}) \to \I_a$ корректно определен.

Пусть система координат
$(\check{x}'_1(0),\check{x}'_2(0),\check{x}'_3(0))$ получается
из системы координат
$(\check{x}_1(0),\check{x}_2(0),\check{x}_3(0))$
преобразованием на $x \in \I_b$. Тогда система координат
$(\check{y}'_1(2\pi),\check{y}'_2(2\pi),\check{y}'_3(2\pi))$, полученная
в результате продолжения
$(\check{x}'_1(0),\check{x}'_2(0),\check{x}'_3(0))$ вдоль $l$
параллельно третьей координате также
получается из системы координат
$(\check{y}_1(2\pi),\check{y}_2(2\pi),\check{y}_3(2\pi))$
преобразованием на $x \in \I_b$.

Пусть стандартизация пути $l$ выбрана двумя способами. При
элементарных перестройках одной стандартизации в другую значение
$\phi_{X^s_{\I_b\circ}}$ сохраняется. В частности, пусть
элементарная перестройка двух стандартизаций пути $l$ состоит в
следующем. Элементарная особенность множества $U$, в которой
$i$--ая критическая третья координата меняется на $-1 \in \I_d$,
заменяется на две элементарные особенности множества $V$, в одной
из которых $i$--ая первая, а в другой $i$--ая вторая координата
являются критическими и изменяются на элемент $-1 \in \I_d$. В
этом случае преобразование $\aleph_{1,2}$ не  изменяется.
Корректность $\phi_{X^s_{\I_b\circ},\ast}$ доказана.

Приступим к определению гомоморфизма $\phi_{X^s_{\I_a\circ},\ast}:
\pi_1(X^s_{\I_a\circ}) \to \I_a$.

 На
элементарном страте глубины $s$ полиэдра $\bar X^s_{\I_a\circ}$,
содержащем точку $pt=l(pt_0)$, $pt_0 \in S^1$--отмеченная точка,
выберем систему координат
$(\check{x_1}(0),\check{x_2}(0),\check{x_3}(0))$. Обозначим
$\prod_{i=1}^r v_{1,2;i}=s(1,2)$, $\prod_{i=1}^r v_{3,1;i}=s(3,1)$
$\prod_{i=1}^r v_{2,3;i}=s(2,3)$. Упорядочим две первые координаты
и подправим координаты умножением на $\pm 1$ так, что
произведения вычетов
\begin{eqnarray}\label{s31Q}
s(3,1)=1, s(3,2)=\i,
\end{eqnarray}
 тогда
\begin{eqnarray}\label{s12Q}
 s(1,2)=+\i,
\end{eqnarray}
  поскольку
$s(3,2)=-s(2,3)$, $s(1,2)s(2,3s(3,1)=1$. (Можно воспользоваться
другим соглашением о выборе вычетов, например, после умножения
набора вторых координат на $\i$, второе уравнение из
$(\ref{s31Q})$ перейдет в уравнение $s(3,2)=1$.)

Пусть $l: S^1 \subset X^s_{\I_a \circ}$ -- произвольный
кусочно-линейный путь, трансверсально пересекающий подпространство
особых стратов $PK^{(s-1)}_{0,\circ} \subset X^s_{\I_b \circ}$ в
конечном числе критических точек, которые обозначим через $\{ t_1,
\dots, t_j\}$. Выберем путь $l$ стандартизованным. Разобъем
множество критических точек на два подмножества  $\{ t_1, \dots,
t_j\} = U \cup V$. Отнесем к подмножеству $U$--критические точки,
в которых одна третья координата из набора третьих координат
меняется на элемент из $\I_a$. Отнесем к подмножеству
$V$--критические точки, в которых одна координата из набора первой
или второй координат меняется на элемент $-1 \in \I_d$.

Обозначим критические точки множества $U$ через $U=\{\tau_1,
\dots, \tau_k\}$. На каждом из отрезков $[\tau_i \tau_j]$  пути
$l$ (возможно, что начальная или конечная точка этого пути
совпадет с отмеченной) продолжим систему координат параллельно
третьей координате. В каждой критической точке множества $V$
преобразуем наборы первых и вторых координат системы, чтобы
произведение вычетов набора третьих координат вдоль пути не
изменилось.

В результате получим преобразование исходной системы координат,
которое обозначим через $\aleph$. При этом $\aleph = \aleph_{1,2}
\oplus \aleph_3$, $\aleph_{1,2} \in \E_a$, $\aleph_3 \in \I_d$.
Преобразование $\aleph_3$ совпадает с естественным.

Определим гомоморфизм $\phi_{X^s_{\I_a\circ},\ast}:
\pi_1(X^s_{\I_a\circ}) \to \I_a$ по формуле
$\phi_{X^s_{\I_a\circ},\ast}([l])=\bar p_{\E_a,
\I_a}(\aleph_{1,2})$, где $\bar p_{\E_a, \I_a}: \E_a \to \I_a$ --
проекция с ядром $aс$. (Напомню, что $\E_a$ является квадратичным
расширением подгруппы $\I_a$, эта группа изоморфна группе $\Z/4
\oplus \Z/2$. Элемент $a^2$ расширяется элементом $c$ четвертого
порядка, $c^2 = a^2$. Элемент $c$ действует умножением на $\i$ на
каждую координату.) (Если вместо равенства $(\ref{s12Q})$
использовать равенство $s(1,2)=1$, то $\aleph_{1,2} \in \E_b$ и
проекцию $\bar p_{\E_a, \I_a}: \E_a \to \I_a$ следует заменить на
проекцию  $p_{\E_b, \I_a}: \E_b \to \I_a$ с ядром $b$.)

Докажем, что гомоморфизм $\phi_{X^s_{\I_a\circ},\ast}:
\pi_1(X^s_{\I_a\circ}) \to \I_a$ корректно определен. Система
координат $(\check{x_1}(0),\check{x_2}(0),\check{x_3}(0))$
выбирается неоднозначно, в частности, координата $\check{x_1}(0)$
определена с точность до преобразования $-1 \in \I_d$. Пусть
$\check{x_1}'(0)$ другая координата такая, что
$\check{x_1}'(0)=-\check{x_1}(0)$. Пусть координата
$\check{x}'_1(2\pi)$ получена из $\check{x_1}(0)$ в результате
естественного преобразования вдоль $l$. Пусть координата
$\check{x}'_2(2\pi)$ получена из $\check{x}_2(0)$ в результате
естественного преобразования вдоль $l$.Тогда, очевидно,
$\check{x}'_1(2\pi) = -\check{x}'_1(2\pi)$. Следовательно,
преобразование $\aleph \in \I_a$ не зависит от выбора
первоначальной системы координат.

При изменении стандартизации пути $l$, преобразование
$\aleph_{1,2}$ не меняется. Гомоморфизм
$\phi_{X^s_{\I_a\circ},\ast}: \pi_1(X^s_{\I_a\circ}) \to \I_a$
корректно определен.

 Поскольку по Лемме $(\ref{etaK30})$  структурное
отображение $\eta_{X^s_{\circ}}$ на  подпространстве
$X^s_{\I_b\circ} \cap X^s_{\I_a\circ} \subset X^s_{\circ}$
гомотопно отображению в $K(\I_d,1) \subset K(\D_4,1)$. классифицирующие отображения
$\phi_{X^s_{\I_a\circ}}$, $\phi_{X^s_{\I_b\circ},\ast}$
согласованы на $X^s_{\I_a\circ} \cap X^s_{\I_b\circ}$ и
определяют отображение $\phi_{X^s_{0\circ}}: X^s_{0\circ} \to
K(\I_a,1)$. Отображение $\phi_{\hat X^s_{0\circ}}: \hat X^s_{0\circ}
\to K(\I_a,1)$ определяется аналогичным способом.

Докажем условие.   Проверим, что отображение $\phi_{X^s_{0\circ}}:
X^s_{0\circ} \to K(\I_a,1)$ продолжается до отображения
$\phi_{Z^{s-1}_0}: Z^{s-1}_0 \to K(\I_a,1)$. Сначала заметим, что
отображение $\phi_{X^s_{0\circ}}: X^s_{0\circ} \to K(\I_a,1)$
продолжается до отображения $\phi_{X^s_{0}}: X^s_{0} \to
K(\I_a,1)$. Докажем, что отображение $\phi_{X^s_{0}}: X^s_{0} \to
K(\I_a,1)$ продолжается до отображения пространства $Z^{s-1}_0$,
т.е. проверим, что существует продолжение на каждый элементарный
страт глубины $s-1$.

Пусть $a \subset PK^{(s-1)}_{0\circ}$--элементарный страт глубины
$s-1$, причем элементарные страты глубины $s$ и $s+1$, лежащие на
его границе, отождествлены с соответствующим подпространством  в
$X^s_{0\circ}$. Для произвольного замкнутого пути $l$, лежащего на
элементарных стратах глубины $s$ и $s+1$ на границе $a$ определено
естественное преобразование системы координат вдоль $l$ и элемент
$\eta_{0,\ast}([l])$, Воспользуемся Леммой $\ref{etaK30}$,
согласно которой $\eta_{0,\ast}([l]) \in \I_d$. Кроме того, вдоль
$l$ определено преобразование $\aleph$ и, очевидно,
$\aleph=\eta_{0,\ast}([l])$. Поэтому отображение
$\phi_{X^s_{0\circ}}$ продолжается до отображения
$\phi_{Z^{s-1}_0}$ (и притом единственным способом). Построение
отображения $\phi_{Z^{s-1}_0}: Z^{s-1}_0 \to K(\I_a,1)$
аналогично.

Проверим, что выполняется граничное условие $ \eta_{PQ^s_0} =
\phi_{X^s_{0}} \vert_{PQ^s_0}: PQ^s_0 \to K(\I_a,1)$. Это
равенство достаточно проверить для диагональной и антидиагональной
компоненты границы пространства $X^s_{0}$. Проверим это для
диагональной компоненты, которая содержится в подпространстве
$X^s_{0,\I_b}$. Для антидиагональной компоненты, которая
содержится в подпространстве $X^s_{0,\I_a}$ проверка аналогична.

Рассмотрим замкнутый путь $l$ на компоненте $PQ^s_{diag}$. Пусть
$\aleph$ -- преобразование системы координат вдоль $l$. В случае
$\theta(3,1)=1$, преобразование координат  $\aleph$ вдоль $l$
совпадает с естественным преобразованием координаты $\check{x}_3$
вдоль $l$, т.к. в этом случае система координат
$(\check{y}_1(2\pi),\check{y}_2(2\pi))$ совпадает с системой
координат $(\check{x}_1(2\pi),\check{x}_2(2\pi))$.  В этом случае
$\aleph=\eta_{0,\ast}([l]) \in \I_d$. Если $\theta(3,1)=-1$, то
система координат $(\check{y}_1(2\pi),\check{y}_2(2\pi))$
отличается от системы координат
$(\check{x}_1(2\pi),\check{x}_2(2\pi))$, а на образующую $-1 \in
\I_d$, но при этом $\aleph$ снова совпадает с естественным
преобразованием системы координат. В этом случае снова имеем
$\aleph=\eta_{0,\ast}([l]) \in \I_d$.

Граничные условия $ \eta_{P \hat Q^s_0} = \phi_{\hat X^s_{0}}
\vert_{P \hat Q^s_0}: P \hat Q^s_0 \to K(\I_a,1) $ проверяются
аналогично.

Лемма $(\ref{barX})$ доказана.


 Остановимся лишь на формулировке
следующей леммы $(\ref{barXQ})$, которая аналогична Лемме
$(\ref{barX})$. Эта лемма будет использована для построения
отображений $\phi_1$, $\hat \phi_1$ в диаграмме $(\ref{16.22})$.

\subsubsection*{Построение вспомогательных пространств $X^s_{1\circ}$, $X^s_{1}$,  $Y^s_{1\circ}$, $Y^s_{1}$,
$Z^s_{1\circ}$, $Z^s_{1}$ }

Это построение полностью аналогично и опускается.

Определено структурные отображения $\zeta_{Y^s_{1\circ}}: Y^s_{1\circ} \to K(\H,1)$,
$\zeta_{\hat Y^{s}_{1\circ}}: \hat Y^{s}_{1\circ} \to K(\G,1)$,
которые получается в результате ограничения структурных отображений
$\zeta_{K_1}$, $\zeta_{\hat K_1}$ при включениях $j_{Y^{s}_{1\circ}}$, $j_{\hat Y^{s}_{1\circ}}$.

Обозначим композицию $p_{X^s_{1\circ}} \circ \zeta_{Y^s_{1\circ}}$ через
$\zeta_{\hat X^s_{1\circ}}: X^s_{1\circ} \to K(\H,1)$, композицию
$p_{\hat X^s_{1\circ}} \circ \zeta_{\hat Y^s_{1\circ}}$ через
$\zeta_{\hat X^{s}_{1\circ}}: \hat X^{s}_{1\circ} \to K(\G,1)$.
Эти отображения также назовем структурными.

Определено пространство $\bar X^{s}_{1\circ}$, которое называется каноническим
2-листное накрывающим над $X^s_{1\circ}$, соответствующее 2-листное накрытие индуцированно из накрытия
$K(\I_b,1) \to K(\D_4,1)$ отображением $\zeta_{X^s_{1\circ}}$. Обозначим это каноническое 2-листное накрытие через
$ p_{\bar X^{s}_{1\circ}}: \bar X^{s}_{1\circ} \to \bar X^s_{1\circ}$. Аналогично определяется пространство
$\tilde X^{s}_{1\circ}$ и 2-накрытие $ p_{\tilde X^{s}_{1\circ}}: \tilde X^{s}_{1\circ} \to \hat X^s_{1\circ}$.

\begin{lemma}\label{barXQ}
В предположении $s = 1 \pmod{2}$ cуществует отображение
$\phi_{X^s_{1\circ}}: X^s_{1\circ} \to K(\Q_a,1)$ $($
соответственно $\phi_{\hat X^s_{1\circ}}: \hat X^s_{1\circ} \to
K(\Q_a,1)$ $)$, удовлетворяющее следующему условию:

---  Отображение $\phi_{X^s_{1\circ}}: X^s_{1\circ} \to K(\Q_a,1)$ продолжается до отображения
$\phi_{Z^{s-1}_1}: Z^{s-1}_1 \to K(\Q_a,1)$ $($отображение
$\phi_{\hat Z^{s-1}_{1\circ}}: \hat Z^{s-1}_{1\circ} \to
K(\Q_a,1)$ продолжается до отображения $\phi_{\hat Z^{s-1}_1}:
\hat Z^{s-1}_1 \to K(\Q_a,1)$$)$, причем для продолженного
отображения выполняются граничные условия:
$$ \zeta_{RQ^s_1} = \phi_{X^s_{1\circ}} \vert_{RQ^s_1}: RQ^s_1 \to K(\Q_a,1), $$
$$ \zeta_{R \hat Q^s_1} = \phi_{\hat X^s_{1\circ}} \vert_{R \hat Q^s_1}: R \hat Q^s_1 \to K(\Q_a,1). $$

\end{lemma}

\subsubsection*{Доказательство Леммы $(\ref{barXQ})$}

Доказательство аналогично доказательству Леммы  $(\ref{barX})$ и
опускается.

\subsubsection*{Отображения $\phi_0$, $\hat \phi_0$ в диаграмме $(\ref{16.2})$}

Построим отображение $\phi_0$ в диаграмме $(\ref{16.2})$.  В
пространстве $RK_{0\circ}$ определим подпространство
$RK^{[2]}_{0\circ} \subset RK_{0\circ}$. Для этого рассмотрим
следующую пару полиэдров (верхняя диаграмма) в вкложение в эту
другой пары полиэдров  (нижняя диаграмма):

\begin{eqnarray}\label{117}
\begin{array}{ccc}
RK_{0}^1 & \subset & RK_{0}^0  \\
\end{array}
\end{eqnarray}

\begin{eqnarray}\label{118}
\begin{array}{ccc}
X^1_{0} & \subset & Z^0_{0}  \\
\end{array}
\end{eqnarray}

Пространство $RK^{[2]}_{0}$ определено как образ пространства
$Z^0_{0}$ из диаграммы $(\ref{118})$ при включении.  Нам
потребуется следующая лемма.

\begin{lemma}\label{RK[2]}
Произвольное $PL$--отображение двумерного диска $f: D^2 \to RK_{0}$ вытесняется в отображение
$f': D^2 \to RK^{[2]}_{0} \subset RK_{0}$ с образом в подпространстве
$RK^{[2]}_{0} \subset RK_{0}$ в результате сколь угодно малой $PL$--деформации
$f \to f'$.
\end{lemma}

\subsubsection*{Доказательство Леммы $\ref{RK[2]}$}

Пространство $RK_{0}$ снабжено естественным отображением в пространство $\delta^{r_0} \times I$
прямого произведение $r_0$--мерного симплекса на отрезок, причем стандартная триангуляция
пространства-образа, задающая стратификацию пространства--образа по коразмерности остовов триангуляции  согласована
при этом отображении со стратификацией пространства-прообраза. Подпространство $RK^{[2]}_{0} \subset RK_{0}$
содержит все страты пространства $RK_{0}$ глубины 0,1,2. Отображение $f: D^2 \to RK_{0}$ можно перевести
сколь угодно малой деформацией в отображение $f': D^2 \to RK^{[2]}_{0}$ так, что проекция образа отображения
$f'$ в пространство
$\delta^{r_0} \times I$ трансверсально к симплексам разбиения и, в частности, не содержит остова коразмерности 3.
Лемма  $\ref{RK[2]}$ доказана.
\[  \]

Определим на подпространстве $RK^{[2]}_{0} \subset RK_{0}$ отображение
\begin{eqnarray}\label{phi[2]}
\phi^{[2]}_0: RK^{[2]}_{0} \to K(\I_a,1).
\end{eqnarray}

На подпространстве $X_0^1 \subset RK^{[2]}_{0}$ отображение
определим как $\phi_{X^1_{0}}: X^1_{0} \to K(\I_a,1)$, построенное
в Лемме $(\ref{barX})$ при $s=1$. По этой же лемме отображение
продолжается до отображения $\phi_{Z^0_{0}}: Z^0_{0} \to
K(\I_a,1)$. Тем самым, отображение $(\ref{phi[2]})$ определено.

По Лемме $\ref{RK[2]}$ отображение $(\ref{phi[2]})$ продолжается
до отображения $\phi_0: RK_0 \to K(\I_a,1)$. По Лемме
$(\ref{barX})$ отображение $(\ref{phi[2]})$ удовлетворяет
сформулированным граничным условиям. Следовательно, отображение
$\phi_0: RK_0 \to K(\I_a,1)$ удовлетворяет граничным условиям,
коорые следут из диаграмм $(\ref{18.1})$, $(\ref{18.2})$.
Построение отображения $\hat \phi_0: R\hat K_0 \to K(\I_a,1)$
полностью аналогично.

\subsubsection*{Отображения $\phi_1$, $\hat \phi_1$ в диаграмме
($\ref{16.22}$)}

Это построение аналогично предыдущему и опускается.
\[  \]

\subsubsection*{Доказательство Леммы $\ref{lemma30}$ }

Рассмотрим отображение $\hat p: S^{n-k}/\i \to J_0$ и обозначим
через $\hat C_{p}$ цилиндр этого отображения. Определены
отображения проекций $\hat p_I: \hat C_{p} \to [0,1]$, $\hat p_J:
\hat C_{p} \to J_0$ и декартово произведение этих отображений,
которое обозначим через $\hat F:  \hat C_{ p} \to J_0 \times
[0,1]$.

Аналогично рассмотрим отображение $p: \RP^{n-k} \to J_0$ и
обозначим через $C_p$ цилиндр этого отображения. Определены
отображения проекций $p_I: C_p \to [0,1]$, $p_J: C_p \to J_0$ и
декартово произведение этих отображений, которое обозначим через
$F: C_{ p} \to J_0 \times [0,1]$.

Определено отображение $r: C_p \to \hat C_{p}$, причем следующие
диаграммы коммутативны:

\begin{eqnarray}
\begin{array}{ccc}
C_p & \longrightarrow & \hat C_{ p} \\
p_I \searrow & & \swarrow \hat  p_I  \\
& I,  & \\
\end{array}
\end{eqnarray}

\begin{eqnarray}
\begin{array}{ccc}
C_p & \longrightarrow & \hat C_{p} \\
p_J \searrow & & \swarrow \hat p_J  \\
& J_0. & \\
\end{array}
\end{eqnarray}

Рассмотрим вложение $I_J: J_0 \times [0,1] \subset \R^n \times
[0,1]$ и определим отображения $I_J \circ \hat F: \hat C_{ p} \to
\R^n \times [0,1]$, $I_J \circ F: C_{p} \to \R^n \times [0,1]$.
Ниже определим отображение $\hat f: \hat C_{p} \to \R^n \times
[0,1]$, которое получено в результате $C^1$--малого возмущения общего
положения отображения $I_J \circ \hat F$. Такое возмущение будем ниже называть деформацией
гомотопии $I_J \circ \hat F$ в гомотопию $\hat f$. При этом гомотопию
$\hat f$ выберем совпадающей на нижнем основании цилиндра $J_0
\subset \hat C_{p}$ с вложением $I_J: J_0 \subset \R^n \times
\{0\}$. Кроме того, потребуем, чтобы композиция $p_{[0,1]} \circ
\hat f : \hat C_{p} \to [0,1]$ совпадала с $\hat p_{I}$, где
$p_{[0,1]}: \R^n \times [0,1] \to [0,1]$ -- проекция на второй
сомножитель. Определено также отображение $f: C_p \to \R^n \times
[0,1]$ как композиция $f = \hat f \circ r$.

Обозначим через $\bar Q \subset \hat C_p$ -- каноническое
накрывающее над полиэдром точек самопересечения гомотопии $\hat
f$, определенное как замыкание по формуле
$$ \bar Q = Cl\{ x \in \hat C_p : \exists y \in \hat C_p, x \ne y, \hat
f(x) = \hat f(y) \}. $$  Определена инволюция $T_{\bar Q}: \bar Q
\to \bar Q$, переставляющая прообразы точек самопересечения.
Инволюция $T_{\bar Q}$ сохраняет значения отображения $p_I$.
Факторпространство полиэдра $\bar Q$ при этой инволции обозначим
через $\hat Q$.

Рассмотрим стратификацию
\begin{eqnarray}
J^{r}_0 \subset \dots \subset J^{2}_0 \subset J^{1}_0 \subset
J^{0}_0 \label{strat}
\end{eqnarray}
 пространства джойна $J^{0}_0=J_0$.
Стратификация ($\ref{strat}$) определяет стратификацию цилиндра по
формуле
 \begin{eqnarray}
\hat C_p^{i} = \hat p_J^{-1}(J^{i}_0)
\end{eqnarray}
и стратификацию многообразия $S^{n-k}/\i$ по формуле
 \begin{eqnarray}
(S^{n-k}/\i)^{i} = \hat p^{-1}(J^{i}_0)
\end{eqnarray}

Определена стратификация
\begin{eqnarray}\bar Q^{r} \subset
\dots \subset \bar Q^{1} \subset \bar Q^{0} \label{strat barQ}
\end{eqnarray}
полиэдра $\bar Q^{0} = \bar Q$ по формуле $\bar Q^{i} = (\bar Q
\cap \hat С_p^{i}) \cap T_{\bar Q}(\bar Q \cap \hat C_p^{i})$.
Фактор стратификации $\ref{strat barQ}$ при инволюции $T_{\bar Q}$
определяет стратификацию
\begin{eqnarray}
\hat Q^{r} \subset \dots \subset \hat Q^{1} \subset \hat Q^{0}
\label{stratQ}
\end{eqnarray}
полиэдра $\hat Q$. Обозначим через  $QJ^{i}_0$ пересечение $\bar
Q^{i} \cap J_0$. При условии, что отображение $\hat f$ общего
положения, стратифицированный полиэдр $QJ^{i}_0$ находится в общем
положении со стратами ($\ref{strat}$) полиэдра $J_0$

При произвольном $t \in (0,+\delta]$ ограничение отображения $
\hat f: \hat C_{p} \to \R^n \times [0,1]$ на $S^{n-k}/\i \times
\{t\}$ обозначим через $\hat d = \hat d(t): S^{n-k}/\i \to \R^n$.
Отображение $d$ двулистно накрывает отображение $\hat d$.
Рассмотрим многообразие $\hat N_{K_0}$  с краем точек
самопересечения отображения $\hat d$. Многообразие $\hat N_{K_0}$
определяется по формуле $\hat N_{K_0} = \hat Q \cap \R^n \times
\{t\}$. Это многообразие снабжено естественной стратификацией
\begin{eqnarray}\label{stra0}
\hat N_{K_0}^{r} \subset \dots \subset \hat N_{K_0}^{0},
\label{strathat}
\end{eqnarray}

\begin{eqnarray}\label{stra00}
\hat W_{K_0}^{(i)} = \hat N_{K_0}^{i} \setminus \hat
N_{K_0}^{i+1},
\end{eqnarray}
которая определяется как подстратификация стратификации
($\ref{stratQ}$). Определено двулистное накрытие $N_{K_0} \to \hat
N_{K_0}$ и стратификация ($\ref{strat0}$) двулистного
накрывающего.

Естественно высказать следующуюгипотезу: применяя теорему
$PL$--трансверсальности можно показать, что произвольный 2-диск в
полиэдре $\hat N_{K_0}^{0}$ снимается с подполиэдра $ \hat
N_{K_0}^{3} \subset \hat N_{K_0}^{0}$ сколь угодно малым
шевелением, поэтому фундаментальные группы пространств $ \hat
N_{K_0}^{0}$ и $ \hat N_{K_0}^{0} \setminus \hat N_{K_0}^{3}$
совпадают (ср. с Леммой $\ref{RK[2]}$).

Теперь определим гомотопию $\hat f$ \footnote{Конструкция была получена на семинаре
С.А.Богатого}. Рассмотрим нормальное расслоение $\nu_{J_0}: E(\nu) \to J_0$ вложения
$i_{J_0}: J_0 \subset \R^n$. Вложение $i_{J_0}$ определено как стандартное вложение $n-k$--мерной сферы в евклидово пространство $\R^n$ с коразмерностью $k$.
Расслоение  $\nu_{J_0}$ является тривиальным $k$--мерным расслоением, причем по условию леммы $k \ge 8$.
Представим расслоение  $\nu_{J_0}$ в виде суммы Уитни тривиального 4-мерного расслоения $\varepsilon_1:
E(\varepsilon_1) \to J_0$, тривиального 2-мерного расслоения $\varepsilon_2: E(\varepsilon_2) \to J_0$ и
 дополнительного тривиального расслоения $\varepsilon_3: E(\varepsilon_3) \to J_0$
 неотрицательной размерности размерности $k-6$.

Переобозначим отображение $(I_J \circ \hat F)$
через $\hat f_{-2}$. Выберем последовательность бесконечно-малых
чисел $\delta_{-1}, \delta_0, \dots, \delta_{r+1}$ (указанную
последовательность можно представлять как убывающую
последовательность положительных чисел, каждое последующее число
выбирается настолько малым, что указанные в построении свойства
аппроксимирующей гомотопии, достигнутые на предыдущем шаге,
сохраняются)  такую, что $\delta_{i+1} = o(\delta_i)$ и построим
последовательность $C^1$--деформаций
\begin{eqnarray}
\hat f_{-2} \mapsto \hat f_{-1} \mapsto \hat f_{0} \mapsto \dots \mapsto \hat f_r \mapsto \hat f
\label{fi}
\end{eqnarray}
  по индукции, задавая $C^1$--калибр деформации на страте
$C_{p}^{(i)}$ бесконечно-малым числом $\delta_i$. Последняя
гомотопия $\hat f$ в этой последовательности будет
$PL$--отображением общего положения. При достаточно малом $t$, $0
< t < \delta_{r+1}$, гомотопия $\hat f$ определяет отображение
$\hat d$.

На предворительном шаге деформация $\hat f_{-2} \mapsto \hat
f_{-1}$ строится вертикальной вдоль подрасслоения $\varepsilon_1$
в нормальном расслоении $\nu_{J_0}$. Калибр рассматриваемой
деформации равен $\delta_{-1}$. Деформация $\hat f_{-2} \mapsto
\hat f_{-1}$ строится такой, что прообраз
$p^{-1}_{J_0}(\alpha^{(s)})$ каждого элементарного страта
$\alpha^{(s)} \subset J_0$ глубины $s$ самопересекается в
тотальном пространстве $E(\varepsilon_1)_{\alpha^{(s)}}$
расслоения $\varepsilon_1$ над $\alpha^{(s)}$ в общем положении
(как погружения общего положения замкнутого $PL$--многообразия).

На нулевом шаге деформация $\hat f_{-1} \mapsto \hat f_0$ строится вертикальной вдоль слоев расслоения
$\varepsilon_2$, неподвижной на страте  $\hat C_p^{(1)}$. Калибр рассматриваемой деформации равен $\delta_0$.
Деформация $\hat f_{-1} \mapsto \hat f_0$ строится такой, что ограничение отображения $\hat f_{0}$
на прообраз элементарного страта $\alpha^{(0)}$ глубины $0$ является погружением открытого многообразия,
множество точек самопересечения которого совпадает с множеством точек тройного самопересечения погружения
$\hat f_{-1}$, ограниченного на $\alpha^{(0)}$ и не содержит точек двойного самопересечения.

Чтобы определить деформацию $\hat f_{-1} \mapsto \hat f_0$ с
указанными свойствами достаточно воспользоваться Леммой
$\ref{lemma32}$, согласно которой, структурная группа накрытия над
многообразием двойных точек самопересечения на каждом элементарном
страте редуцируется к целочисленной. Это означает, что особенность
двойных точек самопересечения вне можества тройных точек
самопересечения распроектируется во вложение в коразмерности 2
вдоль слоев расслоения $\varepsilon_2$.

На  шаге с номером $i$ деформация $\hat f_{i-1} \mapsto \hat f_i$
строится сначала с носителем на $\hat C_p^{(i)}$.
Предварительная деформация, которую обозначим через $\hat f_{i-1}
\to \hat f'_i$, является неподвижной на $\hat C_p^{(i+1)}$. Далее
построенная гомотопия $\hat f'_i$ продолжается до гомотопии
$\hat f_i$ на всем $\hat C_p$ так, что она
совпадает с уже построеным ранее отображением $\hat f_{i-1}$ вне
малой регулярной окрестности (конечного радиуса) $\hat C_p^{(i)}$
и продолжает отображение  $\hat f_{i-1}$ в эту окрестность на все
пространство $\hat C_p$ по общему положению с константой $C^1$--аппроксимации $\delta_i$.

Построение деформации $\hat
f_{i-1} \mapsto \hat f'_i$, $i \ge 1$, осуществляется следующим образом.
Рассмотрим произвольный элементарный страт
$\alpha^{(i)}$ глубины $i$ из $J^{(i)}_0 \setminus J^{(i+1)}_0$ и рассмотрим
регулярную окрестность $U_{\alpha^{(i)}}$ страта $\alpha^{(i)}$ в объединении стратов
$J_0^{(i-1)} \setminus J^{(i)}_0$ глубины $i-1$, примыкающем к страту $\alpha^{(i)}$.
 Окрестность $U_{\alpha^{(i)}}$ представлена как объединение окрестностей $U_{\beta^{i-1}_j, \alpha^{(i)}}$
по окрестностям страта $\alpha^{(i)}$ во всевозможных примыкающих стратах
$\beta^{(i-1)}_j$ глубины $i-1$. Каждый страт $U_{\beta^{(i-1)}_j, \alpha^{(i)}}$
представляет собой прямое произведение $\alpha^{(i)} \times D^2$.
Прообраз $\hat p^{-1}(U_{\beta^{(i-1)}_j, \alpha^{(i)}})$ окрестности
$U_{\beta^{(i-1)}_j, \alpha^{(i)}}$ при отображении $\hat p:
(S^{n-k}/\i)^{(i-1)} \setminus  (S^{n-k}/\i)^{(i)} \to J_0^{(i-1)}
\setminus J_0^{(i)}$ обозначим через $\hat V_{b_j, a}$.
Объединение $\hat V_{b_j, a}$ по всем $b_j$ обозначим через $\hat
V_{a}$. Это объединение является прообразом $U_{\alpha^{(i)}}$ при
отображении $\hat p$.

Прообраз $\hat V_{b_j, a}$ представляет собой прямое произведение
$4^{r-i-1}$-листного накрывающего $\hat a$ над $\alpha^{(i)}$ на открытый
2-диск $D^2$. Снова воспользуемся координатным описанием
отображения $\hat p$, откуда вытекает, что ограниченное $\hat p$
на $\hat V_{b_j, a}$, является прямым произведением стандартного
разветвленного 4-листного над $D^2$ и параметризующего отображения
$\hat b_j \to \beta_j$.

По предположению индукции гомотопия $\hat f_{i-1}$, ограниченная на $\hat V_{a}$ неподвижна на
$\hat a$, причем  для каждого $\hat b_j$ деформация $\hat f_{i-2} \mapsto \hat f_{i-1}$ направлена в трансверсальном направлении к $U_{\beta_j,\alpha}$
и ее $C^1$--калибр существенно превышает $\delta_i$.

Определим гомотопию $\hat f'_i$
  окрестности $\hat V_{\beta} \to
\R^n$, как результат $\delta_i$--малой $C^1$--деформации на каждой окрестности $\hat V_{a}$
в направлении касательного пространства к произвольно выбранному страту $\beta$ глубины $i-1$, примыкающему к $\alpha$.
Гомотопия $\hat f'_i$ строится также как при $i=0$. В результате точки самопересечения элементарного страта
$\alpha^{(i)}$ совпадают с точками тройного самопересечения этого элементарного страна при гомотопии $f_{-1}$.
При этом ограничение гомотопии $\hat f'_i$ на границу окрестности $\hat V_{\alpha}$, в частности, на границу $\partial (\hat a)^{cl}$ центрального подмногообразия $\hat a
\subset \hat V_{b, a}$ страта $\hat a$ неподвижно.

Определим гомотопию $\hat f_i$
в результате продолжения гомотопии  $\hat f'_i$ по общему положению, близкой к гомотопии $\hat f_{i-1}$  с константой $C^1$--аппроксимации $\delta_i$.
Заметим, что отображение, определенное как результат гомотопии $\hat f_i$ при $t < \delta_{i}$, имеет на каждом элементарном страте окрестности остова глубины не менее $i$ пространства $S^{n-k}/\i^{(1)}$ точки самопересечения, совпадающими с точками тройного самопересечения гомотопии $f_{-1}$ которая была определена на предварительном шаге конструкции. Точки, лежащие на  стратах глубины $i+1$ (и выше)
имеют структуру точек самопересечения общего вида. При этом точки самопересечения стратов разной глубины,
у отображения, определенного  в результате гомотопии  $\hat f_i$, вообще говоря, могут присутствовать.

Деформация  $\hat f_{i-1} \to \hat f_i$ и, тем самым,
последовательность деформаций ($\ref{fi}$), построена. На
последнем шаге гомотопия $\hat f_r$ переводится в гомотопию $\hat
f$ общего положения $C^1$--малой деформацией, причем калибр
$\delta_{r+1}$ этой деформации существенно меньше $\delta_r$. Эта
деформация состоит из двух шагов. На первом шаге результатом
деформации является гомотопия $\hat f_i$ общего положения, для
которой множество точек самопересечения каждого элементарного
страта содержатся в множестве точек самопересечения этого же
элементарного страта при погружении $\hat f_{-1}$. Точки
самопересечения стратов глубины большей $i_{0,max}$ отсутствуют по
соображениям размерности. На втором шаге деформации точки
пересечения стратов разной размерности, которые возникли при
построении вне окрестности точек тройного самопересечения
гомотопии $f_{-1}$, вертикально распроектируются вдоль слоя
расслоения $\varepsilon_3$. Указанное распроектирование вдали от
точек тройного самопересечения корректно определено,
посколькуканоническое двулистное накрытие над полиэдром точек
самопересечения между стратами разной глубины, очевидно,
тривиально (пара точек в накрывающем полиэдре естественно
упорядочена  значением глубины стратов).

Отображение $\hat f$
 определяет
искомое отображение $\hat d$ ($d$) при ограничении на $S^{n-k}/\i
\times \{\delta_r\}$ ($\RP^{n-k} \times \{\delta_r\}$).

По построению  стратификации ($\ref{strat0}$), ($\ref{strathat}$)
допускают семейство отображений $\hat t^{(i)}: W^{(i)}_{\hat K_0}
\to R\hat K_0^{(i)}$, ($ t^{(i)}: W^{(i)}_{K_0} \to RK_0^{(i)}$),
$0 \le i \le i_{max,0}$ в пространства разрешения особенностей.

Более того, построенное семейство отображений $\hat t^{(i)}$ ($ t^{(i)}$)
склеивается в одно отображение $\hat t: \hat N_{K_0}
\to R\hat K_0$ ($t: N_{K_0} \to RK_0$).

 Граничные
условия на $N_{Q,diag}$, $N_{Q,antidiag}$ следуют из граничных
условий на $W_{Q,diag}^{(i)}$, $W_{Q,antidiag}^{(i)}$, последние
предусматриваются в Лемме $\ref{lemma28}$. Лемма $\ref{lemma30}$
доказана.

\subsubsection*{Доказательство Леммы $\ref{lemma31}$}

Доказательство аналогично доказательству Леммы $\ref{lemma30}$.
Лемма $\ref{lemma31}$ доказана.

\section{Многообразия с кватернионным оснащением}

Рассмотрим группу $\Q_a$ кватернионов порядка 8.  Определено
представление $\chi_+: \Q_a \to SO(4)$, которое переводит
единичные кватернионы $\i,\j,\k$ в матрицы $(\ref{Q a1}), (\ref{Q
a2}), (\ref{Q a3})$ соответственно.

Представление $\chi_-$ переводит единичные кватернионы $\i,\j,\k$
в матрицы

\begin{eqnarray}\label{Q-1}
 \left(
\begin{array}{cccc}
0 & 1 & 0 & 0 \\
-1 & 0 & 0 & 0 \\
0 & 0 & 0 & -1 \\
0 & 0 & 1 & 0 \\
\end{array}
\right),
\end{eqnarray}

\begin{eqnarray}\label{Q-2}
\left(
\begin{array}{cccc}
0 & 0 & 1 & 0 \\
0 & 0 & 0 & -1 \\
-1 & 0 & 0 & 0 \\
0 & 1 & 0 & 0 \\
\end{array}
\right),
\end{eqnarray}

\begin{eqnarray}\label{Q-3}
\left(
\begin{array}{cccc}
0 & 0 & 0 & 1 \\
0 & 0 & -1 & 0 \\
0 & 1 & 0 & 0 \\
-1 & 0 & 0 & 0 \\
\end{array}
\right).
\end{eqnarray}

Образ группы кватернионов при представлениях $\chi_+$ и $\chi_-$
определяют также подгруппы в группе $\Z/2^{[3]} \subset O(4)$.
Заметить, что внешнее сопряжение подгруппы $\H \subset \Z/2^{[3]}$
на элемент из $\Z/2^{[3]} \setminus \H$, заданный преобразованием
$\e_1 \to -\e_1$, $\e_i \to \e_i$, $2 \le i \le 4$, определяет
автоморфизм $\H \to \H$, при котором подгруппа $Im \chi_+$
переходит подгруппу $Im \chi_-$ единичных кватернионов в $\H$,
представление которой задается матрицами
($\ref{Q-1}$),($\ref{Q-2}$),($\ref{Q-3}$).

Определена прямая сумма представлений $\chi_+ \oplus \chi_-: \Q_a
\to SO(4) \oplus SO(4) \to SO(8)$, это представление обозначим
через $\psi_{+,-}$. Определено представление $\chi_+ \oplus
\chi_+: \Q_a \to SO(4) \oplus SO(4) \to SO(8)$, которое обозначим
через $\psi_{+,+}$.

Определены свободные действия
$p_7: \Q_a \times S^7 \to S^7$, $p_3: \Q_a \times S^3 \to
S^3$, построенные по представлениям
$\psi_{+,+}$, $\chi_{+}$. Соответствующее факторпространства будем
обозначать через
$S^{7}/\Q_a$, $S^{3}/\Q_a$.

Рассмотрим однородные пространства $S^7/\Q_a$, $S^3/\Q_a$,
определенные как факторпространства действий $p_7$, $p_3$.
Определены пара векторных $8$- и пара векторных $4$-х мерных
расслоения $\zeta_{+,+}: E(\zeta_{+,+}) \to S^7/\Q_a$,
$\zeta_{+,-}: E(\zeta_{+,-}) \to S^7/\Q_a$, $\eta_{+}: E(\eta_{+})
\to S^7/\Q_a$, $\eta_{-}: E(\eta_{-}) \to S^7/\Q_a$ со структурной
группой $\Q_a$. При этом 8-мерное расслоения $\zeta_{+,+}$,
$\zeta_{+,-}$ ассоциированы с представлениями $\psi_{+,+}$,
$\psi_{+,-}$ расслоения $\eta_+$, $\eta_-$ ассоциированы с
представлениями $\chi_+$, $\chi_-$ соответственно.

Рассматривая клеточное разбиение пространства $S^7/\Q_a$,
заключаем, что $H^4(S^7/\Q_a;\Z)=H_3(S^7/\Q_a;\Z)=\Z/8$. При этом
образующая $t$ указанной группы определена как образ
фундаментального класса $[S^3/\Q_a]$ при вложении $S^3/\Q_a
\subset S^7/\Q_a$, которое параметризует 3-остов $K^3 \subset
S^7/\Q_a$.

\begin{proposition}\label{prop32}

Пусть $K^7$ -- произвольное ориентированное замкнутое
многообразие, стабильное нормальное расслоение которого изоморфно
расслоению $f^{\ast}(k\zeta_{+,+})$,
 $f: K^7 \to S^7/\Q_a$, $s = 1 \pmod{2}$. Тогда $\deg(f)=0 \pmod{2}$.
\end{proposition}

\subsubsection*{Доказательство Предложения $\ref{prop32}$}

Заметим, что расслоение над $S^7/\Q_a$, полученное из
$k\zeta_{+,+}$ обращением ориентации, совпадает с $k\zeta_{+,-}$.
Рассмотрим ориентированное (несвязное) многообразие $K^7 \cup
-K^7$, с нормальным расслоением $f^{\ast}k\zeta_{+,+} \cup
f^{\ast}k\zeta_{+,-}$.

 Докажем, что  в группе  $H_3(S^7/\Q_a;\Z)$ c образующей, обозначенной через
 $t$,
справедливы равенства:
\begin{eqnarray}\label{p++}
[p_1(k\zeta_{+,+})]^{op}=4t,
\end{eqnarray}

\begin{eqnarray}\label{p+-}
[p_1(k\zeta_{+,-})]^{op}=0.
\end{eqnarray}

Сначала докажем равенство ($\ref{p++}$) при $k=1$. Расслоение
$\zeta_{+,+}$ является комплексным, следовательно, согласно
Следствию 15.5 из книги [M-S] справедливо равенство:
\begin{eqnarray}\label{p1zeta+}
p_1(\zeta_{+,+})=c^2_1(\zeta_{+,+})-2c_2(\zeta_{+,+}).
\end{eqnarray}
 Поскольку
$\zeta_{+,+}=\eta_+ \oplus \eta_+$, то согласно формуле (14.7)
$c_2(\zeta_{+,+})=c_2(\eta_+ \oplus \eta_+) = c^2_1(\eta_+) +
2c_2(\eta_+)$. Поскольку комплексная размерность расслоения
$\eta_+$ равна $2$, $c_2(\eta_+)$ совпадает с классом Эйлера и
равен $t$. Далее $c_1^2(\eta_+ \oplus \eta_+)=4c_1^2(\eta_+)$.
Суммируя доказанные равенства, c учетом $8t=0$, получим $
p_1(\zeta_{+,+})=2c^2_1(\eta_+) + 4t$.

Осталось доказать, что
\begin{eqnarray}\label{eta+}
c^2_1(\eta_+)=0.
\end{eqnarray}

Справедливо равенство $i^{\ast}(c_1(\eta_+))=0$,
где $i: S^3/\Q_a \subset S^7/\Q_a$ естественное вложение, поскольку расслоение
$i^{\ast}(\eta_+)$ тривиально. Гомоморфизм
$$H^2(S^7/\Q_a;\Z)  \stackrel{i^{\ast}}{\longrightarrow} H^2(S^3/\Q_a;\Z)$$
является мономорфизмом, поэтому $c_1(\eta_+)=0$.  Равенство
($\ref{eta+}$) и равенство ($\ref{p++}$) при $k=1$ доказаны.

Докажем равенство ($\ref{p+-}$) при $k=1$. Расслоение
$\zeta_{+,-}$ является комплексным, следовательно, по аналогичным
вычислениям:
\begin{eqnarray}\label{p1zeta-}
p_1(\zeta_{+,-})=c^2_1(\zeta_{+,-})-2c_2(\zeta_{+,-}).
\end{eqnarray}
Поскольку $\zeta_{+,-}=\eta_+ \oplus \eta_-$, то
$c_2(\eta_+ \oplus \eta_-) = c_1(\eta_+)c_1(\eta_-) + c_2(\eta_+)
+ c_2(\eta_-)$. В последней формуле второе и третье слагаемые
сокращаются, поскольку эйлеровы классы расслоений $\eta_+$,
$\eta_-$ противоположны по знаку. Аналогично предыдущему
вычислению, имеем $c_1^2(\eta_+ \oplus \eta_-)=c_1^2(\eta_+) +
c^2_1(\eta_-) + 2c_1(\eta_+)c_1(\eta_-)$. Суммируя доказанные
равенства, получим $p_1(\zeta_{+,-})=c_1^2(\eta_+) +
c_1^2(\eta_-)$. Остается заметить, что  $c_1^2(\eta_-)=0$ по
аналогичным с ($\ref{eta+}$) вычислениям.
 Равенство ($\ref{p+-}$)
при $k=1$ доказано.

Докажем ($\ref{p++}$), ($\ref{p+-}$) при произвольном нечетном
$k$. Сначала докажем, что у комплексных расслоений
$2l\zeta_{+,+}$, $2l\zeta_{+,-}$ характеристический класс $c_1$
равен нулю, а $c_2(2l\zeta_{+,+})= c_2(2l\zeta_{+,+})=4lt$.
Поскольку $2l\zeta_{+,+}=4l\eta_+$, то согласно формуле (14.7) c
учетом равенства ($\ref{eta+}$), получим:
$c_2(2l\zeta_{+,+})=c_2(4l\eta_+) = 4lc_2(\eta_+)=4lt$. Аналогично
получим $c_2(2l\zeta_{+,-})=4lt$. Равенства ($\ref{p++}$),
($\ref{p+-}$) вытекают из уравнений $(\ref{p1zeta+})$,
$(\ref{p1zeta-})$.

Рассмотрим отображение $F=f \cup f: K^7 \cup -K^7 \to S^7/\Q_a$
ориентированного (несвязного многообразия). По построению
отображение $F$ ориентированно бордантно нулю. С другой стороны,
характеристическое число
$$ F_{\ast}([p_1(k\zeta_{+,+})]^{op} + [p_1(k\zeta_{+,+})]^{op}) \in H_3(S^7/Q_a)$$
сохраняется при кобордизме. Это число равно нулю, при $def(f) = 0
\pmod{2}$, и равно 0,  при $deg(f)=1 \pmod{2}$.
 Предложение $\ref{prop32}$ доказано.
\[   \]

При доказательстве Теоремы 5 было использовано следующее
Предложение.

Обозначим через $\Z/2^{[s]}$   группу, полученную в результате
$s$-кратного сплетения элементарной группы $\Z/2$, которая
определяется аналогично случаю $s=3$. Обозначим через $\omega:
E(\omega) \to K(\Z/2^{[s]},1)$ -- универсальное $2^{s-1}$--мерное
расслоение со структурной группой $\Z/2^{[s]}$.

\begin{proposition}\label{prop33}
Пусть $K$ -- произвольный полиэдр размерности $dim(K)=7$,
снабженный отображением $\varphi: K \to K(\Z/2^{[s]},1)$ и пусть
$\varphi^{\ast}(\omega)$ -- обратный образ расслоения $\omega$ при
отображении $\varphi$. Тогда для любого натурального числа $l$, $l
\equiv 0 \pmod{8}$ расслоение $l \varphi^{\ast}(\omega)$ изоморфно
тривиальному.
\end{proposition}

\subsubsection*{Доказательство Предложения $\ref{prop33}$}
Обозначим расслоение $l \varphi^{\ast}(\omega)$ через $\psi_s$.
Докажем утверждение по индукции. При $s=1$ утверждение хорошо
известно.

Воспользуемся тем, что расслоение $\psi_{s-1}$, $s>1$, тривиально.
Докажем, что расслоение $\psi_s$ тривиально. Рассмотрим двулистное
накрытие $\bar K \to K$, индуцированное подгруппой индекса 2
$\Z/2^{[s-1]} \oplus \Z/2^{[s-1]} \subset \Z/2^{[s]}$. Рассмотрим
расслоение $\psi_s^!$ над $\bar K$, определенное как трансфер
расслоения $\Z/2^{[s]}$. Справедлива формула $\psi_s^! =
\psi_{s-1,1} \oplus \psi_{s-1,2}$, где прямые слагаемые
$\psi_{s-1,1}$, $\psi_{s-1,2}$ суть расслоения над $\bar K$ с
меньшими структурными группами. Эти расслоения изоморфны между
собой, причем изоморфизм осуществляется инволюцией накрытия $\bar
K$. По предположению индукции расслоения $\psi_{s-1,1}$,
$\psi_{s-1,2}$ тривиальны. Следовательно, для расслоений над $K$
справедливо равенство $\psi_s = l2^{s-2}\kappa \oplus l2^{s-2}
\varepsilon$, где $\kappa$ -- линейное расслоение, класифицирующее
накрытие $\bar K \to K$, $\varepsilon$ -- тривиальное расслоение.
Расслоение $l\kappa$ и, тем более, расслоение $l2^{s-2}\kappa$,
тривиально, следовательно, расслоение $\psi_s$ тривиально.
Предложение $\ref{prop33}$ доказано.

 \[  \]
 Московская обл. г.Троицк 142190
ИЗМИРАН

 pmakhmet@izmiran.ru
\[  \]

\end{document}